\documentclass{amsart}
\usepackage{amsmath, mathrsfs,mathtools,amssymb,mathabx,accents,braket,dsfont}
\usepackage{xcolor,tikz}
\usepackage{enumerate}
\usepackage[capitalise]{cleveref}

\newcommand{\rr}{\mathbf{r}}
\newcommand{\sbold}{\mathbf{s}}
\newcommand{\tbold}{\mathbf{t}}
\newcommand{\uu}{\mathbf{u}}
\newcommand{\vv}{\mathbf{v}}
\newcommand{\ww}{\mathbf{w}}

\DeclareMathOperator{\Cliq}{Cliq}

\DeclareMathOperator{\Id}{Id}

\DeclareMathOperator{\Comm}{Comm}

\DeclareMathOperator{\Diag}{Diag}

\DeclareMathOperator{\Span}{Span}

\DeclareMathOperator{\Cstar}{C^*}

\DeclareMathOperator{\Sp}{Sp}

\newcommand{\CC}{\mathbb{C}}

\newcommand{\NN}{\mathbb{N}}

\newcommand{\ZZ}{\mathbb{Z}}

\newcommand{\calA}{\mathcal{A}}
\newcommand{\calB}{\mathcal{B}}

\newcommand{\calF}{\mathcal{F}}

\newcommand{\calH}{\mathcal{H}}
\newcommand{\calI}{\mathcal{I}}
\newcommand{\calJ}{\mathcal{J}}
\newcommand{\calK}{\mathcal{K}}
\newcommand{\calL}{\mathcal{L}}
\newcommand{\calM}{\mathcal{M}}
\newcommand{\calN}{\mathcal{N}}

\newcommand{\calP}{\mathcal{P}}
\newcommand{\calQ}{\mathcal{Q}}

\newcommand{\calS}{\mathcal{S}}
\newcommand{\calT}{\mathcal{T}}

\newcommand{\calW}{\mathcal{W}}

\newcommand{\boldA}{\mathbf{A}}
\newcommand{\boldB}{\mathbf{B}}

\newcommand{\boldM}{\mathbf{M}}
\newcommand{\boldN}{\mathbf{N}}

\newtheorem{theorem}{Theorem}[section]
\newtheorem{thmx}{Theorem}

\newtheorem{lemma}[theorem]{Lemma}
\newtheorem{corollary}[theorem]{Corollary}
\newtheorem{example}[theorem]{Example}

\newtheorem{definition}[theorem]{Definition}
\newtheorem{proposition}[theorem]{Proposition}

\theoremstyle{remark}

\newtheorem{remark}[theorem]{Remark}

\usepackage[date=year,backend =bibtex,doi=false,isbn=false,url=false,eprint=true]{biblatex}
\title{The CCAP for graph products of operator algebras}
\date{\today} 
\author[Matthijs Borst]{Matthijs Borst}
\address{TU Delft, EWI/DIAM,
	P.O.Box 5031,
	2600 GA Delft,
	The Netherlands}
\email{m.j.borst@outlook.com}
\begin{document}

	\begin{abstract}
		For a simple graph $\Gamma$ and for unital $\Cstar$-algebras with GNS-faithful states $(\boldA_v,\varphi_v)$ for $v\in V\Gamma$,  we consider the reduced graph product  $(\calA,\varphi)=*_{v,\Gamma}(\boldA_{v},\varphi_v)$ , and show that if every $\Cstar$-algebra $\boldA_{v}$ has the completely contractive approximation property (CCAP) and satisfies some additional condition, then the graph product has the CCAP as well. The additional condition imposed is satisfied in natural cases, for example for the reduced group $\Cstar$-algebra of a discrete group $G$ that possesses the CCAP. 
		
		Our result is an extension of the result of Ricard and Xu in \cite[Proposition 4.11]{ricardKhintchineTypeInequalities2006a} where they prove this result under the same conditions for free products.  Moreover, our result also extends the result of Reckwerdt in \cite[Theorem 5.5]{reckwerdtWeakAmenabilityStable2017}, where he proved for groups that weak amenability with Cowling-Haagerup constant $1$ is preserved under graph products. Our result further covers many new cases coming from Hecke-algebras and discrete quantum groups.
	\end{abstract}
	\maketitle
	\section{Introduction}
	In this paper we look at graph products of operator algebras. These graph products are a generalization of free products, where certain commutation relations are added. 
	The notion of graph products was first introduced for groups, by Green in her thesis \cite{greenGraphProductsGroups1990}. For groups $G_i$ the free product $G=*_{i} G_i$ consists of all reduced words $g_1\cdots g_l$ with $g_j\in G_{i_j}$ and the group operation consists of concatenation, and reduction. For a given simple graph $\Gamma$, and groups $G_v$ for every vertex $v$, the graph product $G_{\Gamma} = *_{v, \Gamma} G_v$ is obtained from the free product by declaring elements $g_1\cdots g_kg_{k+1}\cdots g_l$ and $g_1\cdots g_{k-1}g_{k+1}g_kg_{k+2}\cdots g_l$ to be equal whenever $g_k\in G_{v}$ and $g_{k+1}\in G_{w}$ and $v$ and $w$ share an edge in $\Gamma$. Graph products preserve many interesting properties like: soficity \cite{ciobanuSoficGroupsGraph2014}, residual finiteness \cite{greenGraphProductsGroups1990}, rapid decay \cite{ciobanuRapidDecayPreserved2011} and other properties, see \cite{antolinTitsAlternativesGraph2015,chiswellOrderingGraphProducts2012, hermillerAlgorithmsGeometryGraph1995, hsuLinearResidualProperties1999}. In particular, approximation properties like the Haagerup property \cite{antolinHaagerupPropertyStable2014} and weak-amenability with constant $1$ \cite{reckwerdtWeakAmenabilityStable2017} are also preserved by graph products of groups.\\
	
	Graph products of operator algebras were introduced in \cite{caspersGraphProductsOperator2017a} by Caspers and Fima as a generalization of free products. Their notion of graph products agrees with that for groups in the sense that for discrete groups $G_v$ one has $*_{v,\Gamma}C_r^*(G_v) = C_r^*(*_{v,\Gamma}{G_v})$ and $\overline{*_{v,\Gamma}}\calL(G_v) = \calL(*_{v,\Gamma}{G_v})$. In their paper, they also showed stability of exactness (for $\Cstar$-algebras), Haagerup property, II$_1$-factoriality (for von Neumann algebras) and rapid decay (for certain discrete quantum groups) under graph products. Also, in \cite{caspersConnesEmbeddabilityGraph2016} it was proven that embeddability is preserved under graph products.\\
	
	The notion of weak amenability for groups originates from the work of Haagerup \cite{haagerupExampleNonNuclear1978}, De Canni\`ere-Haagerup \cite{decanniereMultipliersFourierAlgebras1985a} and Cowling-Haagerup \cite{cowlingCompletelyBoundedMultipliers1989}.   The corresponding notion for unital $\Cstar$-algebras is given by the completely bounded  approximation property (CBAP) in the sense that a discrete group is weakly amenable if and only if its reduced group $\Cstar$-algebra possesses the CBAP. We say that a $\Cstar$-algebra $\boldA$ has the CBAP if there exists a net of completely bounded maps $V_n: \boldA \rightarrow \boldA$ that are finite rank, converge to the identity in the point-norm topology and such that $\sup_n \| V_n \|_{cb} \leq \Lambda < \infty$ for some constant $\Lambda$. The minimal such $\Lambda$ is called the Cowling-Haagerup constant.  If the Cowling-Haagerup constant is 1, then we say that $\boldA$ has the completely contractive approximation property (CCAP).  
	
	Weak amenability and the CBAP/CCAP play a crucial role in functional analysis and operator algebras. Already in case of the group $G =\ZZ$ weak amenability allows, in a way, to approximate a Fourier series by its partial sums. In operator space theory the CBAP has led to a deep understanding of several group $\Cstar$- and von Neumann algebras. Already the results by Cowling and Haagerup \cite{cowlingCompletelyBoundedMultipliers1989} allow for the distinction of group von Neumann algebras of lattices in the Lie groups $\Sp(1,n), n \geq 2$. Later, Ozawa and Popa used the (wk-$*$) CCAP in deformation-rigidity theory of von Neumann algebras \cite{ozawaClassIIFactors2010}. Much more recently also graph products have appeared in the deformation-rigidity programme, see e.g.
	\cite{caspersAbsenceCartanSubalgebras2020}, \cite{borstBimoduleCoefficientsRiesz2021}  \cite{chifanCartanSubalgebrasNeumann2021}, \cite{chifanRigidityNeumannAlgebras2022},\cite{dingProperProximalityVarious2022}. This line of investigation, especially beyond the realm of group algebras, motivates the study of the CCAP for general graph products.

	In this paper we are concerned with showing that the CCAP is preserved under graph products. While we are not able to show this in full, we prove this under a mild extra condition on the algebras $(\boldA_{v},\varphi_v)$, similar to the one imposed by \cite{ricardKhintchineTypeInequalities2006a} for proving the same result for free products.  The conditions that we impose are stated in \cref{section:ucp-extension-for-ccap-preserved}, and we abbreviate them by saying that the algebra \textit{has a u.c.p extension for the CCAP}. This condition is satisfied by many natural unital $\Cstar$-algebras, under which finite-dimensional ones (with a GNS-faithful state), reduced $\Cstar$-algebras of discrete groups (with the Plancherel state) that possess the CCAP \cite{ricardKhintchineTypeInequalities2006a}, and reduced $\Cstar$-algebras of compact quantum groups (with the Haar state) whose discrete dual quantum group is weakly amenable with Cowling-Haagerup constant $1$ \cite{freslonNoteWeakAmenability2012}. Our main result is the following:	
	\begin{thmx}\label{thm:main-result}
		Let $\Gamma$ be a simple graph and for $v\in V\Gamma$ let $(\boldA_v,\varphi_v)$ be  unital $\Cstar$-algebras that have a u.c.p. extension for the CCAP. Then the reduced graph product $(\calA,\varphi) = *_{v,\Gamma}(\boldA_{v},\varphi_{v})$ has the CCAP.
	\end{thmx}

Along the way, in \cref{corollary:CCAP-for-finite-dim-algebras}, we also obtain the following result for von Neumann algebras.
\begin{thmx}\label{thm:intro-finite-dim-vNa}
		Let $\Gamma$ be a simple graph and for $v\in V\Gamma$ let $\boldM_{v}$ be a finite-dimensional von Neumann-algebra together with a normal faithful state $\varphi_v$. Then the von Neumann algebraic graph product $(\calM,\varphi)=\overline{*_{v,\Gamma}}(\boldM_v,\varphi_v)$ has the wk-$*$ CCAP.
\end{thmx}
	The method for proving above results is, on a large scale, similar to \cite{ricardKhintchineTypeInequalities2006a}. However, at most points, the proofs get more involved in order to work for graph products. This becomes most clear in  \cref{section:polynomial-growth-word-length}, where we have to use different methods to show the completely boundedness of the word-length projection maps $\calP_d$ that project on  $\calA_d$, the homogeneous subspace of order $d$. For these maps we show for $d\geq 1$ the linear bound $\|\calP_d\|_{cb}\leq C_{\Gamma}d$, where $C_{\Gamma}$ is some constant only depending on the graph $\Gamma$. In \cref{section:graph-products-of-ucp-maps} we show that the graph product map $\theta$ of state-preserving u.c.p. maps $\theta_v$ on unital $\Cstar$-algebras $\boldA_{v}$, is again a state-preserving u.c.p. map on the reduced graph product $\calA$. Together with our bound on $\|\calP_d\|_{cb}$ we are then able to show the preliminary result, \cref{corollary:CCAP-for-finite-dim-algebras}, that, when all $\Cstar$-algebras,  respectively von Neumann algebras, are finite-dimensional, the reduced graph product has the CCAP, respectively the wk-$*$ CCAP. In \cref{section:graph-product-of-cb-maps-on-Ad} we consider the same problem as in \cref{section:graph-products-of-ucp-maps}, but now for state-preserving completely bounded  maps. We show that the graph product map $T$ of state-preserving completely bounded maps $T_v$ defines a completely bounded map, when restricted to a homogeneous subspace $\calA_d$ (i.e. $T_d := T|_{\calA_d}$ is completely bounded). In order to do this we consider the operator spaces $X_d$ from \cite{caspersGraphProductKhintchine2021a} (analogous to \cite{ricardKhintchineTypeInequalities2006a}) and use the Khintchine type inequality \cite[Theorem 2.9]{caspersGraphProductKhintchine2021a} they proved. We moreover construct other operator spaces $\widetilde{X}_d$ and prove the `reversed'  Khintchine type inequality (\cref{thm:theta-d-completely-contractive}). Finally, in \cref{section:ucp-extension-for-ccap-preserved}, using all our previous results, we are then able to show the main result \cref{thm:main-result} (\cref{thm:ucp-extension-gives-CCAP-for-graph-products}).\\
	
	Our results extends \cite{ricardKhintchineTypeInequalities2006a} (as well as \cite{reckwerdtWeakAmenabilityStable2017}) in a natural way, and provides a unified approach to proving the CCAP and wk-$*$ CCAP for various operator algebras. Specifically, 
	\cref{thm:intro-finite-dim-vNa} can be applied to the graph product $*_{v,\Gamma}\calN_{q_v}(W_v)$ of Hecke-algebras of finite Coxeter groups. Such a graph product is itself a Hecke-algebra, and by the result we obtained, possesses the wk-$*$ CCAP. This result is new, and was previously only known, by \cite[Theorem A]{caspersAbsenceCartanSubalgebras2020}, for the case that $W_v$ is right-angled for all $v$. Furthermore, the main theorem, \cref{thm:main-result}, can be applied to give new examples of $\Cstar$-algebras that posses the CCAP, for example the graph product $*_{v,\Gamma}(\boldA_v,\varphi_v)$, where some algebras $\boldA_{v}$ are finite-dimensional, and others are reduced group $\Cstar$-algebras of discrete groups that posses the CCAP. 
	
	\section{Preliminaries}	
	We will use basic notions from  $\Cstar$-algebras and von Neumann algebras, for which we refer to \cite{murphyAlgebrasOperatorTheory1990a}. Also, in \cref{section:graph-product-of-cb-maps-on-Ad}, we will use some theory from operator spaces for which we refer to \cite{effrosOperatorSpaces2000},\cite{pisierIntroductionOperatorSpace2003}. Here, in this section, we shall recall the definitions of weak amenability and of graph products of operator algebras, and establish the notation that we shall use for this throughout the paper. We also state \cref{lemma:partition-action} (see  \cite[Lemma 2.5.]{caspersGraphProductKhintchine2021a}) and prove  \cref{lemma:action-on-word-part} that we shall need later for calculations.

	\subsection*{Weak amenability with Cowling-Haagerup constant $1$}
	We recall the definition of the CCAP for unital $\Cstar$-algebras and the wk-$*$ CCAP for von Neumann algebras. A unital $\Cstar$-algebra $\boldA$ with state $\varphi$ is said to have the CCAP if there exists a net $(V_j)_{j\in J}$ of completely contractive, finite-rank maps on $\boldA$ s.t. $V_j\to \Id$ pointwise in the norm-topology. A von Neumann algebra $\boldM$ is said to have the wk-$*$ CCAP if there exists a net $(V_j)_{j\in J}$ of normal, completely contractive, finite-rank maps on $\boldM$ s.t. $V_j\to \Id$ pointwise in the $\sigma$-weak topology.
	
	\subsection*{Graph products of operator algebras}
	Let $\Gamma$ be a finite graph that is simple (i.e. undirected, no multiple edges, no edges that start and end in the same vertex), with to each vertex $v\in V\Gamma$ associated a unital $\Cstar$-algebra $\boldA_{v}$ together with a state $\varphi_{v}$ on $\boldA_{v}$. Moreover, for $v\in V\Gamma$, let $\pi_v:\boldA_{v}\to \calB(\calH_v)$ be a given faithful representation of $\boldA_v$ on a Hilbert space $\calH_v$ such that for $a\in \boldA_{v}$ we have $\varphi_v(a) = \langle \pi_v(a)\xi_v,\xi_v\rangle$ for some unit vector $\xi_v\in \calH_v$. In the case that the states $\varphi_v$ are GNS-faithful (meaning the GNS-representations are faithful), and we do not specify the representations, we take the GNS-representation $(\calH_{v},\pi_{v},\xi_{v})$ for $\boldA_{v}$ and simply consider $\boldA_v\subseteq \calB(\calH_v)$ as a subalgebra.
	We will moreover denote $\mathring{\boldA}_v = \ker \varphi_{v}$ and  $\mathring{\calH}_{v} := \xi_{v}^{\perp}$. Moreover, for an element $a\in \boldA_v$, we will write $\mathring{a}:= a- \varphi_v(a)\in \mathring{\boldA}_{v}$ and $\widehat{a} :=\pi_v(a)\xi_{v}\in \calH_{v}$. For a vector $\eta\in \calH_{v}$ we will denote $\mathring{\eta} = \eta - \langle \eta,\xi_{v}\rangle \xi_{v}\in \mathring{\calH}_{v}$. We note that $a\in \mathring{\boldA}_{v}$ implies $\widehat{a}\in \mathring{\calH}_v$.
	
	\subsubsection{The Coxeter group}
	We will call a finite sequence $(v_1,\ldots,v_n)$ of elements of $V\Gamma$ a \textit{word}, and we will denote the set of all words by $\calW$. 
	This includes the empty word.
	We equip the set $\calW$ with the equivalent relation generated by
	\begin{align}
	(v_1,\ldots,v_n) &\sim (v_1,\ldots,v_{i-1},v_{i+2},\ldots,v_n) \text{ whenever } v_{i}=v_{i+1}\\
	(v_1,\ldots,v_n) &\sim (v_1,\ldots,v_{i-1},v_{i+1},v_i,v_{i+2},\ldots,v_n) \text{ whenever } (v_{i},v_{i+1})\in E\Gamma.
	\end{align} 
	We will call a word $(v_1,\ldots,v_n)$ \textit{reduced} if it is not equivalent to a word $(v_1',\ldots,v_m')$ with $m<n$. We note that if two reduced words are equivalent, then they must have equal length. Also we note that every word is equivalent to a reduced word.
	We shall now denote $W$ for the set of words $\calW$ modulo the equivalence relation. We equip $W$ with the operation of concatenating tuples, which makes $W$ into a group. We denote $e$ for the identity element in $W$, which is the equivalence class corresponding to the empty word in $\calW$. We note that, in fact, $W$ equals the right-angled Coxeter group whose Coxeter diagram is the graph $\Gamma$ (for references on Coxeter groups, see \cite[Chapter 3]{davisGeometryTopologyCoxeter2008}). For a word $(v_1,\ldots,v_n)\in \calW$ we will write $v_1\cdots v_n$ for the corresponding element in $W$. For every Coxeter element $\ww\in W$, we will fix a reduced element $(w_1,\ldots,w_n)$ in the equivalence class $\ww$. This element will be called the \textit{representative of $\ww$}. Furthermore, we will write $|\ww|$ for the \textit{length} of $\ww$, which we define as the length of its representative. We remark here that $|e|=0$.
	If $\ww_1,\ldots,\ww_n\in W$, we will say that the expression $\ww_1\cdots \ww_n$ is \textit{reduced} if it holds that $|\ww_1|+\ldots+|\ww_n| = |\ww_1\cdots \ww_n|$.  We will say that a word $\ww\in W$ \textit{starts with a word $\uu\in W$} when $|\ww|=|\uu|+|\uu^{-1}\ww|$, and similarly we will say that $\ww$ \textit{ends with a word $\uu\in W$} when $|\ww| = |\ww\uu^{-1}|+|\uu|$. A word $\ww\in W$ with representative $(w_1,\ldots ,w_n)$ will be called a \textit{clique word} when any two letters $w_i$ and $w_j$ with $i\not=j$ share an edge in $\Gamma$.
	For a word $\ww\in W$ we define $\sbold_l(\ww)$ and $\sbold_r(\ww)$ as the maximal clique words that $\ww$ respectively starts with and ends with. We note that $\sbold_l(\ww) = \sbold_r(\ww^{-1})$.

	\subsubsection{The Hilbert spaces}
	For a word $\ww\in W$, $\ww\not=e$ with representative $(w_1,\ldots,w_n)\in \calW$ define the Hilbert spaces
	\begin{align}
	\mathring{\calH}_{\ww} &= \mathring{\calH}_{w_1}\otimes \cdots \otimes 
	\mathring{\calH}_{w_n}
	\end{align}
	We also set
	\begin{align}
	\mathring{\calH}_{e} =  \CC\Omega
	\end{align}
	where the vector $\Omega$ is called the \textit{vacuum} vector.
		For $d\geq 0$ set
	\begin{align}
	\calF_d &= \bigoplus_{\ww\in W, |\ww| = d}\mathring{\calH}_{\ww}
	\end{align}
	and set
	\begin{align}
	\calF = \bigoplus_{\ww\in W}\mathring{\calH}_{\ww}.
	\end{align}	

	\subsubsection{The operator algebras}
	For an element $\ww\in W$, $\ww\not=e$ with representative $(w_1\ldots w_l)\in \calW$ define the algebraic tensor products
	\begin{align}
	\mathring{\boldA}_{\ww} &= \mathring{\boldA}_{w_1}\otimes \cdots \otimes \mathring{\boldA}_{w_l}.
	\end{align}
	 Also define
	\begin{align}
	\mathring{\boldA}_{e} = \calB(\mathring{\calH}_{e}).
	\end{align}
	
	Moreover, for $d\geq 0$ we define the direct sums
	\begin{align}
	\mathring{\boldA}_d &= \bigoplus_{\substack{\ww\in W\\ |\ww|=d}}\mathring{\boldA}_{\ww}
	\end{align}
	Now we set
	\begin{align}
	\boldA &= \bigoplus_{\ww\in W}\mathring{\boldA}_{\ww}
	\end{align}
	
	\subsubsection{Identifying Hilbert spaces and operator algebras}
	Let $(\vv_1,\ldots,\vv_n)\in W^n$ be s.t. $|\vv_1\cdots\vv_n|=|\vv_1|+\ldots+|\vv_n|$. Write $\calI$ for the set of all indices $1\leq i\leq n$ s.t. $\vv_i\not=e$.
	For $i\in \calI$ write $(v_{(i,1)},\ldots,v_{(i,l_i)})\in \calW$ for the representative of $\vv_i$. Also, write $(\widetilde{v}_1,\ldots,\widetilde{v}_l)\in \calW$ for the representative of $\vv:=\vv_1\cdots\vv_n$. By the assumption it holds that $l = \sum_{i\in \calI}l_i $.
	For convenience, we define a bijection $\sigma$ from $\{1,\ldots,l\}$ to $\{(i,j)| i\in \calI, 1\leq j\leq l_i\}$ as 
	$\sigma(m)=(i,j)$ where  $(i,j)$ is uniquely chosen with the property that $m = j+\sum_{k\in I, k<i}l_k$.
	Now, we have by the definitions that $(v_{\sigma(1)},\ldots,v_{\sigma(l)}) \sim (\widetilde{v}_1,\ldots,\widetilde{v}_l)$. 
	Therefore, by \cite[Lemma 2.3]{caspersGraphProductsOperator2017a} we obtain that there is a unique permutation $\pi$ of $\{1,\ldots,l\}$ with the property that
	\begin{align}
	(v_{\sigma(\pi(1))},\ldots,v_{\sigma(\pi(l))}) = (\widetilde{v}_{1},\ldots,\widetilde{v}_{l})
	\end{align}
	and satisfying that if $1\leq i<j\leq l$ are s.t. $v_{\sigma(i)}=v_{\sigma(j)}$, then  $\pi(i)<\pi(j)$.

	We will now define a unitary 	$\calQ_{(\vv_1,\ldots,\vv_n)}:\mathring{\calH}_{\vv_1}\otimes \cdots\otimes \mathring{\calH}_{\vv_n}\to \mathring{\calH}_{\vv_1\cdots\vv_n}$ as follows.
	For $i\in \calI$ choose pure tensors $\eta_i = \eta_{i,1}\otimes \cdots\otimes \eta_{i,l_i}\in \mathring{\calH}_{\vv_i}$ and for $1\leq i\leq n$ with $i\not\in \calI$ denote $\eta_i = \Omega$. We define
	\begin{align}
	\calQ_{(\vv_1,\ldots,\vv_n)}(\eta_1\otimes \cdots\otimes \eta_n) = \begin{cases}
		\eta_{\sigma(\pi(1))}\otimes \cdots\otimes \eta_{\sigma(\pi(l))} & \text{ when } \calI\not=\emptyset\\
		\Omega & \text{ when }\calI=\emptyset
	\end{cases}
	\end{align} and we extend this definition linearly to a bounded map.
	
	Similarly, we define another map 	$\calQ_{(\vv_1,\ldots,\vv_n)}:\mathring{\boldA}_{\vv_1}\otimes\cdots\otimes \mathring{\boldA}_{\vv_n}\to \mathring{\boldA}_{\vv_1\cdots\vv_n}$, denoted by the same symbol, as follows.
	For $i\in \calI$ choose pure tensors $a_i = a_{i,1}\otimes \cdots\otimes a_{i,l_i}\in \mathring{\boldA}_{\vv_i}$ and for $1\leq i\leq n$ with $i\not\in \calI$ denote $a_i = \Id_{\mathring{\calH}_e}$. We define
	\begin{align}
	\calQ_{(\vv_1,\ldots,\vv_n)}(a_1\otimes \cdots\otimes a_n) =
	\begin{cases}
		a_{\sigma(\pi(1))}\otimes \cdots\otimes a_{\sigma(\pi(l))} & \text{ when }\calI\not=\emptyset\\
		\Id_{\calF}	& \text{ when }\calI=\emptyset\\
	\end{cases} 
	\end{align} and we extend this definition to a linear map.
	
	\subsubsection{Defining the graph product}
	For $\uu\in W$ let $\calW^L(\uu)$ be the set of words $\ww\in W$, s.t. $\uu\ww$ is reduced.
	We define
	\begin{align}
	\calH^L(\uu) = \bigoplus_{\ww\in \calW^L(\uu)}\mathring{\calH}_{\ww}.
	\end{align}
	We will now, for $u\in V\Gamma$, define a unitary $U_u: \calH_{u}\otimes \calH^L(u)\to \calF$ by setting
	\begin{align}
	U_{u}|_{\mathring{\calH}_{u}\otimes \mathring{\calH}_{\ww}} &= \calQ_{(u,\ww)} &\quad \text{for } \ww\in \calW^L(u)\\
	U_{u}(\xi_{u}\otimes \eta)&=\eta &\quad \text{ for } \eta\in \calH^L(u).
	\end{align}
	Furthermore, we define for $u\in V\Gamma$ an operator $\lambda_u: \calB(\calH_u)\to \calB(\calF)$ as 
	\begin{align}
	\lambda_u(a) = U_u(a\otimes \Id)U_u^*.
	\end{align}
	
	The definitions of $U_u$ and $\lambda_u(a)$ are the same as in \cite{caspersGraphProductsOperator2017a} and the intuition behind these maps is as follows.
	The unitary $U_u^*$ represents a pure tensor $\eta=\eta_{v_1}\otimes \cdots \otimes \eta_{v_n} \in \mathring{\calH}_{\vv}\subseteq \calF$ by an element in $\calH_{u}\otimes \calH^L(u)$ by either shuffling the indices (when $\vv$ starts with $u$), or tensoring with the vector $\xi_u$ (when $\vv$ does not start with $u$). The operator $\lambda_u(a)$ acts on $\eta\in\calF$ by rearranging the tensor $\eta$ using $U_u^*$, acting with $a$ on the part in $\calH_{u}$, and subsequently using $U_u$ to map the vector back to an element from $\calF$.
	
	This construction also coincides with \cite[Section 1.5]{caspersGraphProductKhintchine2021a} where the shuffling is done implicit by using an equivalence relation (called shuffle equivalence) to identify Hilbert spaces $\mathring{\calH}_{w_1}\otimes \cdots \otimes \mathring{\calH}_{w_n}$ and $\mathring{\calH}_{w_1'}\otimes \cdots \otimes \mathring{\calH}_{w_n'}$ whenever $w_1\cdots w_n=w_1'\cdots w_n'$ are two reduced expressions for the same word. The action is then defined by $a\cdot \eta = \widehat{\mathring{a}}\otimes \eta + \varphi(a)\eta$ when $\vv$ does not start with $u$, and $a\cdot \eta = \mathring{(a\eta_0)}\otimes \eta' + \langle a\eta_0,\xi_u\rangle \eta'$ when $\vv$ starts with $u$ and $\eta$ is shuffle equivalent to $\eta_0\otimes \eta'\in \mathring{\calH}_{u}\otimes \mathring{\calH}_{u\vv}$.\\

	We will define a linear map $\lambda: \boldA\to \calB(\calF)$ for $\ww\in W$ with representative $(w_1,\ldots w_t)\in \calW$ and for a pure tensor $a=a_1\otimes \cdots\otimes a_t\in \mathring{\boldA}_{\ww}$  as 
	\begin{align}
	\lambda(a_1\otimes \cdots\otimes a_t) &= \lambda_{w_1}(a_1)\lambda_{w_2}(a_2)\ldots\lambda_{w_t}(a_t)
	\end{align}
	and we moreover define $\lambda(\Id_{\mathring{\calH}_e})=\Id_{\calF}$.
	We note that $\lambda$ is injective as $\widehat{a}:=\lambda(a)\Omega = \hat{a}_1\otimes \cdots \otimes \hat{a}_n$ for $a = a_1\otimes \cdots \otimes a_n \in \mathring{\boldA}_{\ww}$.	
	We moreover note that for words $\vv_1,\ldots,\vv_n\in W$ with $|\vv_1|+\ldots+|\vv_n| = |\vv_1\cdots\vv_n|$ and elements $a_i\in \mathring{\boldA}_{\vv_i}$ we have
	for $a = \calQ_{(\vv_1,\ldots,\vv_n)}(a_1\otimes \cdots\otimes a_n)$ that $\lambda(a) = \lambda(a_1)\ldots\lambda(a_n)$. \\
	
	Now, we define the \textit{graph product of  unital $\Cstar$-algebras} as
	\begin{align}
	\calA &:= *_{v,\Gamma}(\boldA_{v},\varphi_v) :=  \overline{\lambda(\boldA)}^{\|\cdot\|}
\end{align}
	Also, for $d\geq 0$ we define the \textit{homogeneous subspace of degree $d$}  as
	\begin{align}
	\calA_{d}&:=\overline{\lambda(\boldA_d)}^{\|\cdot\|}.
	\end{align}
	We moreover define a state $\varphi$ on $\calA$ as $\varphi(a)=\langle a\Omega,\Omega\rangle$, so that $\varphi(\Id_{\calF})=1$ and $\varphi(a)=0$ for $a\in \lambda(\mathring{\boldA}_{\ww})$, with $\ww\not=e$.
	We note that for $v\in V\Gamma$ we have that $\boldA_{v}$ is isomorphic to $\lambda(\boldA_{v})\subseteq \calA$, and that  $\varphi_v = \varphi \circ \lambda|_{\boldA_{v}}$. When we are using the GNS-representations, we will call $\calA$ the \textit{reduced} graph products. 

	Similarly, when all $\boldA_{v}$ for $v\in V\Gamma$ are von Neumann algebras, and the states $\varphi_{v}$ are all normal, we define the  \textit{graph product of von Neumann algebras} as
	 \begin{align}
	 	\calM &:= \overline{*_{v,\Gamma}}(\boldA_{v},\varphi_v) :=  \overline{\lambda(\boldA)}^{SOT}
	 \end{align}
 	and the state $\varphi$ is normal on $\calM$ in that case. 	We also define the \textit{homogeneous subspace of degree $d$} as 
	 \begin{align}
	 	\calM_{d}&:=\overline{\lambda(\boldA_d)}^{SOT}
	 \end{align}

	\subsubsection{Creation, annihilation and diagonal operators}
		For $v\in V\Gamma$ denote $P_v\in \calB(\calF)$ for the projection on the complement of $\calH^L(v)$.  Let $\ww\in W$, $\ww\not=e$ and let $a=a_1\otimes \cdots \otimes a_n\in \mathring{\boldA}_{\ww}=\mathring{\boldA}_{w_1}\otimes \cdots \otimes \mathring{\boldA}_{w_n}$. We now define the
	\textit{annihilation operator} $\lambda_{ann}:\boldA\to \calB(\calF)$, the \textit{diagonal operator} $\lambda_{dia}:\boldA\to \calB(\calF)$ and the \textit{creation operator} $\lambda_{cre}:\boldA\to \calB(\calF)$ by
	\begin{align}
	\lambda_{ann}(a_1\otimes \cdots\otimes a_n)&=(P_{w_1}^\perp \lambda(a_1)P_{w_1})(P_{w_2}^\perp\lambda(a_2)P_{w_2})\ldots (P_{w_{n}}^{\perp}\lambda(a_{n})P_{w_n})\\
	\lambda_{dia}(a_1\otimes \cdots \otimes a_n)&=(P_{w_1} \lambda(a_1)P_{w_1})(P_{w_2} \lambda(a_2)P_{w_2})\ldots(P_{w_n} \lambda(a_n)P_{w_n})\\
	\lambda_{cre}(a_1\otimes \cdots \otimes a_n)&=(P_{w_1} \lambda(a_1)P_{w_1}^{\perp})(P_{w_2} \lambda(a_2)P_{w_2}^{\perp})\ldots(P_{w_n} \lambda(a_n)P_{w_n}^{\perp})
	\end{align} and by $\lambda_{ann}(\Id_{\mathring{\calH}_e}) = \lambda_{dia}(\Id_{\mathring{\calH}_e}) = \lambda_{cre}(\Id_{\mathring{\calH}_e})=\Id_{\calF}$  and extended linearly.
	For $\eta\in \mathring{\calH}_{\vv}$ for some $\vv\in W$ and $b\in \mathring{\boldA}_{w}$ we see that $\lambda_{ann}(b)\eta\in \mathring{\calH}_{w\vv}$ when $\vv$ starts with $w$ and that $\lambda_{ann}(b)\eta =0$  when $\vv$ does not start with $w$. Also, we see that $\lambda_{dia}(b)\eta \in \mathring{\calH}_{\vv}$ when $\vv$ starts with $w$ and that $\lambda_{dia}(b)\eta =0$ when $\vv$ does not start with $w$. Similarly, we see that $\lambda_{cre}(b)\eta \in \mathring{\calH}_{w\vv}$ when $\vv$ does not start with $w$ and that $\lambda_{cre}(b)\eta=0$ when $\vv$ starts with $w$.
	Now, using the fact that $\lambda_{ann}(a) = \lambda_{ann}(a_1)\ldots \lambda_{ann}(a_n)$, $\lambda_{dia}(a) = \lambda_{dia}(a_1)\ldots \lambda_{dia}(a_n)$ and $\lambda_{cre}(a) = \lambda_{cre}(a_1)\ldots \lambda_{cre}(a_n)$ we obtain by repetition that
$\lambda_{ann}(a)\eta \in \mathring{\calH}_{\ww\vv}$ and 
$\lambda_{dia}(a)\eta \in \mathring{\calH}_{\vv}$ and 
$\lambda_{cre}(a)\eta \in \mathring{\calH}_{\ww\vv}$ and that
	\begin{itemize}
		\item $\lambda_{ann}(a)\eta=0$ whenever $\vv$ does not start with $\ww^{-1}$ (i.e. $\ww\vv\not\in \calW^L(\ww^{-1})$)
		\item $\lambda_{dia}(a)\eta =0$ whenever  $\vv$ does not start with $w_i$ for some $1\leq i\leq d$ (equivalently: when $\vv$ does not start with $\ww$, or $\ww$ is not a clique word).
		\item $\lambda_{cre}(a)\eta =0$ whenever $\ww\vv$ does not start with $\ww$ (i.e. when $\vv\not\in \calW^L(\ww)$).
	\end{itemize}
	Last, we note that if $a\in \mathring{\boldA}_{\ww}$ and $\eta\in \mathring{\calH}_{\vv}$ are both pure tensors, then $\lambda_{ann}(a)\eta$, $\lambda_{dia}(a)\eta$ and $\lambda_{cre}(a)\eta$ are also pure tensors.\\
	
	Let $(\ww_1,\ww_2,\ww_3)\in W^3$ be s.t. $\ww:=\ww_1\ww_2\ww_3$ is a reduced expression.  We then define a linear map
	$\lambda_{(\ww_1,\ww_2,\ww_3)}:\boldA\to \calB(\calF)$ as follows.
	For a pure tensor $a\in \mathring{\boldA}_{\ww}$, there is a unique tensor $a_1\otimes a_2\otimes a_3\in \mathring{\boldA}_{\ww_1}\otimes \mathring{\boldA}_{\ww_2}\otimes \mathring{\boldA}_{\ww_3}$ s.t. $a = \calQ_{(\ww_1,\ww_2,\ww_3)}(a_1\otimes a_2\otimes a_3)$.
	We then define
	\begin{align}
	\lambda_{(\ww_1,\ww_2,\ww_3)}(a) = \lambda_{cre}(a_1)\lambda_{dia}(a_2)\lambda_{ann}(a_3)
	\end{align}
	Furthermore, we define $\lambda_{(\ww_1,\ww_2,\ww_3)}(a)=0$ for $a\in \mathring{\boldA}_{\ww'}$ with $\ww'\not=\ww_1\ww_2\ww_3$.
	
	The operator $\lambda_{(\ww_1,\ww_2,\ww_3)}(a)$ must be seen as the part of $\lambda(a)$ that acts on a vector precisely by annihilating the $\ww_3$-part, diagonally acting on a $\ww_2$-part, and creating a $\ww_1$-part.
	
	For an element $\ww\in W$, we now define the set of triple splittings
	\begin{align}
	\calS_{\ww}=\Set{(\ww_1,\ww_2,\ww_3)\in W^3| \begin{array}{l}\ww = \ww_1\ww_2\ww_3\\
		\ww_2 \text{ is a clique word}\\ |\ww|=|\ww_1|+|\ww_2|+|\ww_3|
		\end{array}}
	\end{align}
	and also define $\calS = \bigcup_{\ww\in W}\calS_{\ww}$.\\
	
		\begin{remark}\label{remark-set-Sv}
		
	We explain how the definitions of the sets $\calS_{\vv}$ relates to permutations defined in \cite[Definition 2.3]{caspersGraphProductKhintchine2021a}.
	Let $\vv=v_1\cdots v_d\in W$ be a reduced expression, let $0\leq l\leq d$, $0\leq k\leq d-l$ and let $\tbold,\uu_l,\uu_r\in W$ be clique words such that $\uu_l\tbold$, $\tbold\uu_r$ are clique words, $\uu_l\tbold\uu_r$ is reduced, and $|\tbold| = l$
	(in the notation of \cite[Definition 2.3]{caspersGraphProductKhintchine2021a} $\tbold,\uu_l,\uu_r$ correspond to the cliques $\Gamma_0,\Gamma_1,\Gamma_2$, and the conditions we put on $\tbold,\uu_l,\uu_r$ are equivalent to $\Gamma_0\in \Cliq(\Gamma,l)$ and $(\Gamma_1,\Gamma_2)\in \Comm(\Gamma_0)$).
	Then a permutation $\sigma (=\sigma_{l,k,\tbold,\uu_l,\uu_r}^{\vv})$ is defined (if existent) as the permutation such that (1) $\vv = v_{\sigma(1)}\cdots v_{\sigma(d)}$, (2) $v_{\sigma(k+1)}\cdots v_{\sigma(k+l)} = \tbold$, (3) $|v_{\sigma(1)}\cdots v_{\sigma(k)}s| = k-1$ for any letter $s$ of $\uu_l$,  (4) $|v_{\sigma(1)}\cdots v_{\sigma(k)}s|=k+1$ for any letter $s$ such that $s\uu_l\tbold$ is a reduced clique word, (5) $|sv_{\sigma(k+l+1)}\cdots v_{\sigma(d)}| = d-k-l-1$ for any letter $s$ of $\uu_r$, (6) $|s v_{\sigma(k+l+1)}\cdots v_{\sigma(d)}| = d-k-l+1$ for any letter $s$ such that $s\uu_r\tbold$ is a reduced clique word.
	Furthermore $\sigma$ is chosen  such that the expressions $\vv_1:=v_{\sigma(1)}\cdots v_{\sigma(k)}$, $\vv_2:=v_{\sigma(k+1)}\cdots v_{\sigma(k+l)}$ and $\vv_3:=v_{\sigma(k+l+1)}\cdots v_{\sigma(d)}$ are the representatives of their equivalence classes and such that $v_i=v_j$ for $i<j$ implies $\sigma(i)<\sigma(j)$. Such permutation, if existent, is unique.
	
	We make a few remarks. First of all we note that conditions (3)+(4) are equivalent to  $\sbold_r(v_{\sigma(1)}\cdots v_{\sigma(k)}\tbold) = \uu_l\tbold$, and similarly that conditions (5)+(6) are equivalent to 
	$\sbold_l(\tbold v_{\sigma(k+l+1)}\cdots v_{\sigma(d)}) = \uu_r\tbold$.
	Secondly, we note that, when $\sigma$ exists, the obtained triple $(\vv_1,\vv_2,\vv_3)$ lies in $\calS_{\vv}$. 	
	In fact, for $\vv=v_1\cdots v_d\in W$, this correspondence $(l,k,\uu_l,\uu_r,\tbold)\mapsto (\vv_1,\vv_2,\vv_3)$ between tuples $(l,k,\uu_l,\uu_r,\tbold)$ for which $\sigma_{l,k,\tbold,\uu_l,\uu_r}^{\vv}$ exists, and tuples $(\vv_1,\vv_2,\vv_3)$ in $\calS_{\vv}$, is bijective. Indeed, for $(\vv_1,\vv_2,\vv_3)\in \calS_{\vv}$ the tuple $(l,k,\uu_l,\uu_r,\tbold)$ such that the corresponding permutation $\sigma$ satisfies $\vv_1 = v_{\sigma(1)}\cdots v_{\sigma(k)}$, $\vv_2 = v_{\sigma(k+1)}\cdots v_{\sigma(k+l)}$, $\vv_3 = v_{\sigma(k+l+1)}\cdots v_{\sigma(d)}$ is given by $k = |\vv_1|$, $l = |\vv_2|$, $\tbold = \vv_2$, $\uu_l=\sbold_r(\vv_1\tbold)\tbold$, $\uu_r = \sbold_l(\tbold\vv_3)\tbold$.
\end{remark}

The following lemma was essentially proven in \cite[Lemma 2.5, Proposition 2.6]{caspersGraphProductKhintchine2021a}, and tells in what ways an element $a\in \lambda(\boldA)$ can act on a vector.

	\begin{lemma}
	\label{lemma:partition-action}
	We have that
	\begin{align}\label{eq:ways-operator-can-act}
	\lambda =\sum_{\substack{(\ww_1,\ww_2,\ww_3)\in \calS}}\lambda_{(\ww_1,\ww_2,\ww_3)}.
	\end{align}
	Moreover, $\lambda_{(\ww_1,\ww_2,\ww_3)} = 0$ whenever $\ww_2$ is not a clique word.
	In particular, for $\ww\in W$ and $a\in \mathring{\boldA}_{\ww}$ we find
	\begin{align}\label{eq:ways-operator-of-type-w-can-act}
	\lambda(a) = \sum_{(\ww_1,\ww_2,\ww_3)\in \calS_{\ww}}\lambda_{(\ww_1,\ww_2,\ww_3)}(a).
	\end{align}
\begin{proof}
	Let $\ww=w_1\cdots w_d\in W$ and $(\ww_1,\ww_2,\ww_3)\in \calS_{\ww}$ and let $\sigma$ be the corresponding permutation with $\ww_1 = w_{\sigma(1)}\cdots w_{\sigma(k)}$, $\ww_2 = w_{\sigma(k+1)}\cdots w_{\sigma(k+l)}$ and $\ww_3 = w_{\sigma(k+l+1)}\cdots w_{d}$. Then, for $a=a_1\otimes \cdots \otimes a_d\in \mathring{\boldA}_{\ww}$ we have
\begin{align}
	&\lambda_{(\ww_1,\ww_2,\ww_3)}(a) =\\
	&=\lambda_{cre}(a_{\sigma(1)}\otimes \cdots \otimes a_{\sigma(k)})\\ &\cdot\lambda_{dia}(a_{\sigma(k+1)}\otimes \cdots \otimes a_{\sigma(k+l)})\\
	&\cdot\lambda_{ann}(a_{\sigma(k+l+1)}\otimes \cdots \otimes a_{\sigma(d)})\\	 
	&= (P_{w_{\sigma(1)}}\lambda_{w_{\sigma(1)}}(a_{\sigma(1)})P_{w_{\sigma(1)}}^{\perp})\ldots  (P_{w_{\sigma(k)}}\lambda_{w_{\sigma(1)}}(a_{\sigma(k)})P_{w_{\sigma(k)}}^{\perp})\\
	&\cdot (P_{w_{\sigma(k+1)}}\lambda_{w_{\sigma(k+1)}}(a_{\sigma(k+1)})P_{w_{\sigma(k+1)}})\ldots (P_{w_{\sigma(k+l)}}\lambda_{w_{\sigma(k+l)}}(a_{\sigma(k+l)})P_{w_{\sigma(k+l)}})\\
	&\cdot (P_{w_{\sigma(k+l+1)}}^{\perp}\lambda_{w_{\sigma(m+1)}}(a_{\sigma(k+l+1)})P_{w_{\sigma(k+l+1)}})\ldots 
	(P_{w_{\sigma(d)}}^{\perp}\lambda_{w_{\sigma(d)}}(a_{\sigma(d)})P_{w_{\sigma(d)}}).
\end{align}
Equation \eqref{eq:ways-operator-of-type-w-can-act} now follows from \cite[Proposition 2.6]{caspersGraphProductKhintchine2021a} and from the bijective correspondence between the tuples $(l,k,\uu_l,\uu_r,\tbold)$ and the elements in $\calS_{\ww}$ as  described in \cref{remark-set-Sv}.
Equation \eqref{eq:ways-operator-can-act} then follows from linearity and the fact that $\lambda_{(\ww_1,\ww_2,\ww_3)}(b)=0$ whenever  $b\in \mathring{\boldA}_{\ww'}$ with $\ww'\not=\ww$.
Last, we note that by \cite[Lemma 2.5]{caspersGraphProductKhintchine2021a} we have $\lambda_{(\ww_1,\ww_2,\ww_3)}(a)=0$ whenever $\ww_2$ is not a clique word, which completes the proof.
\end{proof}

\end{lemma}

	We now prove the following
	\begin{lemma}\label{lemma:action-on-word-part}
		Let $\vv_1,\vv_2\in W$ with $|\vv_1\vv_2|=|\vv_1|+|\vv_2|$. Let $\eta\in \mathring{\calH}_{\vv_1\vv_2}$ be a pure tensor, and write $\eta = \calQ_{(\vv_1,\vv_2)}(\eta_1\otimes \eta_2)$ for some $\eta_1\otimes \eta_2\in \mathring{\calH}_{\vv_1}\otimes \mathring{\calH}_{\vv_2}$. Let $\ww\in W$ and  let $a\in \mathring{\boldA}_{\ww}$. The following holds
		\begin{enumerate}[(i)]
			\item \label{lemma:item:lambda-ann-dia} If $|\vv_1| = |\ww|+|\ww\vv_1|$ then also $|\ww\vv_1\vv_2| = |\ww\vv_1|+|\vv_2|$ and \begin{align}\lambda_{ann}(a)\eta &= \calQ_{(\ww\vv_1,\vv_2)}(\lambda_{ann}(a)\eta_1\otimes \eta_2)\\
			\lambda_{dia}(a)\eta &= \calQ_{(\vv_1,\vv_2)}(\lambda_{dia}(a)\eta_1\otimes \eta_2).
			\end{align}
			\item \label{lemma:item:lambda-cre} If $|\ww\vv_1\vv_2|=|\ww|+|\vv_1\vv_2|$  then also $|\ww\vv_1| = |\ww|+|\vv_1|$ and
			\begin{align}\lambda_{cre}(a)\eta = \calQ_{(\ww\vv_1,\vv_2)}(\lambda_{cre}(a)\eta_1\otimes \eta_2).
			\end{align}
			\item \label{lemma:item:combo} If $(\ww_1,\ww_2,\ww_3)\in \calS_{\ww}$ and if $|\vv_1| = |\ww_2\ww_3|+|\ww_2\ww_3\vv_1|$ and 
			$|\ww_1\ww_3\vv_1\vv_2| = |\ww_1|+|\ww_3\vv_1\vv_2|$, then also $|\ww_1\ww_3\vv_1\vv_2| = |\ww_1\ww_3\vv_1|+|\vv_2|$ and
			\begin{align}
			\lambda_{(\ww_1,\ww_2,\ww_3)}(a)\eta &= \calQ_{(\ww_1\ww_3\vv_1,\vv_2)}(\lambda_{(\ww_1,\ww_2,\ww_3)}(a)\eta_1\otimes \eta_2).
			\end{align}
		\end{enumerate}
		\begin{proof}
			(i) Assume that $|\vv_1|= |\ww|+|\ww\vv_1|$. Then
			\begin{align}|\vv_1\vv_2|-|\ww|\leq |\ww\vv_1\vv_2| \leq |\ww\vv_1|+|\vv_2|=|\vv_1|+|\vv_2|-|\ww| = |\vv_1\vv_2|-|\ww|.
			\end{align}
			Hence, $|\ww\vv_1\vv_2| = |\ww\vv_1|+|\vv_2|$, which proves the remark. We now prove that the equations by induction to the length $|\ww|$. First of all, it is clear that the statement holds when $\ww =e$, as then $\lambda_{ann}(a)=\lambda_{dia}(a)=a\in \CC\Id_{\mathring{\calH}_e}$. 
			
				Thus assume that $|\ww|\geq 1$ and that the statement holds for $\widetilde{\ww}$ with $|\widetilde{\ww}|\leq |\ww|-1$.
			Write $\ww = \widetilde{\ww}w$ with $\widetilde{\ww}\in W$ and $w\in V\Gamma$ and  s.t. $|\widetilde{\ww}| = |\ww|-1$.
			Then we also have $|\vv_1| = |w| + |w\vv_1|$. 
			Let us write $a = \calQ_{(\widetilde{\ww},w)}(a_1\otimes a_2)$ with $a_1\otimes a_2\in \mathring{\boldA}_{\widetilde{\ww}}\otimes \mathring{\boldA}_{w}$. Then $\lambda_{ann}(a) = \lambda_{ann}(a_1)\lambda_{ann}(a_2)$.

			Now, write $\eta = \calQ_{(w,w\vv_1,\vv_2)}(\eta_{w}\otimes \eta_1' \otimes \eta_2)$ for some $\eta_w\otimes \eta_1'\otimes \eta_2\in \mathring{\calH}_{w}\otimes \mathring{\calH}_{w\vv_1}\otimes \mathring{\calH}_{\vv_2}$ and define
			\begin{align}
			\eta' &= \calQ_{(w\vv_1,\vv_2)}(\eta_1' \otimes \eta_2)\\
			\eta_1 &= \calQ_{(w,w\vv_1)}(\eta_w \otimes \eta_1')
			\end{align}
			so that also $\eta = \calQ_{(w,w\vv_1\vv_2)}(\eta_{w}\otimes \eta')=\calQ_{(\vv_1, \vv_2)}(\eta_1\otimes \eta_2)$.
			
			We now have, using the definitions, that
			\begin{align}
				\lambda_{ann}(a_2)\eta &=P_w^\perp\lambda_w( a_2)P_w\eta=P_w^{\perp}U_w((a_2\eta_w)\otimes \eta')
				= \langle a_2\eta_w,\xi_{w}\rangle\eta'\\
				\lambda_{ann}(a_2)\eta_1 &=P_w^\perp \lambda_w(a_2)P_w\eta_1
				=P_w^{\perp}U_w((a_2\eta_w)\otimes \eta_1')
				= \langle a_2\eta_w,\xi_{w}\rangle\eta_1'
			\end{align}
			and 
			\begin{align}
				\lambda_{dia}(a_2)\eta &=P_wU_w((a_2\eta_w)\otimes \eta')
				=\calQ_{(w,w\vv_1\vv_2)}(\mathring{(a\eta_w)}\otimes \eta' )\\
				\lambda_{dia}(a_2)\eta_1 &=P_wU_w((a_2\eta_w)\otimes \eta_1')
				=\calQ_{(w,w\vv_1)}(\mathring{(a\eta_w)}\otimes \eta_1').
			\end{align}
			Now this means that
			\begin{align}
					\lambda_{ann}(a_2)\eta &=\varphi_w(a_2\eta_w)\eta'\\
					&=\calQ_{(w\vv_1,\vv_2)}(\langle a_2\eta_w,\xi_{w}\rangle\eta_1' \otimes \eta_2)\\
					&=\calQ_{(w\vv_1,\vv_2)}(\lambda_{ann}(a_2)\eta_1\otimes \eta_2)
			\end{align}
			and 
				\begin{align}
				\lambda_{dia}(a_2)\eta &= \calQ_{(w,w\vv_1\vv_2)}(\mathring{(a_2\eta_w)}\otimes \eta' )\\
				&= \calQ_{(w,w\vv_1,\vv_2)}(\mathring{(a_2\eta_w)}\otimes \eta_1' \otimes \eta_2 )\\
				&= \calQ_{(\vv_1,\vv_2)}( \calQ_{(w,w\vv_1)}(\mathring{(a_2\eta_w)}\otimes \eta_1' )\otimes \eta_2 )\\
				&= \calQ_{(\vv_1,\vv_2)}(\lambda_{dia}(a_2)\eta_1 \otimes \eta_2 ).
			\end{align}
			Now, we note that $|w\vv_1| = |\widetilde{\ww}|+|\widetilde{\ww}w\vv_1|$ so that using the induction hypothesis and the fact that $|\widetilde{\ww}|=|\ww|-1$ we find
			\begin{align}
		\lambda_{ann}(a)\eta &= \lambda_{ann}(a_1)\lambda_{ann}(a_2)\eta \\
		&=\lambda_{ann}(a_1)\calQ_{(w\vv_1,\vv_2)}(\lambda_{ann}(a_2)\eta_1 \otimes \eta_2 )\\
			&= \calQ_{(\widetilde{\ww}w\vv_1,\vv_2)}(\lambda_{ann}(a_1)\lambda_{ann}(a_2)\eta_1 \otimes \eta_2 )\\
			&= \calQ_{(\ww\vv_1,\vv_2)}(\lambda_{ann}(a)\eta_1 \otimes \eta_2 ).
				\end{align}
			Similarly
			\begin{align}
				\lambda_{dia}(a)\eta  &= \lambda_{dia}(a_1)\lambda_{dia}(a_2)\eta\\
				&=\lambda_{dia}(a_1)\calQ_{(\vv_1,\vv_2)}(\lambda_{dia}(a_2)\eta_1 \otimes \eta_2 )\\
				&=\calQ_{(\vv_1,\vv_2)}(\lambda_{dia}(a_1)\lambda_{dia}(a_2)\eta_1 \otimes \eta_2 )\\	&=\calQ_{(\vv_1,\vv_2)}(\lambda_{dia}(a)\eta_1 \otimes \eta_2 ).
			\end{align}
			This finishes the induction, and proves the statement.\\
			
			(ii) Assume that $|\ww\vv_1\vv_2| = |\ww|+|\vv_1\vv_2|$. Then 
			\begin{align}
				|\ww\vv_1\vv_2| \leq |\ww\vv_1|+|\vv_2|\leq |\ww|+|\vv_1|+|\vv_2| = |\ww|+|\vv_1\vv_2|
				= |\ww\vv_1\vv_2|.
			\end{align}
			Hence $|\ww\vv_1|=|\ww|+|\vv_1|$, which shows the first remark. 	Again we  prove the equation by induction to the length $|\ww|$. Again, it is clear that the statement holds when $\ww =e$.
			Thus assume that $|\ww|\geq 1$ and that the statement holds for $\widetilde{\ww}$ with $|\widetilde{\ww}|\leq |\ww|-1$.
			Write $\ww = \widetilde{\ww}w$ with $\widetilde{\ww}\in W$ and $w\in V\Gamma$ and  s.t. $|\widetilde{\ww}| = |\ww|-1$.
			Then we also have $|w\vv_1\vv_2| = |w| + |\vv_1\vv_2|$.
			Let us write $a = \calQ_{(\widetilde{\ww},w)}(a_1\otimes a_2)$ with $a_1\otimes a_2\in \mathring{\boldA}_{\widetilde{\ww}}\otimes \mathring{\boldA}_{w}$. Then $\lambda_{cre}(a) = \lambda_{cre}(a_1)\lambda_{cre}(a_2)$.
			
			We now have by definition
			\begin{align}
				\lambda_{cre}(a_2)\eta &=P_w\lambda_w(a_2)P_w^\perp\eta =(P_wU_{w})((a_2 \xi_w) \otimes \eta )
				=\calQ_{(w,\vv_1\vv_2)}(\widehat{a}_2\otimes \eta)\\
				\lambda_{cre}(a_2)\eta_1 &=P_w\lambda_w(a_2)P_w^\perp\eta_1
				=(P_wU_w)((a_2\xi_w)\otimes \eta_1)
					=\calQ_{(w,\vv_1)}(\widehat{a}_2\otimes \eta_1).
			\end{align}
			Now this means that
			\begin{align}
				\lambda_{cre}(a_2)\eta &=\calQ_{(w,\vv_1\vv_2)}(\widehat{a}_2\otimes \eta)\\
				&=\calQ_{(w,\vv_1,\vv_2)}(\widehat{a}_2\otimes \eta_1\otimes \eta_2)\\
				&=\calQ_{(w\vv_1,\vv_2)}(\calQ_{(w,\vv_1)}(\widehat{a}_2\otimes \eta_1)\otimes \eta_2)\\
				&=\calQ_{(w\vv_1,\vv_2)}(\lambda_{cre}(a_2)\eta_1\otimes \eta_2).
			\end{align}
			Now, we note that $|\widetilde{\ww}w\vv_1\vv_2| = |\widetilde{\ww}|+|w\vv_1\vv_2|$ so that using the induction hypothesis and the fact that $|\widetilde{\ww}|=|\ww|-1$ we find
			\begin{align}
				\lambda_{cre}(a)\eta &= \lambda_{cre}(a_1)\lambda_{cre}(a_2)\eta \\
				&=\lambda_{cre}(a_1)\calQ_{(w\vv_1,\vv_2)}(\lambda_{cre}(a_2)\eta_1 \otimes \eta_2 )\\
				&= \calQ_{(\widetilde{\ww}w\vv_1,\vv_2)}(\lambda_{cre}(a_1)\lambda_{cre}(a_2)\eta_1 \otimes \eta_2 )\\
				&= \calQ_{(\ww\vv_1,\vv_2)}(\lambda_{cre}(a)\eta_1 \otimes \eta_2 ).
			\end{align}
			This finishes the induction, and proves the statement.\\		
			
			(iii) Let $(\ww_1,\ww_2,\ww_3)\in \calS_{\ww}$ be s.t $|\vv_1| = |\ww_2\ww_3|+|\ww_2\ww_3\vv_1|$ and $|\ww_1\ww_3\vv_1\vv_2| = |\ww_1|+|\ww_3\vv_1\vv_2|$. We will write $\lambda_{(\ww_1,\ww_2,\ww_3)}(a) = \lambda_{cre}(a_1)\lambda_{dia}(a_2)\lambda_{ann}(a_3)$ for some $a_1\otimes a_2\otimes a_3\in \mathring{\boldA}_{\ww_1}\otimes \mathring{\boldA}_{\ww_2}\otimes \mathring{\boldA}_{\ww_3}$. 
			Now, first, as $|\vv_1| = |\ww_2\ww_3|+|\ww_2\ww_3\vv_1|$, we also have
			\begin{align}
			|\vv_1|&\leq |\ww_3|+|\ww_3\vv_1|\\
			&\leq|\ww_2|+|\ww_3|+|\ww_2\ww_3\vv_1|\\
			&=|\ww_2\ww_3|+|\ww_2\ww_3\vv_1| \\
			&= |\vv_1|
			\end{align}
			and therefore $|\vv_1|= |\ww_3|+|\ww_3\vv_1|$.
			By \eqref{lemma:item:lambda-ann-dia} this gives us 
			\begin{align}
			\lambda_{ann}(a_3)\eta = \calQ_{(\ww_3\vv_1,\vv_2)}(\lambda_{ann}(a_3)\eta_1\otimes \eta_2)
			\end{align}
			and also $|\ww_3\vv_1\vv_2|=|\ww_3\vv_1|+|\vv_2|$.
			Now, we also find 
			\begin{align}
			|\ww_3\vv_1| = |\vv_1|-|\ww_3| = |\ww_2\ww_3|+|\ww_2\ww_3\vv_1|-|\ww_3| = |\ww_2|+|\ww_2\ww_3\vv_1|.
			\end{align}
			Let us set $\vv_1'=\ww_3\vv_1$ and $\vv_2'=\vv_2$, so that $|\vv_1'\vv_2'|=|\vv_1'|+|\vv_2'|$ and $|\vv_1'|=|\ww_2|+|\ww_2\vv_1'|$.
			Moreover set $\eta' = \lambda_{ann}(a_3)\eta$ and $\eta_1' = \lambda_{ann}(a_3)\eta_1$ and $\eta_2'=\eta_2$.
			Now $\eta'= \calQ_{(\vv_1',\vv_2')}(\eta_1'\otimes \eta_2')$ and we see that the conditions for applying \eqref{lemma:item:lambda-ann-dia} are satisfied. This thus gives us that
			\begin{align}
			\lambda_{dia}(a_2)\lambda_{ann}(a_3)\eta = \calQ_{(\ww_3\vv_1,\vv_2)}(	\lambda_{dia}(a_2)\lambda_{ann}(a_3)\eta_1\otimes \eta_2).
			\end{align}
			Now, set $\widetilde{\vv}_1 = \vv_1'=\ww_3\vv_1$ and $\widetilde{\vv}_2=\vv_2'=\vv_2$  so that again $|\widetilde{\vv}_1\widetilde{\vv}_2| = |\widetilde{\vv}_1|+|\widetilde{\vv}_2|$. Also we get $|\ww_1\widetilde{\vv}_1\widetilde{\vv}_2|=|\ww_1\ww_3\vv_1\vv_2|=|\ww_1\ww_3\vv_1|+|\vv_2|=|\ww_1\widetilde{\vv}_1|+|\widetilde{\vv}_2|$.
			Also set  $\widetilde{\eta} = 	\lambda_{dia}(a_2)\lambda_{ann}(a_3)\eta$ and 
			$\widetilde{\eta}_1 = \lambda_{dia}(a_2)\lambda_{ann}(a_3)\eta_1$ and 
			$\widetilde{\eta}_2 = \eta_2$
			Then $\widetilde{\eta}= \calQ_{(\vv_1',\vv_2')}(\widetilde{\eta}_1\otimes \widetilde{\eta}_2)$ and all conditions for applying \eqref{lemma:item:lambda-cre} are satisfied. By \eqref{lemma:item:lambda-cre} we thus get
			\begin{align}
			\lambda_{cre}(a_1)\lambda_{dia}(a_2)\lambda_{ann}(a_3)\eta = \calQ_{(\ww_1\ww_3\vv_1,\vv_2)}(\lambda_{cre}(a_1)	\lambda_{dia}(a_2)\lambda_{ann}(a_3)\eta_1\otimes \eta_2)
			\end{align}
			and moreover $|\ww_1\ww_3\vv_1| = |\ww_1|+|\ww_3\vv_1|$. The previous equation is precisely what we needed to show, and we moreover obtain			
			$|\ww_1\ww_3\vv_1\vv_2| = |\ww_1|+|\ww_3\vv_1\vv_2|= |\ww_1|+|\ww_3\vv_1|+|\vv_2| = |\ww_1\ww_3\vv_1|+|\vv_2|$, which proves the statement.
		\end{proof}
	\end{lemma}

\section{polynomial growth of word-length projections}\label{section:polynomial-growth-word-length}
In this section we shall fix a simple \textbf{finite} graph $\Gamma$, together with  unital $\Cstar$-algebras $\boldA_v$ for $v\in \Gamma$ and states $\varphi_v$ on $\boldA_v$ for which the GNS representation is faithful. We shall look at the reduced graph product $(\calA,\varphi) = *_{v,\Gamma}(\boldA_{v},\varphi_v)$ and investigate for $d\geq 0$ the natural projections $\calP_d:\calA\to \calA_d$.
	The main result of this section, \Cref{thm:projection-maps-linear-growth}, is that these maps are completely bounded, and that we can obtain a bound on $\|\calP_{d}\|_{cb}$ that depends only linearly on $d$. To prove this, we can not use the same method as \cite{ricardKhintchineTypeInequalities2006a}, since that relies on the fact that each element either does not act diagonally on a pure tensor $\eta\in \mathring{\calH}_{\vv}\subseteq \calF$, or acts diagonally on $\eta$ on precisely one letter. This holds true for elements in the free product, but not generally for elements in the graph product, as they may act diagonally on any clique. Therefore, we will instead introduce completely contractive maps $H_{\tau}$ (and completely bounded maps $\widetilde{H}_{\rho}$) and write $\calP_d$ as linear combination of these.
	
\subsection{The maps $H_\tau$}
We introduce some extra notation.
Let $W$ be the right-angled Coxeter group associated to the graph $\Gamma$.  Recall, for a word $\ww\in W$ we defined $\sbold_l(\ww)$ and $\sbold_r(\ww)$ as the maximal clique words that $\ww$ respectively starts with and ends with.
 For a word $\uu\in W$, $n\geq 0$,  we define
 \begin{align}
 \calW^L(\uu) &= \{\ww\in W: |\uu\ww|=|\uu|+|\ww|\}\\
 \calW^R(\uu) &= \{\ww\in W: |\ww\uu|=|\ww|+|\uu|\}\\
 \widetilde{\calW}^L(\uu) &= \{\ww\in \calW^L(\uu): \sbold_l(\uu\ww)=\sbold_l(\uu)\}\\
 \widetilde{\calW}^R(\uu) &= \{\ww\in \calW^R(\uu): \sbold_r(\ww\uu)=\sbold_r(\uu)\}\\
 \widetilde{\calW}_n^L(\uu) &= \{\ww\in \widetilde{\calW}^L(\uu): |\ww|=n\}\\
 \widetilde{\calW}_n^R(\uu) &= \{\ww\in \widetilde{\calW}^R(\uu): |\ww|=n\}.
 \end{align}
 Now, let $\uu\in W$ and let $\uu_L,\uu_R\in W$ be s.t. $|\uu|=|\uu\uu_L^{-1}|+|\uu_L|$ and $|\uu|=|\uu_R|+|\uu_R^{-1}\uu|$, i.e. $\uu_L$ is some word that $\uu$ ends with and $\uu_R$ is some word that $\uu$ starts with. Then we have for $\ww_L\in \calW^L(\uu)$ and $\ww_R\in \calW^R(\uu)$ that $\uu_L\ww_L$ and $\ww_R\uu_R$ are reduced expressions.
 Let $n\geq 0$. We define
 \begin{align}
 \calH^L(\uu,\uu_L) &= \bigoplus_{\ww\in \calW^L(\uu)}\mathring{\calH}_{\uu_L\ww}&
 \calH^R(\uu,\uu_R) &= \bigoplus_{\ww\in \calW^R(\uu)}\mathring{\calH}_{\ww\uu_R}
 \\
 \calF^L(\uu,\uu_L) &= \bigoplus_{\ww\in \widetilde{\calW}^L(\uu)}\mathring{\calH}_{\uu_L\ww}&
 \calF^R(\uu,\uu_R) &= \bigoplus_{\ww\in \widetilde{\calW}^R(\uu)}\mathring{\calH}_{\ww\uu_R}
 \\
 \calF_n^L(\uu,\uu_L) &= \bigoplus_{\substack{\ww\in \widetilde{\calW}_n^L(\uu)}}\mathring{\calH}_{\uu_L\ww}&
 \calF_n^R(\uu,\uu_R) &= \bigoplus_{\substack{\ww\in \widetilde{\calW}_n^R(\uu)}}\mathring{\calH}_{\ww\uu_R}.
 \end{align}
 For $\uu\in W$ and $n\geq 0$ we moreover define
 \begin{align}
 \calF_n^M(\uu) = \bigoplus_{\substack{\ww_1 \in \widetilde{\calW}_{n}^R(\uu)\\ 
 		\ww_2 \in \calW^L(\uu)}}\mathring{\calH}_{\ww_1\uu\ww_2}.
 \end{align}
We note that for  $\ww_1 \in \widetilde{\calW}_{n}^R(\uu)$ and 
$\ww_2 \in \calW^L(\uu)$ we have that  $\ww_1\uu\ww_2$ is a reduced expression. Indeed, it is clear that $\ww_1\uu$ and $\uu\ww_2$ are reduced by definition. Now, since moreover $\sbold_r(\ww_1\uu)=\sbold_r(\uu)$, we have that no letter from $\ww_1$ can cancel out a letter of $\ww_2$, so that the expression is reduced.
 \begin{definition}
 	Let $\uu\in W$ and let $\rr\in W$ be any clique word that $\uu$ ends with. Then $\uu\rr$ is a word in $W$ that $\uu$ starts with, and $|\uu\rr| +|\rr| = |\uu|$. For $n\geq 0$ we define a partial isometry $V_{n}^{\uu,\rr}:\calF\otimes \calF\to \calF$ with initial subspace $\calF_n^R(\uu,\uu\rr)\otimes \calH^L(\uu,\rr)$ and final subspace $\calF_n^M(\uu)$ as
 	\begin{align}
 	V_{n}^{\uu,\rr}|_{\mathring{\calH}_{\vv_r\uu\rr}\otimes \mathring{\calH}_{\rr\vv_{tail}}} = \calQ_{(\vv_r\uu\rr, \rr\vv_{tail})}&\quad \text{ for } \vv_r\in \widetilde{\calW}_n^R(\uu), \vv_{tail}\in \calW^L(\uu).
 	\end{align}
 \end{definition}
 We note that this is well-defined.
 Indeed, as just pointed out, for  $\vv_r\in \widetilde{\calW}_n^R(\uu)$ and $\vv_{tail}\in \calW^L(\uu)$ we have that $\vv_r\uu\vv_{tail}$ is reduced. Therefore, we get $|\vv_r\uu\vv_{tail}|\leq |\vv_r\uu\rr|+|\rr\vv_{tail}| \leq |\vv_r|+|\uu\rr|+|\rr|+|\vv_{tail}|=|\vv_r|+|\uu|+|\vv_{tail}| = |\vv_r\uu\vv_{tail}|$.
 This shows that $|\vv_r\uu\rr| + |\rr\vv_{tail}| = |\vv_r\uu\vv_{tail}|$, so that $\calQ_{(\vv_r\uu\rr,\rr\vv_{tail})}$ is well-defined.
 
 \begin{definition}\label{definition:clique-pairs}
 	We denote 
 	\begin{align}
 	\calT &= \Set{(\uu_l,\uu_r,\tbold) \in W^3| \begin{array}{l}
 		\uu_l\tbold, \tbold\uu_r \text{ clique words}, \\
 		\uu_l\tbold\uu_r \text{ reduced}\\
 		\end{array}}.
 	\end{align}
 	We remark that it follows from the definition that $\uu_l,\uu_r$ and $\tbold$ must also be clique words and that $\uu_l\uu_r$ must be reduced.
 \end{definition}

 \begin{definition}
 	Let $(\uu_l,\uu_r,\tbold)\in \calT$.
 	Also let $\rr\in W$ be a sub-clique word of $\tbold$ and let $n_l,n_r\geq 0$. For the tuple $\tau = (n_l,n_r,\uu_l,\uu_r,\tbold,\rr)$ define a map $H_{\tau}:\calB(\calF)\to \calB(\calF)$ as
 	\begin{align}
 	H_{\tau}(a) = V_{n_l}^{(\uu_l\tbold),\rr}(a\otimes \Id_{\calF})\left(V_{n_r}^{(\uu_r\tbold),\rr}\right)^*.
 	\end{align}
 \end{definition}	
	It is clear that $H_{\tau}$ is completely contractive.
 
 \begin{example}\label{example:projection-length-zero}
 	We note that the partial isometry $V_{0}^{e,e}:\calF\otimes \calF\to \calF$ has initial subspace $\calF_0^R(e,e)\otimes \calH^L(e,e) = \CC\Omega\otimes \calF$ and final subspace $\calF_0^M(e) = \calF$ and that on $\CC\Omega\otimes \calF$ it is given by
 	$V_{0}^{e,e}(z\Omega\otimes \eta) = z\eta$ for $z\in \CC$, $\eta\in \calF$. Setting  $\tau = (0,0,e,e,e,e)$ and letting $a\in \boldA$ be a pure tensor $a=a_1\otimes \cdots \otimes a_t$, we can for $\eta\in \calF$ calculate 
 	$H_\tau(\lambda(a))\eta = V_{0}^{e,e}(\lambda(a)\Omega \otimes \eta)$. Now, if $\lambda(a)\Omega \not\in \CC\Omega$, then we get $H_{\tau}(\lambda(a))\eta =0$. On the other hand, if $\hat{a}=\lambda(a)\Omega \in \CC\Omega$, then we must have that $\lambda(a) \in \CC\Id_{\calF}$ and we get $H_{\tau}(a)\eta = a\eta$. We conclude that $\calP_{0} = H_{(0,0,e,e,e,e)}$ and $\|\calP_{0}\|_{cb}=1$.
 \end{example}
 Similarly to \Cref{example:projection-length-zero}, we aim to write $\calP_{d}$ for $d\geq 1$ as a linear combination of $H_{\tau}$'s for different tuples $\tau$, in order to give a bound on $\|\calP_{d}\|_{cb}$. To achieve this, we introduce some convenient notation.

	\begin{definition}
		Let $\calH_{1}$ and $\calH_2$ be closed subspaces of $\calF$. For an operator $b\in \calB(\calF)$ we define a closed subspace $\calJ_b(\calH_{1},\calH_{2})$ of $\calF$ as
		\begin{align}
		\calJ_{b}(\calH_1,\calH_2) = \{\eta\in \calH_{1}|  b\eta\in \calH_{2}\}.
		\end{align}
	\end{definition}

	\begin{proposition}\label{prop:formula-for-H-tau}
		Let $(\uu_l,\uu_r,\tbold)\in \calT$. Also let $\rr\subseteq \tbold$ be a sub-clique, and let $n_l,n_r\geq 0$. Set $\tau = (n_l,n_r,\uu_l,\uu_r,\tbold,\rr)$. For $\ww\in W$ and $\omega=(\ww_1,\ww_2,\ww_3)\in \calS_{\ww}$ and for pure tensor $a=a_1\otimes \cdots \otimes a_t\in \mathring{\boldA}_{\ww}$ we have that
		\begin{align}
		H_{\tau}(\lambda_{\omega}(a)) = \lambda_{\omega}(a)P_{a}(\tau,\omega)
		\end{align}
		where $P_{a}(\tau,\omega)$ is the projection in $\calB(\calF)$ on the closed subspace spanned by
		\begin{align}
		\bigcup_{
			\substack{\vv_l\in \widetilde{\calW}_{n_l}^R(\uu_l\tbold), 
			\vv_r\in \widetilde{\calW}_{n_r}^R(\uu_r\tbold)\\
			\vv_{tail}\in \calW^L(\uu_l\tbold)\cap \calW^L(\uu_r\tbold)\\
		|\vv_r\uu_r\tbold\rr|= |\ww_2\ww_3|+|\ww_2\ww_3\vv_r\uu_r\tbold\rr|\\
		|\ww_1\ww_3\vv_r\uu_r\tbold\vv_{tail}| = |\ww_1|+|\ww_3\vv_r\uu_r\tbold\vv_{tail}|}}
		\calJ_{\lambda_{\omega}(a)}(\mathring{\calH}_{\vv_r\uu_r\tbold\vv_{tail}},\mathring{\calH}_{\vv_l\uu_l\tbold\vv_{tail}}).
		\end{align}
		
		\begin{proof} We show that the identity holds on pure tensors.
		First, let $\vv\in W$ and let $\eta\in \mathring{\calH}_{\vv}\subseteq \calF$ be a pure tensor  s.t. $\lambda_{\omega}(a)P_a(\tau,\omega)\eta = 0$. If $\eta\perp \calF_{n_r}^M(\uu_r\tbold)$, then clearly $(V_{n_r}^{\uu_r\tbold,\rr})^*\eta =0$ so that $H_{\tau}(\lambda_{\omega}(a))\eta=0= \lambda_{\omega}(a)P_{a}(\tau,\omega)\eta$, and we are done.
		Thus, assume that $\eta\in \calF_{n_r}^M(\uu_r\tbold)$ and $\eta\not=0$, so that $\eta\in \mathring{\calH}_{\vv_r\uu_r\tbold\vv_{tail}}$ for some $
		\vv_r\in \widetilde{\calW}_{n_r}^R(\uu_r\tbold)$, 
		$\vv_{tail}\in \calW^L(\uu_r\tbold)$. Let us write $V_{n_r}^{\uu_r\tbold,\rr*}\eta = \eta_1\otimes \eta_2$ with $\eta_1\in \mathring{\calH}_{\vv_r\uu_r\tbold\rr}$, $\eta_2\in  \mathring{\calH}_{\rr\vv_{tail}}$.	Then $H_{\tau}(\lambda_{\omega}(a))\eta = V_{n_l}^{\uu_l\tbold,\rr}(\lambda_{\omega}(a)\eta_1\otimes \eta_2)$.
		We can assume that $0\not=\lambda_{\omega}(a)\eta_1\in \calF_{n_l}^R(\uu_l\tbold,\uu_l\tbold\rr)$ and $\eta_2\in \calH^L(\uu_l\tbold,\rr)$ since otherwise we find directly $H_{\tau}(\lambda_{\omega}(a))\eta =0$. 
		Now we thus have that $\lambda_{\omega}(a)\eta_1\in \mathring{\calH}_{\vv_l\uu_l\tbold\rr}$ for some $\vv_l\in \widetilde{\calW}_{n_l}^R(\uu_l\tbold)$ and that $\eta_2\in \mathring{\calH}_{\rr\vv_{tail}'}$ for some $\vv_{tail}'\in \calW_{n_r}^L(\uu_l\tbold)$. 
		
		As $\eta_2$ is non-zero, and as $\eta_2\in \mathring{\calH}_{\rr\vv_{tail}}\cap \mathring{\calH}_{\rr\vv_{tail}'}$ we find that $\vv_{tail}=\vv_{tail}'\in \calW^L(\uu_l\tbold)\cap \calW^L(\uu_r\tbold)$.
		Also, since $\eta_1\in \mathring{\calH}_{\vv_r\uu_r\tbold\rr}$ we find that  $\lambda_{\omega}(a)\eta_1\in \mathring{\calH}_{\ww_1\ww_3\vv_r\uu_r\tbold\rr}$. Now, we already had $\lambda_{\omega}(a)\eta_1\in \mathring{\calH}_{\vv_l\uu_l\tbold\rr}$ and by the assumption that $\lambda_{\omega}(a)\eta_1$ is non-zero, we thus find $\vv_l\uu_l\tbold\rr = \ww_1\ww_3\vv_r\uu_r\tbold\rr$.
		 Moreover, as $\lambda_\omega(a)\eta_1$ is non-zero, we must have that $|\vv_r\uu_r\tbold\rr|=|\ww_2\ww_3|+|\ww_2\ww_3\vv_r\uu_r\tbold\rr|$ and $|\ww_1\ww_3\vv_r\uu_r\tbold\rr| = |\ww_1|+|\ww_3\vv_r\uu_r\tbold\rr|$
		
		Set $\vv_1 = \vv_r\uu_r\tbold\rr$ and $\vv_2 = \rr\vv_{tail}$, so that $|\vv_1\vv_2|=|\vv_1|+|\vv_2|$, and by the above
		\begin{align}\label{eq:annihilation-reduced-word}
			 |\vv_1|&=|\ww_2\ww_3|+|\ww_2\ww_3\vv_1|\\ |\ww_1\ww_3\vv_1|&=|\ww_1|+|\ww_3\vv_1|
		\end{align}
		Moreover, we now find
		\begin{align}
		|\ww_1\ww_3\vv_1\vv_2|&\leq |\ww_1|+|\ww_3\vv_1\vv_2|\\
		&\leq |\ww_1|+|\ww_3\vv_1|+|\vv_2|\\
		&= |\ww_1\ww_3\vv_1|+|\vv_2|\\
		&=|\ww_1\ww_3\vv_r\uu_r\tbold\rr| +|\rr\vv_{tail}| \\
		&=|\vv_l\uu_l\tbold\rr|+|\rr\vv_{tail}|\\ 
		&= |\vv_l\uu_l\tbold\vv_{tail}|\\
		&= |\ww_1\ww_3\vv_r\uu_r\tbold\vv_{tail}|\\
		&= |\ww_1\ww_3\vv_1\vv_2|.
		\end{align}
		This shows that 
		\begin{align}\label{eq:creation-is-reduced}
			|\ww_1\ww_3\vv_1\vv_2|=|\ww_1|+|\ww_3\vv_1\vv_2|
		\end{align}
		Now as $\eta\in \mathring{\calH}_{\vv_1\vv_2}$, and as all conditions of  \Cref{lemma:action-on-word-part}\eqref{lemma:item:combo} are satisfied, this gives us
		\begin{align}
		H_{\tau}(\lambda_{\omega}(a))\eta &= V_{n_l}^{\uu_l\tbold,\rr}(\lambda_\omega(a)\eta_1\otimes \eta_2)\\
		&= \calQ_{(\ww_1\ww_3\vv_1,\vv_2)}(\lambda_{\omega}(a)\eta_1\otimes \eta_2)\\
		&=\lambda_{\omega}(a)\calQ_{(\ww_1\ww_3\vv_1,\vv_2)}(\eta_1\otimes \eta_2)\\
		&=\lambda_{\omega}(a)\eta.
		\end{align}
		Moreover we find $\lambda_{\omega}(a)\eta \in \mathring{\calH}_{\ww_1\ww_3\vv_1\vv_2} = \mathring{\calH}_{\vv_l\uu_l\tbold\vv_{tail}}$.
		However, this shows that
		$\eta \in  \calJ_{\lambda_{\omega(a)}}(\mathring{\calH}_{\vv_l\uu_l\tbold\vv_{tail}}, \mathring{\calH}_{\vv_r\uu_r\tbold\vv_{tail}})$. By all the conditions we have shown for $\vv_l,\vv_r,\vv_{tail}$, and as we have shown that $|\vv_1|=|\ww_2\ww_3|+|\ww_2\ww_3\vv_1|$ (\Cref{eq:annihilation-reduced-word}) and  $|\ww_1\ww_3\vv_1\vv_2|=|\ww_1|+|\ww_3\vv_1\vv_2|$ (\Cref{eq:creation-is-reduced}) it follows that $P_{a}(\tau,\omega)\eta = \eta$. We conclude that 
		$H_{\tau}(\lambda_\omega(a))\eta = \lambda_{\omega}(a)\eta = \lambda_{\omega}(a)P_{a}(\tau,\omega)\eta$.\\
				
		Alternatively, let $\eta\in \mathring{\calH}_{\vv}\subseteq \calF$ be a pure vector s.t. $\lambda_{\omega}(a)P_{a}(\tau,\omega)\eta\not=0$. Then we must have that $P_{a}(\tau,\omega)\eta = \eta$ and moreover that $\lambda_{\omega}(a)\eta$ is non-zero. We thus get that $\eta \in \calJ_{\lambda_{\omega}(a)}(\mathring{\calH}_{\vv_r\uu_r\tbold\vv_{tail}},\mathring{\calH}_{\vv_l\uu_l\tbold\vv_{tail}})$ with 
		$\vv_l\in \widetilde{\calW}_{n_l}^R(\uu_l\tbold)$,
		$\vv_r\in \widetilde{\calW}_{n_r}^R(\uu_r\tbold)$,
		$\vv_{tail}\in \calW^L(\uu_r\tbold)\cap \calW^L(\uu_l\tbold)$
		and so that 	
		\begin{align}
		|\vv_r\uu_r\tbold\rr|&= |\ww_2\ww_3|+|\ww_2\ww_3\vv_r\uu_r\tbold\rr|\\
		|\ww_1\ww_3\vv_r\uu_r\tbold\vv_{tail}| &
		= |\ww_1|+|\ww_3\vv_r\uu_r\tbold\vv_{tail}|.
		\end{align}
		Set $\vv_1 = \vv_r\uu_r\tbold\rr$ and $\vv_2 = \rr\vv_{tail}$, so that $|\vv_1\vv_2| = |\vv_1|+|\vv_2|$. Moreover the above equations state that $|\vv_1| = |\ww_2\ww_3|+|\ww_2\ww_3\vv_1|$
		and $|\ww_1\ww_3\vv_1\vv_2|=|\ww_1|+|\ww_3\vv_1\vv_2|$. 
		As $\eta\in \mathring{\calH}_{\vv_r\uu_r\tbold\vv_{tail}}\subseteq \calF^M_{n_r}(\uu_r\tbold)$, we can write
		$V_{n_r}^{\uu_r\tbold,\rr*}\eta = \eta_1\otimes \eta_2\in \mathring{\calH}_{\vv_r\uu_r\tbold\rr}\otimes \mathring{\calH}_{\rr\vv_{tail}}=\mathring{\calH}_{\vv_1}\otimes \mathring{\calH}_{\vv_2}$.  
		By the above properties we get from \Cref{lemma:action-on-word-part}\eqref{lemma:item:combo} that
		\begin{align}
		\lambda_{\omega}(a)\eta = \calQ_{(\ww_1\ww_3\vv_1,\vv_2)}(\lambda_{\omega}(a)\eta_1\otimes \eta_2)\in \mathring{\calH}_{\ww_1\ww_3\vv_1\vv_2}.
		\end{align}
		However, we also know that $	\lambda_{\omega}(a)\eta\in \mathring{\calH}_{\vv_l\uu_l\tbold\vv_{tail}}$. Therefore, as $\lambda_{\omega}(a)\eta$ is non-zero we find $\vv_l\uu_l\tbold\vv_{tail} = \ww_1\ww_3\vv_1\vv_2 = \ww_1\ww_3\vv_r\uu_r\tbold\vv_{tail}$. We thus find 
		$\vv_l\uu_l\tbold\rr = \ww_1\ww_3\vv_r\uu_r\tbold\rr=\ww_1\ww_3\vv_1$, and hence $\lambda_{\omega}(a)\eta_1\in \mathring{\calH}_{\ww_1\ww_3\vv_1} = \mathring{\calH}_{\vv_l\uu_l\tbold\rr}\subseteq \calF_{n_l}^R(\uu_l\tbold,\uu_l\tbold\rr)$. Note moreover that $\eta_2\in \calH^L(\uu_l\tbold,\rr)$ by the assumption on $\vv_{tail}$.
		
		Hence, as $\lambda_\omega(a)\eta_1\otimes \eta_2\in \calF_{n_l}^R(\uu_l\tbold,\uu_l\tbold\rr)\otimes \calH^L(\uu_l\tbold,\rr)$ we find that
		\begin{align}
		H_{\tau}(\lambda_{\omega}(a))\eta &= V_{n_l}^{\uu_l\tbold,\rr}(\lambda_{\omega}(a)\eta_1\otimes \eta_2)\\ 
		&=\calQ_{(\ww_1\ww_3\vv_1,\vv_2)}(\lambda_{\omega}(a)\eta_1\otimes \eta_2)\\
		&=\lambda_{\omega}(a)\eta\\
		&=\lambda_{\omega}(a)P_a(\tau,\omega)\eta
		\end{align}
		which proves the statement.
		\end{proof}
	\end{proposition}

\subsection{The maps $\widetilde{H}_\rho$}
We shall now introduce other maps, $\widetilde{H}_{\rho}$, that are linear combinations of the maps $H_{\tau}$ for different $\tau$'s, and that satisfy a nice equation. We use these maps to show that $\calP_d$ is completely bounded, and give a bound on $\|\calP_d\|_{cb}$.
\begin{definition}\label{def:sets-Swrho}
	Let $n_l,n_r\geq 0$ and $(\uu_l,\uu_r,\tbold)\in \calT$.  For $\ww\in W$ and for the tuple $\rho = (n_l,n_r,\uu_l,\uu_r,\tbold)$ define the set
	\begin{align}
	\calS_{\ww}(\rho) &= \Set{(\ww_1,\ww_2,\ww_3)\in \calS_{\ww}| \begin{array}{l}
		\ww_1 = \vv_l\uu_l, \ww_2 =\tbold\text{ and } \ww_3 =\uu_r^{-1}\vv_r^{-1} \\
		\text{ for some } \vv_l\in \widetilde{\calW}_{n_l}^R(\uu_l\tbold), \vv_r\in \widetilde{\calW}_{n_r}^R(\uu_r\tbold)\end{array}}.
	\end{align}
Also denote $|\rho| := n_l + |\uu_l| + |\tbold| + |\uu_r| + n_r$.
\end{definition}

\begin{remark}\label{remark:partition-of-calS_w}
	We note that we can partition $\calS_{\ww}$ as $\{\calS_{\ww}(\rho)\}_{|\rho|=|\ww|}$ where we run over all tuples $\rho = (n_l,n_r,\uu_l,\uu_r,\tbold)$ for $n_l,n_r\geq 0$, $(\uu_l,\uu_r,\tbold)\in \calT$ with  $|\rho| = |\ww|$.
	Indeed, if $(\ww_1,\ww_2,\ww_3)\in \calS_{\ww}(\rho)$ then $\ww_1 = \vv_l\uu_l$, $\ww_2 = \tbold$, $\ww_3 = \uu_r^{-1}\vv_r^{-1}$ for some $\vv_l \in \widetilde{\calW}_{n_l}^R(\uu_l\tbold)$ and $\vv_r \in \widetilde{\calW}_{n_r}^R(\uu_r\tbold)$ and we obtain that
	$\tbold=\ww_2$, $\uu_l = (\uu_l\tbold)\tbold =  \sbold_r(\vv_l\uu_l\tbold)\tbold= \sbold_r(\ww_1\ww_2)\ww_2$ and $\uu_r = (\uu_r\tbold)\tbold =  \sbold_r(\vv_r\uu_r\tbold)\tbold= \sbold_r(\ww_3^{-1}\ww_2)\ww_2$ and $n_l = |\ww_1|-|\uu_l| = |\ww_1| - |\sbold_r(\ww_1\ww_2)\ww_2|$ and $n_r = |\ww_3| - |\uu_r|= |\ww_3| - |\sbold_r(\ww_3^{-1}\ww_2)\ww_2|$. Since we can retrieve $\rho$ from $(\ww_1,\ww_2,\ww_3)$, this shows the sets $\calS_{\ww}(\rho)$ are disjoint.
	
	Now let $(\ww_1,\ww_2,\ww_3)\in \calS_{\ww}$ and set $\tbold := \ww_2$, $\uu_l := \sbold_r(\ww_1\tbold)\tbold$, $\uu_r := \sbold_r(\ww_3^{-1}\tbold)\tbold$.
	Then  $\uu_l\tbold$ and $\tbold\uu_r$ are clique words and
	\begin{align}
		|\ww| &\leq |\ww_1\ww_2\sbold_r(\ww_1\ww_2)| + |\sbold_r(\ww_1\ww_2)\ww_2\sbold_l(\ww_2\ww_3)| + |\sbold_l(\ww_2\ww_3)\ww_2\ww_3| \\
		&= (|\ww_1\ww_2| - |\sbold_r(\ww_1\ww_2)|) + |\uu_l\tbold\uu_r| + (|\ww_2\ww_3| - |\sbold_l(\ww_2\ww_3)|)\\
		&=|\ww|+|\uu_l\tbold\uu_r|  - |\sbold_r(\ww_1\ww_2)| +|\ww_2| - |\sbold_l(\ww_2\ww_3)|\\
		&=|\ww|+|\uu_l\tbold\uu_r|  - |\sbold_r(\ww_1\ww_2)\ww_2| -|\ww_2| - |\sbold_l(\ww_2\ww_3)\ww_2|\\
		&=|\ww|+|\uu_l\tbold\uu_r| - |\uu_l|  -|\tbold|  - |\uu_r|\\
		&\leq |\ww|.
	\end{align}
	Thus all inequalities must be equalities and we get $|\uu_l\tbold\uu_r| = |\uu_l|+|\tbold|+|\uu_r|$ so $\uu_l\tbold\uu_r$ is reduced. This shows $(\uu_l,\uu_r,\tbold)\in \calT$. 
	Now, set
	$n_l := |\ww_1| - |\uu_l|\geq 0$, $n_r := |\ww_3| - |\uu_r|\geq 0$.  Then we have $\vv_l:= \ww_1\uu_l^{-1}\in \widetilde{\calW}_{n_l}^R(\uu_l\tbold)$ and $\vv_r := \ww_3^{-1}\uu_r^{-1}\in \widetilde{\calW}_{n_r}^R(\uu_r\tbold)$. 
	Set $\rho=(n_l,n_r,\uu_l,\uu_r,\tbold)$ and observe that  $|\rho| = n_l + |\uu_l\tbold\uu_r| + n_r = |\ww_1|+|\ww_2|+|\ww_3| = |\ww|$.	
	Now, as $\ww_1 = \vv_l\uu_l$, $\ww_2 =\tbold$ and $\ww_3 = \uu_r^{-1}\vv_r^{-1}$ we obtain $(\ww_1,\ww_2,\ww_3)\in \calS_{\ww}(\rho)$. This proves the claim.
	
\end{remark}

	\begin{proposition}
		For $n_l,n_r\geq 0$ and $(\uu_l,\uu_r,\tbold)\in \calT$ define for the tuple $\rho=(n_l,n_r,\uu_l,\uu_r,\tbold)$ an operator $\widetilde{H}_{\rho}: \calB(\calF)\to \calB(\calF)$ as 
		\begin{align}
		\widetilde{H}_{\rho} =  \sum_{\rr\subseteq \tbold}(-1)^{|\rr|}H_{(n_l,n_r,\uu_l,\uu_r,\tbold,\rr)}.
		\end{align}
		Then we have for $\ww\in W$, $\omega\in \calS_{\ww}$ and $a\in \boldA$ that
		\begin{align}\label{eq:formula-for-H-rho}
		\widetilde{H}_{\rho}(\lambda_{\omega}(a)
		) = \begin{cases}
		\lambda_{\omega}(a) & \text{ if } \omega\in \calS_{\ww}(\rho)\\
		0 & \text{ else}
		\end{cases}.
		\end{align}
		
		\begin{proof}
			Let $\ww\in W$, $\omega \in \calS_{\ww}$ and let $a=a_1\otimes \cdots \otimes a_t\in \boldA$ be a pure tensor.
			By  \Cref{prop:formula-for-H-tau} we have
			\begin{align}
			\widetilde{H}_{\rho}(\lambda_{\omega}(a)) &= \sum_{\rr\subseteq \tbold}(-1)^{|\rr|}\lambda_{\omega}(a)P_{a}((\rho,\rr),\omega).
			\end{align}
			Let $\vv\in W$ and let $\eta\in\mathring{\calH}_{\vv}\subseteq\calF$ be a pure tensor. If $\lambda_{\omega}(a)\eta=0$, then it is clear that $\widetilde{H}_{\rho}(\lambda_\omega(a))\eta =0$, so that \Cref{eq:formula-for-H-rho} applied to $\eta$ holds in either case. Thus assume $\lambda_{\omega}(a)\eta\not=0$.	Let 
			$\calI_{\eta,\omega}$ be the set of all sub-clique words $\rr\subseteq \tbold$ s.t. $P_a((\rho,\rr),\omega)\eta=\eta$,
			that is 
			\begin{align}
			\calI_{\eta,\omega} = \{\rr\subseteq \tbold| P_a((\rho,\rr),\omega)\eta\not=0\}.
			\end{align}
			We prove the proposition using the following steps.\\
			
			1) We prove that $\calI_{\eta,\omega}$ is closed under taking sub-cliques. Let $\rr_1\subseteq \rr_2\subseteq \tbold$, and suppose that $\rr_2\in \calI_{\eta,\omega}$. Then 	we must have $\eta\in \calJ_{\lambda_\omega(a)}(\mathring{\calH}_{\vv_r\uu_r\tbold\vv_{tail}},\mathring{\calH}_{\vv_l\uu_l\tbold\vv_{tail}})$ with $\vv_l\in \calW_{n_l}^R(\uu_l\tbold)$, $\vv_r\in \calW_{n_r}^R(\uu_r\tbold)$ and $\vv_{tail}\in \calW^L(\uu_l\tbold)\cap \calW^L(\uu_r\tbold)$, and $|\vv_r\uu_r\tbold\rr_2|= |\ww_2\ww_3|+|\ww_2\ww_3\vv_r\uu_r\tbold\rr_2|$
			and $|\ww_1\ww_3\vv_r\uu_r\tbold\vv_{tail}| = |\ww_1|+|\ww_3\vv_r\uu_r\tbold\vv_{tail}|$
			
			Now this means that also 
			\begin{align}
			|\vv_r\uu_r\tbold| &\leq |\vv_r\uu_r\tbold\rr_1|+|\rr_1| \\
			&\leq |\ww_1\ww_2|+|\ww_1\ww_2\vv_r\uu_r\tbold\rr_1| +|\rr_1|\\
			&\leq |\ww_1\ww_2| + |\ww_1\ww_2\vv_r\uu_r\tbold\rr_2|+|\rr_2\rr_1|+|\rr_1|\\
			&= |\vv_r\uu_r\tbold\rr_2| + |\rr_2| \\
			&= |\vv_r\uu_r\tbold|
			\end{align}
			and therefore $|\vv_r\uu_r\tbold\rr_1| = |\ww_1\ww_2|+|\ww_1\ww_2\vv_r\uu_r\tbold\rr_1|$. This shows $P_a((\rho,\rr_1),\omega)\eta=\eta$, hence $\rr_1\in \calI_{\eta,\omega}$.\\
			
			2) We prove that $\calI_{\eta,\omega}$ is closed under taking unions. Let $\rr_1,\rr_2\subseteq \tbold$ be sub-cliques with $\rr_1,\rr_2\in \calI_{\eta,\omega}$.
			Then $P_a((\rho,\rr_1),\omega)\eta = P_a((\rho,\rr_2),\omega)\eta=\eta$. Moreover, by previous step we moreover have $P_{a}((\rho,e),\omega)\eta=\eta$. 
			We must now have $\eta\in \calJ_{\lambda_\omega(a)}(\mathring{\calH}_{\vv_r\uu_r\tbold\vv_{tail}},\mathring{\calH}_{\vv_l\uu_l\tbold\vv_{tail}})$ with $\vv_l\in \widetilde{\calW}_{n_l}^R(\uu_l\tbold)$, $\vv_r\in \widetilde{\calW}_{n_r}^R(\uu_r\tbold)$ and $\vv_{tail}\in \calW^L(\uu_l\tbold)\cap \calW^L(\uu_r\tbold)$, and 
			$|\ww_1\ww_3\vv_r\uu_r\tbold\vv_{tail}| = |\ww_1|+|\ww_3\vv_r\uu_r\tbold\vv_{tail}|$
			and moreover
			\begin{align}
			\label{eq:prop_Htilde:1} |\vv_r\uu_r\tbold|&= |\ww_2\ww_3|+|\ww_2\ww_3\vv_r\uu_r\tbold|\\
			|\vv_r\uu_r\tbold\rr_1|&= |\ww_2\ww_3|+|\ww_2\ww_3\vv_r\uu_r\tbold\rr_1|\\
			|\vv_r\uu_r\tbold\rr_2|&= |\ww_2\ww_3|+|\ww_2\ww_3\vv_r\uu_r\tbold\rr_2|.
			\end{align}
			Now we note that also $|\vv_r\uu_r\tbold|= |\vv_r\uu_r\tbold\rr_1|+|\rr_1| = |\vv_r\uu_r\tbold\rr_2| + |\rr_2|$, hence
			\begin{align}
				|\ww_2\ww_3\vv_r\uu_r\tbold| = |\ww_2\ww_3\vv_r\uu_r\tbold\rr_1| +|\rr_1| = |\ww_2\ww_3\vv_r\uu_r\tbold\rr_2| + |\rr_2|.
			\end{align}
			As $\rr_1,\rr_2$ are cliques, this implies $\rr_1,\rr_2\subseteq \sbold_r(\ww_2\ww_3\vv_r\uu_r\tbold)$ so that for $\rr=\rr_1\cup \rr_2$ it holds that $\rr\subseteq \sbold_r(\ww_2\ww_3\vv_r\uu_r\tbold)$.
			But this implies 
			\begin{align}
			|\ww_2\ww_3\vv_r\uu_r\tbold| = |\ww_2\ww_3\vv_r\uu_r\tbold\rr| +|\rr|.
			\end{align}
			Now, as also $|\vv_r\uu_r\tbold| = |\vv_r\uu_r\tbold\rr|+|\rr|$ we find using \eqref{eq:prop_Htilde:1} that 
			$|\vv_r\uu_r\tbold\rr|=|\ww_2\ww_3|+|\ww_2\ww_3\vv_r\uu_r\tbold\rr|$. It now directly follows that $P((\rho,\rr),\omega)\eta = \eta$. This shows that $\rr\in \calI_{\eta,\omega}$, and thus that $\calI_{\eta,\omega}$ is closed under taking unions.\\
			
			3) We prove the equation $\widetilde{H}_{\rho}(\lambda_{\omega}(a))\eta = \mathds{1}(\calI_{\eta,\omega}=\{e\})\lambda_{\omega}(a)\eta$.
			Here $\mathds{1}(\calI_{\eta,\omega}=\{e\})$ denotes $1$ whenever $\calI_{\eta,\omega}=\{e\}$ is satisfied, and $0$ otherwise.
			In the case that $\calI_{\eta,\omega}$ is empty we directly find $	\widetilde{H}_{\rho}(\lambda_{\omega}(a))\eta=0$, so that the equation is satisfied. Thus assume that $\calI_{\eta,\omega}$ is non-zero. Then as $\calI_{\eta,\omega}$ is closed under taking unions, there exists a maximal element $\rr_{\eta,\omega}\in \calI_{\eta,\omega}$. However, since $\calI_{\eta,\omega}$ is also closed under taking sub-cliques, we then find
			$\calI_{\eta,\omega} = \{\rr\subseteq \rr_{\eta,\omega}\}$.
			We conclude that 
			\begin{align}
			\widetilde{H}_{\rho}(\lambda_{\omega}(a))\eta &= \sum_{\rr\subseteq \tbold}(-1)^{|\rr|}\lambda_{\omega}(a)P_{a}((\rho,\rr),\omega)\eta\\
			&=\sum_{\rr\subseteq \rr_{\eta,\omega}}(-1)^{|\rr|}\lambda_{\omega}(a)\eta\\
			&=\mathds{1}(\rr_{\eta,\omega}=e)\lambda_{\omega}(a)\eta\\
			&=\mathds{1}(\calI_{\eta,\omega}=\{e\})\lambda_{\omega}(a)\eta.
			\end{align}
			
			4) We will now show, for a pure tensor $\eta\in\mathring{\calH}_{\vv} \subseteq \calF$ with $\lambda_\omega(a)\eta \not=0$, that $\calI_{\eta,\omega}=\{e\}$ if and only if  $\omega\in \calS_{\ww}(\rho)$.
			First, suppose that $\omega\in \calS_{\ww}(\rho)$. Then we can write $\omega = (\ww_1,\ww_2,\ww_3)$, where $\ww_1 = \vv_l\uu_l$ and $\ww_2= \tbold$ and $\ww_3 = \uu_r^{-1}\vv_r^{-1}$ for some $\vv_l\in \widetilde{\calW}_{n_l}^R(\uu_l\tbold)$ and $\vv_r\in \widetilde{\calW}_{n_r}^R(\uu_r\tbold)$.
			Then as $\lambda_{\omega}(a)\eta\not=0$, we must have that $\eta\in \calJ_{\lambda_{\omega}(a)}({\mathring{\calH}_{\vv_r\uu_r\tbold\vv_{tail}},\mathring{\calH}_{\vv_l\uu_l\tbold\vv_{tail}}})$ for some $\vv_{tail}\in \calW^L(\uu_l\tbold)\cap \calW^L(\uu_r\tbold)$. It is clear that 
			\begin{align}
			|\ww_1\ww_3\vv_r\uu_r\tbold\vv_{tail}| &= |\vv_l\uu_l\tbold\vv_{tail}| \\
			&=|\vv_l\uu_l| + |\tbold\vv_{tail}|\\
			&= |\ww_1|+|\ww_3\vv_r\uu_r\tbold\vv_{tail}|.
			\end{align}
			Moreover, as $\ww_2\ww_3\vv_r\uu_r\tbold\subseteq \tbold$ it is also clear that $|\vv_r\uu_r\tbold| = |\ww_2\ww_3|+|\ww_2\ww_3\vv_r\uu_r\tbold|$. This shows that $P_a((\rho,e),\omega)\eta = \eta$, hence $e\in \calI_{\eta,\omega}$.
			
			Now let  $\rr\subseteq \tbold$ be a sub-clique with $\rr\not=e$. Then we have $\ww_2\ww_3\vv_r\uu_r\tbold\rr=\rr$. Hence, we have
			\begin{align}
			|\vv_r\uu_r\tbold\rr|+|\rr| &= |\vv_r\uu_r\tbold| \\
			&= |\ww_2\ww_3|+|\ww_2\ww_3\vv_r\uu_r\tbold|\\ 
			&=|\ww_2\ww_3|+|\ww_2\ww_3\vv_r\uu_r\tbold\rr| - |\rr|.
			\end{align}
			Now as $\rr\not=e$ we have $|\rr|\geq 1$, which shows that 
			$|\vv_r\uu_r\tbold\rr|\not=|\ww_2\ww_3|+|\ww_2\ww_3\vv_r\uu_r\tbold\rr|$. This proves that $P_a((\rho,\rr),\omega)\eta=0$. Thus $\rr\not\in \calI_{\eta,\omega}$. This shows $\calI_{\eta,\omega}=\{e\}$.
		
			Now, let $\omega\in \calS_{\ww}$ for some $\ww\in W$ be s.t. $\calI_{\eta,\omega}=\{e\}$.
			Then $P((\rho,e),\omega)\eta=\eta$. Hence
			$\eta \in \calJ_{\lambda_{\omega}(a)}(\mathring{\calH}_{\vv_r\uu_r\tbold\vv_{tail}},\mathring{\calH}_{\vv_l\uu_l\tbold\vv_{tail}})$ for some  $\vv_l\in \widetilde{\calW}^R(\uu_l\tbold)$, $\vv_r\in \widetilde{\calW}^R(\uu_r\tbold)$ and $\vv_{tail} \in \calW^L(\uu_l\tbold)\cap \calW^L(\uu_r\tbold)$ and 
			$|\ww_1\ww_3\vv_r\uu_r\tbold\vv_{tail}|=|\ww_1|+|\ww_3\vv_r\uu_r\tbold\vv_{tail}|$ and $|\vv_r\uu_r\tbold| = |\ww_2\ww_3| + |\ww_2\ww_3\vv_r\uu_r\tbold|$. Now as also $\lambda_{\omega}(a)\eta\in \mathring{\calH}_{\ww_1\ww_3\vv_r\uu_r\tbold\vv_{tail}}$, and as $\lambda_{\omega}(a)\eta\not=0$, we have that $\ww_1\ww_3\vv_r\uu_r\tbold\vv_{tail} = \vv_l\uu_l\tbold\vv_{tail}$. Hence, $\ww_1\ww_3 = \vv_l\uu_l\uu_r^{-1}\vv_{r}^{-1}$. Now, as $P_a((\rho,\rr),\omega)\eta=0$ for all $\rr\subseteq \tbold$ with $\rr\not=e$, we must have that $\sbold_r(\ww_2\ww_3\vv_r\uu_r\tbold)\cap \tbold=e$. However, multiplying $\ww_2\ww_3$ with $\vv_r\uu_r\tbold$ removes all letters from $\ww_2\ww_3$. This means that $\sbold_r(\ww_2\ww_3\vv_r\uu_r\tbold)\subseteq \sbold_r(\vv_r\uu_r\tbold)=\sbold_r(\uu_r\tbold)$.
			Now we also have
			\begin{align}
				|\vv_l\uu_l\tbold|&\leq |\ww_2\ww_1^ {-1}| + |\ww_2\ww_1^{-1}\vv_l\uu_l\tbold|\\
				&= |\ww_2\ww_1^ {-1}|+|\vv_r\uu_r\tbold| - |\ww_2\ww_3|\\
				&\leq |\ww_2\ww_1^ {-1}| + |\ww_3\vv_r\uu_r\tbold|-|\ww_2|\\
				&=|\ww_1| + |\ww_3\vv_r\uu_r\tbold|\\
				&=|\ww_1\ww_3\vv_r\uu_r\tbold|\\
				&=|\vv_l\uu_l\tbold|
			\end{align}
			so that $|\vv_l\uu_l\tbold|=|\ww_2\ww_1^ {-1}| + |\ww_2\ww_1^{-1}\vv_l\uu_l\tbold|$. 
			Now this means that $\sbold_r(\ww_2\ww_1^{-1}\vv_l\uu_l\tbold)\subseteq \sbold_r(\vv_l\uu_l\tbold)=\sbold_r(\uu_l\tbold)$.
			Hence, as $\ww_2\ww_3\vv_r\uu_r\tbold = \ww_2\ww_1^{-1}\vv_l\uu_l\tbold$, we find
			$\sbold_r(\ww_2\ww_3\vv_r\uu_r\tbold)\subseteq \uu_l\tbold\cap \uu_r\tbold = \tbold$.
			However, as also $\sbold_r(\ww_2\ww_3\vv_r\uu_r\tbold)\cap \tbold = e$, we conclude that $\sbold_r(\ww_2\ww_3\vv_r\uu_r\tbold)=e$, that is $\ww_2\ww_1^{-1}\vv_l\uu_l\tbold = \ww_2\ww_3\vv_r\uu_r\tbold = e$.
			But this means that $\ww_3^{-1}\ww_2 = \vv_r\uu_r\tbold$
			and  $\ww_1\ww_2 = \vv_l\uu_l\tbold$. From this it follows that $\ww_2 \subseteq \sbold_r(\vv_l\uu_l\tbold) \cap \sbold_r(\vv_r\uu_r\tbold) = \tbold$. Now, we can not have that $\ww_2\subseteq \tbold$ strictly, as this would mean that $\ww_3$ starts with a part of $\tbold$ that $\ww_1$ ends with, which would contradict the fact that $\ww_1\ww_2\ww_3$ is reduced. Thus we now find $\ww_2 = \tbold$ and then also $\ww_1 = \vv_l\uu_l$ and $\ww_3 = \uu_r^{-1}\vv_r^{-1}$.
			This means that $\omega\in \calS_{\ww}(\rho)$.\\
			
			5) We now conclude the proof of the proposition as we have shown for $\ww\in W$, $\omega\in \calS_{\ww}$, pure tensor $a=a_1\otimes \cdots\otimes a_t\in \boldA$ and pure tensor $\eta\in \mathring{\calH}_{\vv}\subseteq \calF$ with $\lambda_{\omega}(a)\eta\not=0$ that
			\begin{align}
				\widetilde{H}_{\rho}(\lambda_{\omega}(a))\eta
				&= \begin{cases}
				\lambda_{\omega}(a)\eta & \calI_{\eta,\omega}=\{e\}\\
				0 & else
				\end{cases}
				&= \begin{cases}
				\lambda_{\omega}(a)\eta & \omega\in \calS_{\ww}(\rho)\\
				0 & else
				\end{cases}.
			\end{align}
			Now, as noted earlier, the equation is also satisfied when $\eta$ is a pure tensor with $\lambda_{\omega}(a)\eta =0$.
			Therefore, by linearity and continuity, the equation in the proposition holds for all $\eta\in \calF$. By linearity of $\widetilde{H}_{\rho}$ and $\lambda_{\omega}$ the equation also holds for all $a\in \boldA$. This proves the statement.
		\end{proof}
	\end{proposition}

We now prove our main theorem of this section, that shows that $\|\calP_{d}\|_{cb}$ is polynomially bounded in $d$.
	\begin{theorem}\label{thm:projection-maps-linear-growth}
	For $d\geq 0$ we have  (on $\calA$) that		
	\begin{align}
	\calP_{d} &= \sum_{\substack{(\uu_l,\uu_r,\tbold)\in \calT\\
			0\leq n\leq d-|\uu_l\tbold\uu_r|}}
	\sum_{\rr\subseteq \tbold}(-1)^{|\rr|} H_{(n,d-n-|\uu_l\tbold\uu_r|,\uu_l,\uu_r,\tbold,\rr)}.
	\end{align}
	Moreover, for $d\geq 1$ we get the linear bound $\|\calP_d\|_{cb} \leq C_{\Gamma} d$, where $C_{\Gamma}$ denotes the constant
	\begin{align} C_{\Gamma} = \sum_{(\uu_l,\uu_r,\tbold)\in \calT}2^{|\tbold|}.
	\end{align} 
	\begin{proof}
		For $d\geq 0$ define
		\begin{align}
		\calT_{d} &= \{\rho = (n_l,n_r,\uu_l,\uu_r,\tbold)\in \ZZ_{\geq 0}^2\times \calT| |\rho|=d\}.
		\end{align}
		We recall for $\ww\in W$ that $\{\calS_{\ww}(\rho)\}_{\rho\in \calT_{|\ww|}}$ is a partition of $\calS_\ww$ by \cref{remark:partition-of-calS_w}.
		Fix some $a \in \boldA$. For $d\geq 0$ we find using \Cref{lemma:partition-action} that
		\begin{align}
		\calP_{d}(\lambda(a))&=  
		\sum_{\ww\in W, |\ww|=d}\sum_{\omega\in \calS_{\ww}} \lambda_{\omega}(a)\\
		&= 	\sum_{\rho\in \calT_d}
	 	\sum_{\ww\in W, |\ww|=d}\sum_{\omega\in \calS_{\ww}(\rho)} \lambda_{\omega}(a)\\
		&= 	\sum_{\rho\in \calT_d}\widetilde{H}_{\rho}\left(
		\sum_{\ww\in W}\sum_{\omega\in \calS_{\ww}
		} \lambda_{\omega}(a)\right)\\
		&=\sum_{\rho\in \calT_d}\widetilde{H}_{\rho}(
		\lambda(a))\\
		&=\sum_{\substack{(\uu_l,\uu_r,\tbold)\in \calT\\
		0\leq n\leq d-|\uu_l\tbold\uu_r|}}\sum_{\rr\subseteq \tbold}(-1)^{|\rr|}H_{(n,d-|\uu_l\tbold\uu_r|-n,\uu_l,\uu_r,\tbold,\rr)}(\lambda(a)).
		\end{align}
		Therefore, the equation holds on $\lambda(\boldA)$ and hence, by continuity, on $\calA$.
		
		Now let $d\geq 1$, we show that the bound holds. We note first that by definition $V_{n}^{e,e}=0$ for $n\geq 1$. This implies directly that $H_{(n,d-n-|\uu_l\tbold\uu_r|,\uu_l,\uu_r,\tbold,e)} = 0$ for $0\leq n\leq d-|\uu_l\tbold\uu_r|$ whenever $(\uu_l,\uu_r,\tbold)=(e,e,e)$. Therefore we find
		\begin{align}
		\|\calP_{d}\|_{cb} &\leq  \sum_{\substack{(\uu_l,\uu_r,\tbold)\in \calT\setminus\{(e,e,e)\}\\
				0\leq n\leq d-|\uu_l\tbold\uu_r|}}
		\sum_{\rr\subseteq \tbold} \|H_{(n,d-n-|\uu_l\tbold\uu_r|,\uu_l,\uu_r,\tbold,\rr)}\|_{cb}\\
		&\leq \sum_{\substack{(\uu_l,\uu_r,\tbold)\in \calT\setminus\{(e,e,e)\}\\
				0\leq n\leq d-|\uu_l\tbold\uu_r|}}
		2^{|\tbold|}\\
		&\leq \left(\sum_{\substack{(\uu_l,\uu_r,\tbold)\in \calT}}
		2^{|\tbold|}\right)d.
		\end{align}		
	\end{proof}
\end{theorem}

\section{Graph products of state-preserving u.c.p maps}\label{section:graph-products-of-ucp-maps}
In \cref{subsection:graph-products-of-ucp-maps} we show that the graph product of state-preserving u.c.p maps extends to a state-preserving u.c.p map.
Thereafter, in \cref{subsection:CCAP-for-finite-dimensional-algebras}, we use this to obtain the result that the graph product of finite-dimensional algebras with GNS-faithful states is weakly amenable with constant $1$. 

\subsection{Graph products of state-preserving ucp maps}\label{subsection:graph-products-of-ucp-maps}
Let $\Gamma$ be a graph, and for $v\in V\Gamma$ let $\theta_v:\boldA_{v}\to \boldB_{v}$ be state-preserving maps between unital $\Cstar$-algebras (with states s.t. the GNS representation is faithful). Let $(\calA,\varphi)=*_{v,\Gamma}(\boldA_{v},\varphi_v)$ and $(\calB,\psi)= *_{v,\Gamma}(\boldB_{v},\psi_v)$ be their reduced graph products. As $\theta_v$ is state preserving it maps $\mathring{\boldA}_v$ to $\mathring{\boldB}_v$. We can look at the map  $\theta:\lambda(\boldA)\to \lambda(\boldB)$ for $a_1\otimes \cdots\otimes a_s \in \mathring{\boldA}_{v_1}\otimes \cdots \otimes \mathring{\boldA}_{v_s}$ for a reduced word $v_1\cdots v_s$ given as
\begin{align}
	\theta(\lambda(a_1\otimes \cdots\otimes a_s)) = 
	\lambda(\theta_{v_1}(a_1)\otimes \cdots \otimes \theta_{v_s}(a_s))
\end{align}
and we set $\theta(\Id)=\Id$. We denote this map by $\theta=*_{v,\Gamma}\theta_v$ and call it the graph product map. The map is clearly state-preserving.  
To prove the main theorem, we need the result that the graph product map $\theta = *_{v,\Gamma}\theta_v$ of state-preserving u.c.p. maps $\theta_v$ extends to a bounded map on the graph product, and that it is again u.c.p. This result was already proven by Blanchard-Dykema in \cite{blanchardEmbeddingsReducedFree2001} for the case of free products.  For graph products the result has been proven by Caspers-Fima in \cite[Proposition 3.30]{caspersGraphProductsOperator2017a} in the setting of von Neumann algebras.

\begin{proposition}\cite[Proposition 3.30]{caspersGraphProductsOperator2017a}
	\label{prop:graph-product-of-ucp-maps-vNa}
	Let $\Gamma$ be a simple graph and for $v\in V\Gamma$, let $\theta_v:\boldM_{v}\to \boldN_{v}$ be state-preserving normal u.c.p. maps between von Neumann algebras $\boldM_v$ and $\boldN_v$ that have faithful normal states. Let $(\calM,\varphi) = \overline{*_{v,\Gamma}} (\boldM_{v},\varphi_v)$ and $(\calN,\psi)=\overline{*_{v,\Gamma}}(\boldN_{v},\psi_v)$ be the  von Neumann algebraic graph products. Then there exists a unique normal u.c.p. map $\theta: \calM\to \calN$ s.t. for all pure tensors $a_1\otimes \cdots \otimes a_s\in \boldM_{v_1}\otimes \cdots \otimes \boldM_{v_s}$ we have
	\begin{align}
		\theta(\lambda(a_1\otimes \cdots \otimes a_s)) = \lambda(\theta_{v_1}(a_1)\otimes \cdots \otimes \theta_{v_s}(a_s)).
	\end{align}
	The map $\theta$ will be denoted as $\theta = *_{\Gamma}\theta_v$
	
\end{proposition}

We give here a proof for the case of $\Cstar$-algebras.
\begin{proposition}\label{prop:graph-products-ucp-maps-Cstar}
	For $v\in V\Gamma$  Let $\theta_v:\boldA_v\to \boldB_v$ be state-preserving, unital completely positive maps between unital $\Cstar$-algebras $(\boldA_v,\varphi_v)$ and $(\boldB_v,\psi_v)$, and assume $\varphi_v$ and $\psi_v$ are GNS-faithful. Then the graph product map $\theta = *_{v,\Gamma}\theta_v$ extends to a state-preserving unital completely positive map between the reduced graph products $\calA$ and $\calB$.
	\begin{proof}
		We will use the notation $\calH_{\vv}^{\calA}$, $\mathring{\calH}_{\vv}^{\calA}$, $\calF^{\calA}$, $\lambda^{\calA}$, $\Omega^{\calA}$, et cetera,  corresponding to the reduced graph product $(\calA,\varphi):= *_{v,\Gamma}(\boldA_{v},\varphi_v)$, and use similar notation for the reduced graph product $(\calB,\psi):= *_{v,\Gamma}(\boldB_{v},\psi_v)$.
		By the Stinespring’s dilation theorem we can write $\theta_v(a) = V_v^*\pi_v(a)V_v$ for some Hilbert space $\widehat{\calH}_v$ and unital $*$-homomorphism $\pi_v:\boldA_v\to \calB(\widehat{\calH}_v)$ of $\boldA_v$ and some isometry $V_v\in \calB(\calH_{v}^{\calB},\widehat{\calH}_v)$.
		We note that for $a\in \boldA_{v}$ we have $\varphi_v(a) = \psi_v(\theta_v(a)) = \langle \theta_v(a)\xi_v^{\calB},\xi_v^{\calB}\rangle = \langle \pi_v(a)\widehat{\xi}_v,\widehat{\xi}_v\rangle$ with $\widehat{\xi}_v = V_v\xi_v^{\calB}$. Also  $\pi_v$ is faithful, as $\pi_{v}(a) = 0$ implies for $b\in \boldA_{v}$ that $0 = \|\pi_v(a)\pi_v(b)\widehat{\xi}_v\|^2 = \|\pi_v(ab)\widehat{\xi}_v\|^2 = \langle \pi_v(b^*a^*ab)\widehat{\xi}_v,\widehat{\xi}_v\rangle = \varphi_v(b^*a^*ab)$, which implies $a =0$ since $\varphi_v$ is GNS-faithful. By these properties we conclude that we can construct the graph product of the $\boldA_{v}$'s w.r.t. the representations $\pi_v$. To distinguish the notation from the other graph products we use \textit{hat}-notation like $\widehat{\calH}_{\vv}$, $\mathring{\widehat{\calH}}_{\vv}$, $\widehat{\calF}$,$\widehat{\lambda}$, $\widehat{\Omega}$. Define a contraction  $V:\calF^{\calB}\to \widehat{\calF}$ for $\eta = \eta_1\otimes \cdots \otimes \eta_l \in \mathring{\calH}_{\vv}^{\calB}$ as  
		\begin{align}
			V|_{\mathring{\calH}_{\vv}}(\eta_1\otimes \cdots \otimes \eta_l) = V_{v_1}\eta_1\otimes \cdots \otimes V_{v_l}\eta_l
		\end{align}
		and $V(\Omega^{\calB}) = \widehat{\Omega}$.
		We note that $\eta_i \in \mathring{\calH}_{v_i}^{\calB}$ implies $\langle V\eta_i,\widehat{\xi}_{v_i}\rangle = \langle V\eta_i,V\xi_{v_i}^{\calB}\rangle = \langle \eta_i,\xi_{v_i}^{\calB}\rangle = 0$ and hence $V\eta_i\in \mathring{\widehat{\calH}}_{v_i}$. This shows that $V$ is well-defined.\\
		
		By \cite[Proposition 3.12]{caspersGraphProductsOperator2017a}, we know that there is a state-preserving, unital $*$-homomorphism $\pi:\calA\to \calB(\widehat{\calF})$ that for $a = a_1\otimes \cdots\otimes a_l\in \mathring{\boldA}_{\vv}$ is given by
		\begin{align}
			\pi(\lambda^{\calA}(a_1\otimes \cdots \otimes a_l)) = \widehat{\lambda}(\pi_{v_1}(a_1)\otimes \cdots \otimes \pi_{v_l}(a_l))
		\end{align}		
		We will now show that $\theta(\lambda^{\calA}(a)) = V^*\pi(\lambda^{\calA}(a))V$ for $a\in \boldA$, which then shows that $\theta$ can be extended to a u.c.p. map on $\calA$.

		Let $\eta = \eta_1\otimes \cdots \otimes \eta_l\in \mathring{\calH}_{\vv}^{\calB}$ for some $\vv\in W$ and let $a\in \mathring{\boldA}_{v}$ for some $v\in V\Gamma$. 
		We will calculate $\widehat{\lambda}_{v}(\pi_v(a))V$.
		First suppose that $v\vv$ is reduced. We have $\langle (I-V_vV_v^*)\pi_v(a)\widehat{\xi}_v,\widehat{\xi}_v\rangle  =\langle \pi_v(a)\widehat{\xi}_v,0\rangle = 0$ so that
		\begin{align}
			\widehat{\lambda}_v((I-V_vV_v^*)\pi_v(a))V\eta &= \widehat{U}_v((I-V_vV_v^*)\pi_v(a)\otimes \Id_{\widehat{\calF}})(\widehat{\xi}_v\otimes V\eta)\\
			&=\widehat{U}_v((I-V_vV_v^*)\pi_v(a)\widehat{\xi}_v\otimes V\eta)\\
			&=\widehat{\calQ}_{(v,\vv)}(((I-V_vV_v^*)\pi_v(a)\widehat{\xi}_v\otimes V\eta).
		\end{align}
		Also we have $\langle V_vV_v^*\pi_v(a)\widehat{\xi}_v,\widehat{\xi}_v\rangle = \varphi_v(a)=0$ and so we find 
		\begin{align}
			\widehat{\lambda}_v(V_vV_v^*\pi_v(a))V\eta &= \widehat{U}_v(V_vV_v^*\pi_v(a)\otimes \Id_{\widehat{\calF}})(\widehat{\xi}_v\otimes V\eta)\\
			&=\widehat{U}_v(V_vV_v^*\pi_v(a)\widehat{\xi}_v\otimes V\eta)\\
			&=\widehat{\calQ}_{(v,\vv)}((V_vV_v^*\pi_v(a)\widehat{\xi}_v)\otimes V\eta)\\			&=\widehat{\calQ}_{(v,\vv)}((V_v\theta_v(a)\xi_v^{\calB})\otimes V\eta)\\
			&=V\calQ_{(v,\vv)}^{\calB}((\theta_v(a)\xi_v^{\calB})\otimes \eta)\\
			&=V\lambda_{v}^{\calB}(\theta_v(a))\eta.
		\end{align}
		
		Now, on the other hand suppose that $\vv$ starts with $v$. Then we can write $\eta = \calQ_{(v,v\vv)}^{\calB}(\eta_0\otimes \eta')$ for some $\eta_0\in \mathring{\calH}_{v}^{\calB}$ and $\eta'\in \mathring{\calH}_{v\vv}^{\calB}$ and we have $V\eta = \widehat{\calQ}_{(v,v\vv)}(V_{v}\eta_0\otimes V\eta')$. Again $\langle (I-V_vV_v^*)\pi_v(a)V_v\eta_0,\widehat{\xi}_v\rangle = 0$ and so
		\begin{align}
			\widehat{\lambda}_{v}((I-V_vV_v^*)\pi_v(a))V\eta 
			&=\widehat{U}_{v}((I-V_vV_v^*)\pi_v(a)\otimes \Id_{\calF})\widehat{U}_v^*V\eta\\
			&=\widehat{U}_{v}((I-V_vV_v^*)\pi_v(a)\otimes \Id_{\calF})(V\eta_0\otimes V\eta')\\
			&=\widehat{U}_{v}\left(((I-V_vV_v^*)\pi_v(a)V_{v}\eta_0)\otimes V\eta'\right)\\
			&=\widehat{\calQ}_{(v,v\vv)}\left(((I-V_vV_v^*)\pi_v(a)V_{v}\eta_0)\otimes V\eta'\right).
		\end{align}
		Furthermore, we have
		\begin{align}
			\widehat{\lambda}_{v}(V_vV_v^*\pi_v(a))V\eta 
			&=\widehat{U}_{v}(V_vV_v^*\pi_v(a)\otimes \Id_{\calF})\widehat{U}_v^*V\eta\\
			&=\widehat{U}_{v}(V_vV_v^*\pi_v(a)\otimes \Id_{\calF})(V\eta_0\otimes V\eta')\\
			&=\widehat{U}_{v}\left(V_vV_v^*\pi_v(a)V_{v}\eta_0)\otimes V\eta'\right)\\
			&=\widehat{U}_{v}\left((V_v\theta_v(a)\eta_0)\otimes V\eta'\right)\\
			&=VU_{v}^{\calB}\left((\theta_v(a)\eta_0)\otimes \eta'\right)\\
			&=VU_v^{\calB}\left(\theta_v(a)\otimes \Id_{\calF}\right)(U_{u}^{\calB})^*\eta\\
			&=V\lambda^{\calB}(\theta_v(a))\eta.
		\end{align}
		
		Now, when $a = a_1\otimes\cdots \otimes a_k\in \mathring{\boldA}_{\ww}$, then we have
		\begin{align}
			V^*\pi(\lambda^{\calA}(a))&V\eta =  V^*\widehat{\lambda}(\pi_{w_1}(a_1))\ldots \widehat{\lambda}(\pi_{w_{k-1}}(a_{k-1}))\widehat{\lambda}(V_{w_k}V_{w_k}^*\pi_{w_k}(a_{k}))V\eta \\
			&+V^*\widehat{\lambda}(\pi_{w_1}(a_1))\ldots \widehat{\lambda}(\pi_{w_{k-1}}(a_{k-1}))\widehat{\lambda}((I-V_{w_k}V_{w_k}^*)\pi_{w_k}(a_{k}))V\eta \\
			&= V^*\widehat{\lambda}(\pi_{w_1}(a_1))\ldots \widehat{\lambda}(\pi_{w_{k-1}}(a_{k-1}))\widehat{\lambda}(V_{w_k}V_{w_k}^*\pi_{w_k}(a_{k}))V\eta \\
			&=V^*\pi(\lambda^{\calA}(a_1\otimes \cdots \otimes a_{k-1}))V\lambda^{\calB}(\theta_{w_k}(a_{k}))\eta.
		\end{align}
		Note here that the reason why we can remove the second summand is because one tensor leg of  $\widehat{\lambda}((I-V_{w_k}V_{w_k}^*)\pi_{w_k}(a_{k}))V\eta$ is of the form $(I-V_{w_k}V_{w_k}^*)\pi_{w_k}(a_k)V_{w_k}\eta_0$ for some $\eta_0\in \calH_{w_k}^{\calB}$. This tensor leg is not changed by the operator $\pi(\lambda^{\calA}(a_1\otimes \cdots \otimes a_{k-1}))$ as it may not act on the same letter. Now after the application of $V^*$ we obtain for this tensor leg that $V_{w_k}^*(I-V_{w_k}V_{w_k}^*)\pi_{w_k}(a_k)V_{w_k}\eta_0 =0$, so that this term vanishes.
		
		By what we showed, it now follows from induction to the tensor length $k$ that $V^*\pi(\lambda^{\calA}(a))V = \theta(\lambda^{\calA}(a))$ for all $a\in \boldA$. This then shows the statement.
	\end{proof}

\end{proposition}

\subsection{CCAP for reduced graph products of finite dimensional algebras}	
\label{subsection:CCAP-for-finite-dimensional-algebras}
We now state the following generalization of \cite[Proposition 3.5.]{ricardKhintchineTypeInequalities2006a} to graph products. The proof uses \Cref{thm:projection-maps-linear-growth} and  \Cref{prop:graph-product-of-ucp-maps-vNa} and \Cref{prop:graph-products-ucp-maps-Cstar} and goes analogously to \cite[Proposition 3.5.]{ricardKhintchineTypeInequalities2006a}.
\begin{proposition}\label{prop:semi-group-ucp-maps}
	Let $\Gamma$ be a finite simple graph. For $v\in V\Gamma$ let $\boldA_v$ be a unital $\Cstar$-algebra together with a GNS-faithful state $\varphi_v$. Let $(\calA,\varphi):= *_{v,\Gamma}(\boldA_v,\varphi_v)$ be the reduced graph product. For $d\geq 0$ let $\calP_d:\calA\to \calA_d$ be the natural projection. Let $0\leq r\leq 1$, $n\in\NN$ and define
	\begin{align}
		\calT_r = \sum_{k=0}^{\infty}r^ k\calP_{k} &\quad \calT_{r,n} = \sum_{k=0}^n r^k\calP_k.
	\end{align}
	Then $\calT_r$ and $\calT_{r,n}$ are completely bounded with 
	\begin{align}\label{eq:cb-estimate}
		\|\calT_r\|_{cb}\leq 1 &\quad \text{ and } \quad \|\calT_{r}-\calT_{r,n}\|_{cb}\leq \frac{C_{\Gamma}nr^ n}{(1-r)^ 2}.
	\end{align}
	The maps $\calT_{e^{-t}}$ for $t\geq 0$ form a one-parameter semi-group of unital completely positive maps on $\calA$ preserving the state $\varphi$.
	Moreover, the sequence $(\calT_{1-\frac{1}{\sqrt{n}},n})_{n\geq 1}$ tends pointwise to the identity of $\calA$ and $\lim_{n\to\infty}\|\calT_{1-\frac{1}{\sqrt{n}},n}\|_{cb} = 1$.
	\begin{proof}
		
	For $v\in V\Gamma$ we define a state-preserving u.c.p map $U_{r,v}:\boldA_{v}\to \boldA_{v}$ as $U_{r,v}(a) = ra + (1-r)\varphi_v(a)\Id_{\calH_{v}}$. It can be seen that $*_{v,\Gamma} U_{r,v} = \calT_r$ on $\lambda(\boldA)$ and by \cref{prop:graph-products-ucp-maps-Cstar} this map extends to a state-preserving u.c.p map on $\calA$. Thus $\|\calT_r\|_{cb}= 1$. Furthermore,
	\begin{align}
		\|\calT_{r} -\calT_{r,n}\|_{cb}\leq \sum_{k=n+1}^\infty r^k \|\calP_k\|_{cb}
		\leq C_{\Gamma}\sum_{k=n}^\infty kr^{k}=  C_{\Gamma}r\frac{d}{dr}\left(\frac{r^n}{1-r}\right)
	\end{align}
	Therefore, as $\frac{d}{dr}\left(\frac{r^n}{1-r}\right) = n r^{n-1}(1-r)^{-1} + r^n(1-r)^{-2} \leq  n r^{n-1}(1-r)^{-2}$ this proves \eqref{eq:cb-estimate}.  It is furthermore clear that $(\calT_{e^{-t}})_{t\geq 0}$ forms a semi-group since $\calP_{m}\calP_n = 0$ when $n\not=m$.
 	By \eqref{eq:cb-estimate} and by the triangle inequality we have
	$\|\calT_{1-\frac{1}{\sqrt{n}},n}\|_{cb}\leq 1+C_{\Gamma}n^2(1-\frac{1}{\sqrt{n}})^n\to 1$ as $n\to\infty$ which shows $\lim\limits_{n\to\infty}\|\calT_{1-\frac{1}{\sqrt{n}},n}\|_{cb}=1$ since the maps $\calT_{1-\frac{1}{\sqrt{n}},n}$ are unital. Moreover, on $\lambda(\boldA)$ it is clear that $(\calT_{1-\frac{1}{\sqrt{n}},n})_{n\geq 1}$ tends pointwise to the identity. Therefore, as $(\calT_{1-\frac{1}{\sqrt{n}},n})_{n\geq 1}$ is uniformly bounded it follows by density that this holds true on $\calA$ as well.

	\end{proof}
	
\end{proposition}

\begin{corollary}\label{corollary:CCAP-for-finite-dim-algebras}
	For $v\in V\Gamma$ let $\boldA_{v}$ be a finite-dimensional $\Cstar$-algebras together with a GNS-faithful state $\varphi_v$. Then the reduced graph product $\calA$ has the CCAP. Similarly, for finite dimensional von Neumann algebras $\boldM_{v}$ together with normal faithful states $\varphi_v$, we have that the graph product $\calM$ has the wk-$*$ CCAP.
\end{corollary}

We give an application of this result to Hecke-algebras (for references on Hecke-algebras see \cite[Chapter 19]{davisGeometryTopologyCoxeter2008}). Let $W$ be a Coxeter group generated by some set $S$ and let $q= (q_s)_{s\in S}$ be a Hecke tuple (i.e. $q_s>0$ for all $s\in S$ and $q_s = q_t$ whenever $s$ and $t$ are conjugate in $W$). We denote $\calN_{q}(W)$ for the Hecke algebra corresponding to $W$ and $q$ . Our application uses the following proposition which asserts that we can decompose Hecke algebras as graph products. This result for right-angled Coxeter groups is stated in \cite[Corollary 3.4]{caspersAbsenceCartanSubalgebras2020}.
	
\begin{proposition}\label{prop:graph-product-decomposition-Hecke-algebra}
	Let $\Gamma$ be a graph, and for $v\in V\Gamma$ let $W_{v}$ be a Coxeter group generated by a set $S_{v}$ and let $q_v=(q_{v,s})_{s\in S_{v}}$ be a Hecke-tuple. Set $W = *_{v,\Gamma} W_{v}$ and $q:= *_{v,\Gamma}q_{v} = (q_{v,s})_{v\in V\Gamma, s\in S_v}$. Then we get a graph product decomposition $\calN_{q}(W) = \overline{*_{v,\Gamma}}\calN_{q_v}(W_v)$.
	\begin{proof}
		This follows from \cite[Proposition 3.22]{caspersGraphProductsOperator2017a} by considering the natural embeddings $\pi_{v}:\calN_{q_v}(W_v)\to \calN_{q}(W)$ that send generators to generators.
	\end{proof}
\end{proposition}  The following was already known from \cite[Theorem A]{caspersAbsenceCartanSubalgebras2020}, but we believe our approach is more conceptual.
\begin{example} Let $W$ be a right-angled Coxeter group generated by a finite set $S$, and $q=(q_v)_{v\in S}$ a Hecke-tuple. Then as $W = *_{v,\Gamma}(\ZZ/2\ZZ)$ for some (finite) graph $\Gamma$,  we can by \cref{prop:graph-product-decomposition-Hecke-algebra} write $\calN_q(W) = \overline{*_{v,\Gamma}}\calN_{q_v}(\ZZ/2\ZZ)$. As $\calN_{q_v}(\ZZ/2\ZZ)$ is finite dimensional we obtain by \cref{corollary:CCAP-for-finite-dim-algebras} that $\calN_q(W)$ has the wk-$*$ CCAP.
\end{example}
The result for the following example is new.
\begin{example} Let $\Gamma$ be a finite simple graph, and for $v\in V\Gamma$ let $W_{v}$ be a finite Coxeter group generated by some set $S_v$ and let $q_{v} = (q_{v,s})_{s\in S_v}$ be a Hecke-tuple for $W_v$. Then if we let $W = *_{v, \Gamma}W_v$ and $q = *_{v,\Gamma} q_v := (q_{v,s})_{v\in V\Gamma, s\in S_{v}}$, we have by \cref{prop:graph-product-decomposition-Hecke-algebra} that $\calN_{q}(W) = \overline{*_{v,\Gamma}}\calN_{q_v}(W_v)$. Since $\calN_{q_v}(W_{v})$ is finite dimensional we obtain by \cref{corollary:CCAP-for-finite-dim-algebras} that $\calN_{q}(W)$ possesses the wk-$*$ CCAP.
\end{example}

\section{Graph product of completely bounded maps on $\calA_d$}
\label{section:graph-product-of-cb-maps-on-Ad}
The main result of this section is \cref{thm:graph-product-cb-maps}, which shows that the graph product of completely bounded maps $T_v$ defines a completely bounded map $T_d$ on the homogeneous subspace $\calA_d$ of degree $d$. The proof of this results follows the lines of \cite{ricardKhintchineTypeInequalities2006a} (where they use the different convention $\langle \hat{a},\hat{b}\rangle = \varphi(a^*b)$), and uses the construction of the operator space $X_d$ as in \cite{caspersGraphProductKhintchine2021a} and another operator space $\widetilde{X}_d$, to extend it to graph products.

\subsection{Free products and operator spaces}

When given a finite graph $\Gamma$ and  algebras $(\boldA_{v},\varphi_v)$ we will denote the reduced free product of the algebras as $(\calA^{f},\varphi^f) = *_{v} (\boldA_v,\varphi_v)$. Let $\Gamma^f$ be the graph with vertex set $V\Gamma^f = V\Gamma$ and no edges. Note that the free product is simply the reduced graph product corresponding to $\Gamma^f$. For the graph product corresponding to $\Gamma^f$ we will use notation using superscript $f$, that is we will write $W^f$, $\lambda^f$, $P_{v}^f$, $\calF^f$,$\mathring{\calH}_{\vv}^f$, $\mathring{\boldA}_{\ww}^f$, et cetera. We remark that
$\calF\subseteq \calF^f$ and $\boldA\subseteq \boldA^f$ as linear subspaces and that $\boldA_{v}=\boldA_{v}^f$ for $v\in V\Gamma$.
For $\ww\in W\setminus\{e\}$ with representative $(w_1,\ldots, w_n)$ we will define $\calH_{\ww} = \calH_{w_1}\otimes \cdots \calH_{w_n}$ and $\boldA_{\ww} = \boldA_{w_1}\otimes \cdots \otimes \boldA_{w_n}$, and we define $\calH_{e} = \CC\Omega$ and $\boldA_{e} = B(\calH_e)$. 
Define  a subspace $L_1$ of $\calB(\calF^f)$ by the closed linear span
\begin{align}
	L_1 = \overline{\Span}\{P_v^f\lambda_v^f(a)P_v^{f\perp}| v\in V\Gamma, a\in \mathring{\boldA}_{v}^f\}, \quad K_1 = L_1^*.
\end{align}
For a Hilbert space $\calH$ denote $\calH_{C}$, $\calH_{R}$ respectively for the column and row Hilbert space, see \cite{pisierIntroductionOperatorSpace2003}. In \cite[Lemma 2.3 and Corollary 2.4]{ricardKhintchineTypeInequalities2006a} it is shown that 
\begin{align}\label{eq:L1K1-complete-isomorphism}
	L_1 \simeq (\oplus_{v\in V\Gamma}\mathring{\calH}_v)_{C}, \quad K_1 \simeq (\oplus_{v\in V\Gamma}\mathring{\calH}_v^{op})_R
\end{align}
completely isometrically, and that the maps $\theta_1:\calA_1^f\to L_1$ and $\rho_1:\calA_1^f\to K_1$ given for $a\in \mathring{\boldA}_v$ by
$\theta_1(\lambda_v^f(a)) = P_{v}^f\lambda_v^f(a)P_{v}^{f \perp}$ and $\rho_1(\lambda_v^f(a))=P_{v}^{f \perp}\lambda_{v}^f(a)P_{v}^f$ are completely contractive. We denote $\otimes_h$ for the Haagerup tensor product, see \cite[Chapter 9]{effrosOperatorSpaces2000}. We denote $L_d = L_1^{\otimes_h d}$ and $K_d = K_1^{\otimes_h d}$ for the $d$-fold tensor product and we write $\theta_1^{\otimes d}$ for the map $\calA_{d}^f\to L_{d}$ defined for $b = b_1\otimes \cdots \otimes b_{d}\in \boldA_{d}$ by $\theta_1^{\otimes d}(\lambda^{f}(b)) = \theta_1(\lambda^{f}(b_1))\otimes_h \cdots \otimes_h \theta_1(\lambda^{f}(b_d))$  and we write  $\rho_1^{\otimes d}$ for the map $\calA_d^{f}\to K_d$ defined similarly.

We introduce notation similar to \cite[ Section 2]{caspersGraphProductKhintchine2021a}.
Let $\ww\in W^f$ s.t. in the graph product $\ww$ is equivalent to some clique word $\vv_{\Gamma_0}$ for some clique $\Gamma_0\subseteq \Gamma$ (which we will denote by $\ww \equiv \vv_{\Gamma_0}$). Let $a=a_1\otimes \cdots\otimes a_d\in \boldA_{\ww}^f$.
We define an operator $\Diag_{\ww}(a):\calF^f\to \calF^f$ on $\mathring{\calH}_{\vv}^f$ for $\vv\in W^f$ with $|\vv|= |\ww|+|\ww^{-1}\vv|$ as
\begin{align}
	\Diag_{\ww}(a)|_{\mathring{\calH}_{\vv}^f} = P_{v_1}a_1P_{v_1}\otimes \cdots \otimes P_{v_d}a_dP_{v_d}\otimes \Id_{\mathring{\calH}_{v_{d+1}}}\otimes \cdots \otimes \Id_{\mathring{\calH}_{v_{|\vv|}}}
\end{align}
and we define $\Diag_{\ww}(a)|_{\mathring{\calH}_{\vv}^f} = 0$ if $\vv\in W^f$ is not of the given form. 
Extending this, we obtain a linear map $\Diag_{\ww}:\boldA_{\ww}^f\to \calB(\calF^f)$. For a clique $\Gamma_0$ in $\Gamma$, we now define the operator space $\boldA_{\Gamma_0} = \Span\{\Diag_{\ww}(\boldA_{\ww}^f)| \ww\in W^f, \ww\equiv \vv_{\Gamma_0}\}$. Also, for $\ww\in W^f$ we consider $\boldA_{\ww}^f$ as an operator space by the embedding $\boldA_{\ww}^f\subseteq \calB(\calH_{\ww}^f)$.

\begin{proposition}
	For a clique $\Gamma_0$ and a word $\ww\in W^f$ with $\ww\equiv \vv_{\Gamma_0}$ we have that the map $\Diag_{\ww}:\boldA_{\ww}^f\to \boldA_{\Gamma_0}$ is completely contractive.
	\begin{proof}
		We define a map $V_{\ww}:\calF^f\to \calH_{\ww}^f\otimes \calF^f$ as 
		\begin{align}
			V_{\ww}|_{\mathring{\calH}_{\vv}^f} :=\calQ_{(\ww,\ww^{-1}\vv)}^{f*}
		\end{align}
		whenever $\vv\in W^f$ is s.t. $|\vv|=|\ww|+|\ww^{-1}\vv|$  and set $V_{\ww}|_{\mathring{\calH}_{\vv}^f} = 0$ when $\vv$ is not of this form.  We then obtain that
		\begin{align}
			\Diag_{\ww}(a) = V_{\ww}^*(a\otimes \Id_{\calF})V_{\ww}
		\end{align}
		which shows the statement.
				
	\end{proof}
	
\end{proposition}

As in \cite{ricardKhintchineTypeInequalities2006a} and \cite{caspersGraphProductKhintchine2021a} we define operator spaces $X_d$ and additionally we will define other operator spaces $\widetilde{X}_d$.  For $\tbold\in W$ a clique word, denote $\Gamma_{\tbold}$ for the clique in $\Gamma$. We now set
\begin{align}
	X_d&=\bigoplus_{\substack{n_l,n_r\geq 0,\\ (\uu_l,\uu_r,\tbold)\in \calT\\ n_l+|\uu_l\tbold\uu_r|+n_r = d}}L_{n_l+|\uu_l|}\otimes_h \boldA_{\Gamma_{\tbold}} \otimes_h K_{n_r+|\uu_r|}\\
	\widetilde{X}_d&=\bigoplus_{\substack{n_l,n_r\geq 0,\\ (\uu_l,\uu_r,\tbold)\in \calT\\ n_l+|\uu_l\tbold\uu_r|+n_r = d}}L_{n_l+|\uu_l|}\otimes_h \boldA_{\tbold} \otimes_h K_{n_r+|\uu_r|}
	\end{align} 
equipped with the sup-norm. We remark here that the operator space structure on $\boldA_{\tbold}$ is given by the inclusion $\boldA_\tbold = \boldA_{\tbold'}^f\subseteq \calB(\calH_{\tbold'}^f)$ where $\tbold'\in W^f$ is the representant of $\tbold$.
Also, recall that $\calT$ was defined in \cref{definition:clique-pairs}  and that in \cref{def:sets-Swrho}  for a tuple $\rho = (n_l,n_r,\uu_l,\uu_r,\tbold)$ with $n_l,n_r\geq 0$, $(\uu_l,\uu_r,\tbold)\in \calT$ we defined $|\rho| = n_l + |\uu_l| + |\tbold| +|\uu_r| +n_r$. 
By the above, we can find a completely contractive map $D_d:\widetilde{X}_d\to X_d$ by defining $D_d = (D_{\rho})_{\rho, |\rho|=d}$ where  $D_{\rho} = (\Id_{L_{n_l+|\uu_l|}} \otimes \Diag_{\tbold'}\otimes \Id_{K_{n_r+|\uu_r|}})$ for $\rho=(n_l,n_r,\uu_l,\uu_r,\tbold)$.

We now define two linear maps $\widetilde{\Theta}_d:\calA_{d}\to \widetilde{X}_d$ and 
$j_d:\boldA_d\to X_d$ as follows. Fix a tuple  $\rho = (n_l,n_r,\uu_l,\uu_r,\tbold)$, $|\rho|=d$. We denote $\widetilde{n}_l = n_l + |\uu_l|$ and $\widetilde{n}_r = n_r + |\uu_r|$.
Let $a\in \mathring{\boldA}_{\ww}$ be a pure tensor with $\ww\in W$.
Suppose that $\ww = \vv_l\uu_l\tbold\uu_r^{-1}\vv_r^{-1}$ for some $\vv_l\in \widetilde{\calW}_{n_l}^R(\uu_{l}\tbold)$ and $\vv_r\in \widetilde{\calW}_{n_r}^R(\uu_{r}\tbold)$. We can then write $a = \calQ_{(\vv_l\uu_l,\tbold,\uu_r^{-1}\vv_r^{-1})}(a_1\otimes a_2\otimes a_3)$ for some $a_1\in \mathring{\boldA}_{\vv_l\uu_l}$, $a_2\in \mathring{\boldA}_{\tbold}$ and $a_3\in \mathring{\boldA}_{\uu_r^{-1}\vv_r^{-1}}$.
We then defined
\begin{align}
	\widetilde{\Theta}_d(\lambda(a))_{\rho} &= \theta_{1}^{\otimes\widetilde{n_l}}(\lambda^f(a_1))\otimes a_2\otimes \rho_1^{\otimes\widetilde{n_r}}(\lambda^f(a_3))\\
	j_d(a)_{\rho} &= \theta_{1}^{\otimes\widetilde{n_l}}(\lambda^f(a_1))\otimes \Diag_{\tbold'}(a_2)\otimes \rho_1^{\otimes\widetilde{n_r}}(\lambda^f(a_3)).
\end{align}
In the case that $\ww$ is not of the given form we define $\widetilde{\Theta}_d(\lambda(a))_{\rho}=0$ and $j_d(a)_{\rho}=0$. This is extended linearly and we set $\widetilde{\Theta}_d(\lambda(a))=(\widetilde{\Theta}_d(\lambda(a))_{\rho})_{\rho}$ and $j_d(a)=(j_d(a)_{\rho})_{\rho}$. We moreover define the map $\Theta_d:=D_d\circ \widetilde{\Theta}_d$ and see that  $j_d = \Theta_d \circ \lambda|_{\boldA_d}$.  We note that the definition of $j_d$  agrees with that in \cite[Equation (2.16)]{caspersGraphProductKhintchine2021a}, and that, in the case of dealing with free products, the map $\Theta_d$ compares with a restriction of the map $\Theta_d$ in \cite{ricardKhintchineTypeInequalities2006a}. In \cite[Equation (2.24)]{caspersGraphProductKhintchine2021a} a completely bounded map $\pi_d:E_d\to \calB(\calF)$ was defined, where $E_d := j_d(\boldA_d)\subseteq X_d$,  and that satisfied $\pi_d\circ j_d = \lambda|_{\boldA_d}$. For $d\geq 1$ the norm bound $\|\pi_d\|_{cb} \leq (\#\Cliq(\Gamma))^3d$ holds by \cite[Theorem 2.9]{caspersGraphProductKhintchine2021a}, where $\#\Cliq(\Gamma)$ denotes the number of cliques in the graph $\Gamma$. We get the following commuting diagram:

\begin{figure}[h!]
	
	\begin{tikzpicture}[baseline]
		\node at (-1.5, -1)			(Xdt) {$\widetilde{X}_d$};
		\node at (-1.5, 0) 	  	(Ed) {$E_d$};
		\node at (-3, 2)   		(Ad) {$\boldA_{d}$};
		\node at (-3, 0)   	(Ad2) {$\calA_d$};
		\node at (-1.5,2)   		(Adf) {$\boldA_{d}^f$};
		\node at (0,0)  	 	(Xd) {$X_d$};
		\node at (-2.3,2)  	 	(inclusion1) {$\subseteq$};
		\node at (-0.8,0)  	 	(inclusion2) {$\subseteq$};
		
		\draw[<-] (Xdt) -- node[xshift =-0.2cm, yshift=-0.2cm] {$\widetilde{\Theta}_d$} (Ad2);
		\draw[->] (Xdt) -- node[xshift =0.2cm, yshift=-0.2cm] {$D_d$} (Xd);

		\draw[->] (Ad) -- node[right] {$\lambda$} (Ad2);
		\draw[<-] (Ad2) -- node[above] {$\pi_d$} (Ed);
		
		\draw[->] (Ad) -- node[xshift =0.2cm, yshift=0.2cm] {$j_d$} (Ed);
	
	\end{tikzpicture}
\end{figure}

For a clique word $\tbold \in W$ with representative $(t_1,\ldots, t_{|\tbold|})$ we define a  unitary $U: \calH_{\tbold}\to \bigoplus_{\rr\subseteq \tbold}\mathring{\calH}_{\rr}$ in a natural way. Let $\eta =\eta_1\otimes \cdots \eta_{|\tbold|}\in \calH_{\tbold}$ be a tensor with either $\eta_i \in \mathring{\calH}_{t_i}$ or $\eta_i\in \CC\xi_{t_i}$. For $1\leq i\leq |\tbold|$ set $r_i:=t_i$ when $\eta_i\in \mathring{\calH}_{t_i}$ and $r_i = e$ when $\eta_i\in \CC\xi_{t_i}$. Then $\rr:=r_1\cdots r_{|\tbold|}$ is a subword of $\tbold$ since $\tbold$ is a clique word. Using the identification $\CC\xi_{t_i} \simeq \mathring{\calH}_{e}$ given by $\xi_{t_i}\to \Omega$ we can define $U(\eta) = \calQ_{(r_1,\ldots,r_{|\tbold|})}(\eta)\in \mathring{\calH}_{\rr}$. This extends linearly to a unitary. We remark that for $a\in \mathring{\boldA}_{\tbold}$ we have $U^*\lambda(a)U = a$. Indeed, it can be checked that for $a_i\in \mathring{\boldA}_{t_i}$ we have $U^*\lambda(a_i)U = \Id_{\calH_{t_1}}\otimes \cdots \Id_{{\calH}_{t_{i-1}}} \otimes a_i \otimes \Id_{\calH_{t_{i+1}}}\otimes \cdots \otimes \Id_{\calH_{t_{|\tbold|}}}$ so that the statement follows as $\lambda(a_1\otimes \cdots \otimes a_n) = \lambda(a_1)\cdots \lambda(a_n)$.

\begin{theorem}\label{thm:theta-d-completely-contractive}
	The map $\widetilde{\Theta}_d$ is completely contractive.
	
	\begin{proof}
		Choose $d\geq 0$. Fix a tuple $\rho = (n_l,n_r,\uu_l,\uu_r,\tbold)$ with $|\rho|=d$ and write $\widetilde{n_l} = n_l+|\uu_l|$, $\widetilde{n_r} = n_r + |\uu_r|$.
		We define two partial isometries 
		\begin{align}
			J_{\rho}^{L}&:{\calF^f}^{\otimes \widetilde{n_l}}\otimes \calH_{\tbold} \to  {\calF^f}^{\otimes \widetilde{n_l}}\otimes \calF\\ 
			J_{\rho}^{R}&:\calH_{\tbold}\otimes {\calF^f}^{\otimes \widetilde{n_r}} \to \calF\otimes {\calF^f}^{\otimes \widetilde{n_r}}
		\end{align}as follows.
		Let $\rr_l\subseteq \tbold$, let $\eta = \eta_1\otimes \cdots  \otimes \eta_{\widetilde{n_l}}\otimes \eta_0\in {\calF^f}^{\otimes \widetilde{n_l}}\otimes (U^*\mathring{\calH}_{\rr_l})$ be a pure tensor and denote $\eta_0':=U\eta_0 \in \mathring{\calH}_{\rr_l}$.  If for $i\geq 1$ we can write $\eta_i = \eta_i' \otimes \widetilde{\eta}_i$ for some $\eta_i' \in \mathring{\calH}_{v_i}$ and $\widetilde{\eta}_i \in \calF^f$ for which $(v_1,\ldots ,v_{n})$ is the representative of $\vv_l\uu_l$ for some $\vv_l\in \widetilde{\calW}_{n_l}^R(\uu_l\tbold)$, then we define
		\begin{align}
			J_{\rho}^{L}\eta &= \widetilde{\eta}_1 \otimes \cdots\otimes \widetilde{\eta}_{\widetilde{n_l}}\otimes \calQ_{(v_1,\ldots,v_{\widetilde{n_l}},\rr_l)}(\eta_1' \otimes \cdots\otimes \eta_{\widetilde{n_l}}'\otimes \eta_0')\in {\calF^{f}}^{\otimes\widetilde{n_l}}\otimes \mathring{\calH}_{\vv_l\uu_l\rr_l}
		\end{align}
		and we define $J_{\rho}^L$ as $0$ on the complement of all such tensors.
		Similarly, let $\rr_r\subseteq \tbold$ let $\eta = \eta_0 \otimes \eta_1\otimes \cdots \otimes \eta_{\widetilde{n_r}}\in (U^*\mathring{\calH}_{\rr_r})\otimes {\calF^f}^{\otimes \widetilde{n_r}}$, denote $\eta_0' := U\eta_0\in \mathring{\calH}_{\rr_r}$ and suppose that for $i\geq 1$ we can write $\eta_i = \eta_i' \otimes \widetilde{\eta}_i$ for some $\eta_i' \in \mathring{\calH}_{v_i}$ and $\widetilde{\eta}_i \in \calF^f$ for which  $(v_1,\ldots,v_{n})$ is the representative of $\uu_r^{-1}\vv_r^{-1}$ for some $\vv_r\in \widetilde{\calW}_{n_r}^R(\uu_r\tbold)$ we define
		\begin{align}
			J_{\rho}^{R}\eta &= \calQ_{(v_{\widetilde{n_r}},\ldots,v_1,\rr_r)}(\eta_{\widetilde{n_r}}' \otimes \cdots\otimes \eta_1'\otimes \eta_0') \otimes \widetilde{\eta}_1 \otimes .. \otimes \widetilde{\eta}_{\widetilde{n_r}}\in \mathring{\calH}_{\vv_r\uu_r\rr_r}\otimes \calF^{f\otimes \widetilde{n_r}}
		\end{align} and we define $J_{\rho}^R$ as $0$ on the complement of all such tensors.
	
		We shall show that 
		\begin{align}\label{eq:dilation_Theta_d}
			\widetilde{\Theta}_d(\lambda(a))_{\rho} = (J_{\rho}^{L*}\otimes \Id_{\calF^f}^{\otimes \widetilde{n_r}})(\Id_{\calF^f}^{\otimes \widetilde{n_l}}\otimes \lambda(a)\otimes \Id_{\calF^f}^{\otimes \widetilde{n_r}})(\Id_{\calF^f}^{\otimes \widetilde{n_l}}\otimes J_{\rho}^{R})
		\end{align}
	which then shows the statement.

		Let $\ww\in W$, $|\ww|=d$, let  $a\in \mathring{\boldA}_{\ww}$ be a pure tensor, let $\omega = (\ww_1,\ww_2,\ww_3)\in \calS_{\ww}$, $\vv_l\in \widetilde{\calW}_{n_l}^R(\uu_l\tbold)$, $\vv_r\in \widetilde{\calW}_{n_r}^R(\uu_r\tbold)$ and $\rr_l,\rr_r\subseteq \tbold$. Now let $\eta\in \mathring{\calH}_{\vv_r\uu_r\rr_r}$ be a pure tensor, in which case $\lambda_{\omega}(a)\eta$ is also a pure tensor. Suppose that $\lambda_{\omega}(a)\eta \in \mathring{\calH}_{\vv_l\uu_l\rr_l}$ and that it is non-zero, so  that $\vv_l\uu_l\rr_r$ and $\vv_r\uu_r\uu_r$ start with $\ww_1\ww_2$ and $\ww_3^{-1}\ww_2$ respectively and so that $\ww_1\ww_3\vv_r\uu_r\rr_r = \vv_l\uu_l\rr_l$. Then put $\ww_{tail}:= \ww_2\ww_3\vv_r\uu_r\rr_r =
		\ww_2\ww_1^{-1}\vv_l\uu_l\rr_l$ so that 
	$\ww_1\ww_2\ww_{tail}$ and $\ww_3^{-1}\ww_2\ww_{tail}$ are reduced expressions for $\vv_l\uu_l\rr_l$ and $\vv_r\uu_r\rr_r$ respectively. We claim that $\sbold_r(\ww_2\ww_{tail})\supseteq \sbold_r(\ww_1\ww_2\ww_{tail})\cap \sbold_r(\ww_3^{-1}\ww_2\ww_{tail})$. Indeed, let $v$ be a letter in $\sbold_r(\ww_1\ww_2\ww_{tail})$ that is not in $\sbold_r(\ww_2\ww_{tail})$. Then $v$ is a letter at the end of $\ww_1$ that commutes with $\ww_2$. If $v$ is at the same time a letter in $\sbold_r(\ww_3^{-1}\ww_2\ww_{tail})$ then $v$ is also a letter at the end of $\ww_3^{-1}$, i.e. a letter at the start of $\ww_3$. But this would contradict the fact that $\ww_1\ww_2\ww_3$ is reduced. Thus we established the inclusion and obtain  $$\sbold_r(\ww_2\ww_{tail})\supseteq \sbold_r(\ww_1\ww_2\ww_{tail})\cap \sbold_r(\ww_3^{-1}\ww_2\ww_{tail}) =  \sbold_r(\vv_l\uu_l\rr_l)\cap \sbold_r(\vv_r\uu_r\rr_r) \supseteq \rr_l\cap \rr_r$$ so that $|\ww_2\ww_{tail}|\geq |\rr_l\cap \rr_r|.$
	Now, combining all this, we find
	\begin{align}
		d + |\rr_l\cap \rr_r| + |\ww_{tail}|&\leq
		|\ww_1\ww_2\ww_3| + |\ww_2\ww_{tail}|+|\ww_{tail}|\\ 
		&=|\ww_1| + 2|\ww_2| + 2|\ww_{tail}| + |\ww_3|\\
		&=|\ww_1\ww_2\ww_{tail}| +|\ww_3^{-1}\ww_2\ww_{tail}|\\
		&= |\vv_l|+|\uu_l| + |\rr_l| + |\rr_r| +|\uu_r| + |\vv_r|\\
		&= d + |\rr_l|+|\rr_r| - |\tbold|\\
		&\leq d + |\rr_l| +|\rr_r| - |\rr_l\cup \rr_r| \\
		&= d +|\rr_l\cap \rr_r| 
	\end{align}
	We conclude that all the above inequalities must be equalities, in particular  $|\ww_{tail}| = 0$, $|\tbold| = |\rr_l\cup \rr_r|$ and $|\ww_2\ww_{tail}| = |\rr_l\cap \rr_r|$. This means $\tbold = \rr_l\cup \rr_r$ and $\ww_2 = \ww_2\ww_{tail}$. Now as also $\ww_2 = \ww_2\ww_{tail} \supseteq \rr_l\cap \rr_r$ we conclude that $\ww_2 = \rr_l \cap \rr_r$. Set $\tbold_l := \rr_l\ww_2 = (\rr_l\cap \tbold\rr_r)$, $\tbold_m := \rr_l\cap \rr_r$ and $\tbold_r := \ww_2\rr_r = (\tbold\rr_l\cap \rr_r)$. Then, as we know $\vv_l\uu_l\rr_l = \ww_1\ww_2\ww_{tail} = \ww_1\ww_2$ and $\vv_r\uu_r\rr_r = \ww_3^{-1}\ww_2\ww_{tail} = \ww_3^{-1}\ww_2$, we then obtain that $\ww_1 = \vv_l\uu_l\tbold_l$ and $\ww_3 = \tbold_r\uu_r^{-1}\vv_r^{-1}$. Hence, $\omega$ is of the form $\omega = (\vv_l\uu_l\tbold_l, \tbold_m, \tbold_r\uu_r^{-1}\vv_r^{-1})$. We note that $\tbold_l,\tbold_m, \tbold_r$ are disjoint subcliques of $\tbold$ with
	$\tbold_l\tbold_m\tbold_r = \tbold$. In particular we find that the assumption implies  $\ww = \vv_l\uu_l\tbold\uu_r^{-1}\vv_r^{-1}$. For a closed subspace $\calK\subseteq \calF$ denote $P_{\calK}$ for the orthogonal projection on $\calK$.	We conclude that
	\begin{align}
		P_{\mathring{\calH}_{\vv_l\uu_l\rr_l}}\lambda(a)	P_{\mathring{\calH}_{\vv_r\uu_r\rr_r}}
		= \lambda_{(\vv_l\uu_l\tbold_l,\tbold_m,\tbold_r\uu_r^{-1}\vv_r^{-1})}(a)	P_{\mathring{\calH}_{\vv_r\uu_r\rr_r}}		
	\end{align} and moreover that this expression is zero whenever $a\not\in \mathring{\boldA}_{\vv_l\uu_l\tbold\uu_r^{-1}\vv_r^{-1}}$.
This shows that for $a\in \mathring{\boldA}_{\ww}$ with $\ww$ not of the form $\ww=\vv_l\uu_l\tbold\uu_r^{-1}\vv_r^{-1}$ for any  $\vv_l\in \widetilde{\calW}_{n_l}^R(\uu_l\tbold)$, $\vv_r \in \widetilde{\calW}_{n_l}^R(\uu_r\tbold)$, the right-hand side of \eqref{eq:dilation_Theta_d} is zero. In this case also the left-hand side is zero by definition of $\Theta_d(\lambda(a))_{\rho}$ so that we get equality.

Let $\vv\in W$. We define.
\begin{align}
	\calK_{\rho,\vv}^L&= \bigoplus_{\rr_l\subseteq \tbold}\mathring{\calH}_{\vv\uu_l\rr_l} & 	\calK_{\rho,\vv} ^R&= {\bigoplus_{\rr_r\subseteq \tbold}\mathring{\calH}_{\vv\uu_r\rr_r}} \\
	\calK_{\rho}^L &= \bigoplus_{\vv_l\in \widetilde{\calW}_{n_l}^R(\uu_l\tbold)}\calK_{\rho,\vv_l} ^L & \calK_{\rho}^R &= \bigoplus_{\vv_r\in \widetilde{\calW}_{n_r}^R(\uu_r\tbold)}\calK_{\rho,\vv_r} ^R.
\end{align}

Let us now assume $a\in \mathring{\boldA}_{\ww}$ with $\ww = \vv_l\uu_l\tbold\uu_r^{-1}\vv_r^{-1}$ for some $\vv_l \in \widetilde{\calW}_{n_l}^R(\uu_l\tbold)$, $\vv_r \in \widetilde{\calW}_{n_r}^R(\uu_r\tbold)$ and write  $a = \calQ_{(\vv_l\uu_l,\tbold,\uu_r^{-1}\vv_r^{-1})}(a_1\otimes a_2\otimes a_3)$ for some $a_1\in \mathring{\boldA}_{\vv_l\uu_l}$, $a_2\in \mathring{\boldA}_{\tbold}$ and $a_3\in \mathring{\boldA}_{\uu_r^{-1}\vv_r^{-1}}$. Note that in such case the words $\vv_l,\vv_r$ are uniquely determined. By the above, we now find
\begin{align}
	&P_{\calK_{\rho}^L}\lambda(a)P_{\calK_{\rho}^R} =\\
	&= P_{\calK_{\rho,\vv_l}^L}\lambda(a)P_{\calK_{\rho,\vv_r}^R}\\	&=P_{\calK_{\rho,\vv_l}^L}\left(\sum_{\substack{\tbold_l,\tbold_m,\tbold_r\\ \text{partition of } \tbold}}\lambda_{(\vv_l\uu_l\tbold_l,\tbold_m,\tbold_r\uu_r^{-1}\vv_r^{-1})}(a)\right)P_{\calK_{\rho, \vv_r}^R}\\
	&=P_{\calK_{\rho,\vv_l}^L}\left(\sum_{\substack{\tbold_l,\tbold_m,\tbold_r\\ \text{partition of } \tbold}}\lambda_{(\vv_l\uu_l,e,e)}(a_1)\lambda_{(\tbold_l,\tbold_m,\tbold_r)}(a_2)\lambda_{(e,e,\uu_r^{-1}\vv_r^{-1})}(a_3)\right)P_{\calK_{\rho, \vv_r}^R}\\
	&\stackrel{\cref{lemma:partition-action}}{=}P_{\calK_{\rho,\vv_l}^L}\lambda_{(\vv_l\uu_l,e,e)}(a_1)\lambda(a_2)\lambda_{(e,e,\uu_r^{-1}\vv_r^{-1})}(a_3)P_{\calK_{\rho, \vv_r}^R}\\
	&=P_{\calK_{\rho,\vv_l}^L}\lambda_{(\vv_l\uu_l,e,e)}(a_1)(Ua_2U^*)\lambda_{(e,e,\uu_r^{-1}\vv_r^{-1})}(a_3)P_{\calK_{\rho, \vv_r}^R}
\end{align}
where we use that $\lambda(a_2)|_{\mathring{\calH}_{\rr}} = Ua_2U^*$ for $\rr\subseteq \tbold$.
Now, a calculation shows that
\begin{align}\label{eq:calculation1}
	(U^*\lambda_{(e,e,\uu_r^{-1}\vv_r^{-1})}(a_3)P_{\calK_{\rho,\vv_r}^R} \otimes \Id)J_{\rho}^R &= (\Id_{\calH_{\tbold}}\otimes \rho_1^{\otimes\widetilde{n_r}}(\lambda^f(a_3)))\\
	\label{eq:calculation2} J_{\rho}^{L*}(\Id\otimes P_{\calK_{\rho,\vv_l}^L}\lambda_{(\vv_l\uu_l,e,e)}(a_1)U) &= (\theta_1^{\otimes\widetilde{n_l}}(\lambda^f(a_1))\otimes \Id_{\calH_{\tbold}})
\end{align}
We describe the calculation for \eqref{eq:calculation1} (the calculation for \eqref{eq:calculation2} is similar by taking adjoints and using that $\theta_1^{\otimes \widetilde{n_l}}(\lambda^f(a_1))^* = \rho_1^{\otimes \widetilde{n_l}}(\lambda^f(a_1^*))$).
Let $\eta =\eta_0 \otimes \eta_1\otimes \cdots \eta_{\widetilde{n_r}} \in (U^*\mathring{\calH_{\rr_r}}) \otimes {\calF^f}^{\otimes \widetilde{n_r}}$ for some $\rr_r\subseteq \tbold$ and so that $\eta_i$ is a pure tensor for $i=0,\ldots, \widetilde{n_r}$. Assume that for $i=1,\ldots, \widetilde{n_r}$ we can write $\eta_i = \eta_{i}' \otimes \widetilde{\eta_{i}}$ with $\eta_i'\in \mathring{\calH}_{v_i}$ and $\widetilde{\eta}_i \in \calF^f$ for which $(v_1,\ldots, v_{\widetilde{n_r}})$ is the representative of $\uu_r^{-1}\vv_r^{-1}$. Indeed, if $\eta$ is not of this form then both $(P_{\calK_{\rho,\vv_r}^R}\otimes \Id)J_{\rho}^R\eta =0$ and $(\Id_{\calH_\tbold}\otimes \rho_1^{\otimes \widetilde{n_r}}(\lambda^f(a_3)))\eta = 0$  which gives the equality.
Now by definition $J_{\rho}^R\eta = \zeta_1\otimes \zeta_2$ where $\zeta_1:=\calQ_{(v_{\widetilde{n_r}},\ldots, v_{1},\rr_r)}(\eta_{\widetilde{n_r}}'\otimes \cdots \otimes \eta_{1}'\otimes U\eta_0)\in \mathring{\calH}_{\vv_r\uu_r\rr_r}$ and $\zeta_2 := \widetilde{\eta_1}\otimes \cdots \otimes \widetilde{\eta}_{\widetilde{n_r}}$. Now
\begin{align}
	(\lambda_{(e,e,\uu_r^{-1}\vv_r^{-1})}&(a_3)P_{\calK_{\rho,\vv_r}^R}\otimes \Id)J_{\rho}^R\eta =(\lambda_{(e,e,\uu_r^{-1}\vv_r^{-1})}(a_3)P_{\calK_{\rho,\vv_r}^R}\zeta_1) \otimes \zeta_2 \\
	&= (\lambda_{(e,e,\uu_r^{-1}\vv_r^{-1})}(a_3)\zeta_1) \otimes \zeta_2\\
	&=\varphi(\lambda_{(e,e,\uu_r^{-1}\vv_r^{-1})}(a_3)\calQ_{(v_{\widetilde{n_r}},\ldots, v_{1})}(\eta_{\widetilde{n_r}}'\otimes \cdots \otimes \eta_{1}'))(U\eta_0)\otimes \zeta_2\\
	&=(U\eta_0)\otimes (\rho_1^{\otimes\widetilde{n_r}}(\lambda^f(a_3))\eta_1\otimes \cdots \eta_{\widetilde{n_r}})\\
	&= (U\otimes \rho_1^{\otimes\widetilde{n_r}}(\lambda^f(a_3)))\eta
\end{align}
This shows equality \eqref{eq:calculation1}.
Hence, combining \eqref{eq:calculation1} and \eqref{eq:calculation2} we obtain
\begin{align}
	\widetilde{\Theta}_d(\lambda(a))_{\rho}&=\theta_1^{\otimes \widetilde{n_l}}(\lambda^f(a_1))\otimes a_2\otimes \rho_1^{\otimes \widetilde{n_r}}(\lambda^f(a_3))\\
	&=(J_{\rho}^{L*}\otimes \Id)(\Id\otimes P_{\calK_{\rho,\vv_l}^L}\lambda_{(\vv_l\uu_l,e,e)}(a_1)Ua_2\otimes \rho_1^{\otimes \widetilde{n_r}}(\lambda^f(a_3)))\\
	&=(J_{\rho}^{L*}\otimes \Id)(\Id\otimes P_{\calK_{\rho}^L}\lambda(a)P_{\calK_{\rho}^R}\otimes \Id)(\Id\otimes J_{\rho}^{R})\\
	&=(J_{\rho}^{L*}\otimes \Id)(\Id\otimes \lambda(a)\otimes \Id)(\Id\otimes J_{\rho}^{R})
\end{align}
This shows the equality holds for all $a\in \boldA_d$, and hence, by density it holds on $\calA_{d}$. This then finishes the proof.	
	\end{proof}
	
\end{theorem}

\begin{theorem}\label{thm:graph-product-cb-maps}
	Fix $d\geq 1$ and for $v\in V\Gamma$,  let $T_{v}:\boldA_{v}\to \boldA_{v}$ be a state-preserving completely bounded map that naturally extends to a bounded map on $L_2(\boldA_{v},\varphi_v)$ and $L_2(\boldA_{v}^{op},\varphi_v)$. Then, for the reduced graph product, the map $T_d:\calA_d\to \calA_d$ defined for $a=a_1\otimes \cdots \otimes a_d\in \mathring{\boldA}_{\vv}\subseteq \mathring{\boldA}_{d}$ as
	\begin{align}
		T_d(\lambda(a_1\otimes \cdots \otimes a_d)) = \lambda(T_{v_1}(a_1)\otimes \cdots \otimes T_{v_d}(a_d))
	\end{align} 
	admits a completely bounded extension on $\calA_d$ with
	\begin{align}
		\|T_d\|_{cb}\leq (\#\Cliq(\Gamma))^3 d \cdot \left(\max_{v} C(T_v)\right)^d.
	\end{align}
	where 
	\begin{align}
		C(T_v) := \max\{\|T_v\|_{cb},\|T_v\|_{\calB(L_2(\boldA_{v},\varphi_v))},\|T_{v}\|_{\calB(L_2(\boldA_{v}^{op},\varphi_v))}\}.
	\end{align}
	 We will denote this map as $T_{d}:=*_{v,\Gamma} T_{v}$. Moreover, if $(S_v)_{v\in V\Gamma}$ are maps satisfying the same conditions as $(T_v)_{v\in V\Gamma}$ then
		\begin{align}\|T_d - S_d\|_{cb}\leq (\#\Cliq(\Gamma))^3 d^2 \left( \max_{v}\max\{C(T_{v}),C(S_{v})\}\right)^{d-1}\max_{v}C(T_{v}-S_v).
		\end{align}
	\begin{proof}
		Fix $d\geq 1$ and suppose first that for all $1\leq i\leq d$ we are given maps $E_{v,i}:\boldA_{v}\to \boldA_{v}$ satisfying the assumptions of the theorem for $T_{v}$.
		Now for $1\leq i\leq d$ the direct sum $\bigoplus_{v\in V\Gamma}E_{v,i}$ extends to a bounded map on $(\oplus_{v\in V\Gamma} \mathring{\calH}_{v})_{C}$. Moreover, by \cite[Theorem 3.4.1]{effrosOperatorSpaces2000} this map is in fact completely bounded  with the same norm. Hence by \eqref{eq:L1K1-complete-isomorphism} the map $E_{L,i} := (\bigoplus_{v\in V\Gamma} E_{v,i})$ is completely bounded on $L_1$ with norm $\|E_{L,i}\|_{cb}\leq \max_{v\in V\Gamma}\|E_{v,i}\|_{\calB(L_2(\boldA_{v},\varphi_v))}$. Similarly we obtain that $E_{R,i}:=(\bigoplus_{v\in V\Gamma} T_{v,i})$ is completely bounded on $K_1$ with norm $\|E_{R,i}\|_{cb}\leq \max_{v\in V\Gamma}\|E_{v,i}\|_{\calB(L_2(\boldA_{v}^{op},\varphi_v))}$. Now, fix a tuple $\rho = (n_l,n_r,\uu_l,\uu_r,\tbold)$ and denote $\widetilde{n}_l = n_l +|\uu_l|$ and $\widetilde{n}_r= n_r+|\uu_r|$. Then by \cite[Proposition 9.2.5]{effrosOperatorSpaces2000} we obtain that
		$$\Pi_{\rho}[(E_{v,i})_{v,i}]:=E_{L,1}\otimes \cdots \otimes E_{L,\widetilde{n_l}}\otimes E_{t_1,\widetilde{n_l}+1}\otimes \cdots \otimes E_{t_{|\tbold|},\widetilde{n_l}+|\tbold|} \otimes E_{R,\widetilde{n_l}+|\tbold|+1}\otimes \cdots \otimes E_{R,d}$$ 
		is a completely bounded map on $L_{\widetilde{n_l}}\otimes_h \boldA_{\tbold}\otimes_h K_{\widetilde{n_r}}$ with norm 
		\begin{align}
			\|\Pi_\rho[(E_{v,i})_{v,i}]\|_{cb} \leq \prod_{i=1}^{\widetilde{n_l}}\|E_{L,i}\|_{cb}\prod_{i=1}^{|\tbold|}\|E_{t_{i},\widetilde{n_l}+i}\|_{cb}\prod_{i=\widetilde{n_l}+|\tbold|+1}^{d}\|E_{R,i}\|_{cb}
			&\leq \prod_{i=1}^d \max_{v}C(E_{v,i}).
		\end{align} Now let the maps $(T_v)$ be given and set $T_{\rho} = \Pi_{\rho}[(T_{v})_{v,i}]$ (i.e. taking $E_{v,i} = T_{v}$ for all $i$). Hence, we get a completely bounded map $\widetilde{T}_d = (T_{\rho})_{\rho}$ on $\widetilde{X}_d$. Denote $T_d'$ for the natural product map on $\boldA_{d}$ that is given by $T_{v_1}\otimes \cdots \otimes T_{v_d}$ on $\mathring{\boldA}_{\vv}$ for $\vv=v_1\cdots v_d$.
		We then find 
		\begin{align}
		T_d\circ \lambda|_{\boldA_d} = \lambda \circ T_d'|_{\boldA_d} =\pi_d\circ j_d\circ T_d'|_{\boldA_d} = \pi_d \circ D_d\circ \widetilde{T}_d\circ \widetilde{\Theta}_d\circ \lambda|_{\boldA_d}.
		\end{align}
		This shows that $T_d$ extends to a completely bounded map on $\calA_d$. The norm-bound now follows from the bound $\|\pi_d\|_{cb}\leq (\#\Cliq(\Gamma))^3 d$, the bound on $\|\widetilde{T}_d\|_{cb}$ and the fact that $D_d$ and $\widetilde{\Theta}_d$ are completely contractive.\\
		
		Now suppose we are given maps $(T_v)_{v\in V}$ and $(S_v)_{v\in V}$ satisfying the assumptions of the theorem. Set $S_{\rho} := \Pi_{\rho}[(S_v)]$ and $\widetilde{S}_{d} := (S_{\rho})_{\rho}$. Set $E_{v,i,j} = T_{v}$ for $i<j$, set $E_{v,i,j} = T_{v} - S_{v}$ for $i=j$ and set $E_{v,i,j} = S_{v}$ for $i>j$. Then by cancellation it follows that $\Pi_{\rho}[(T_v)] - \Pi_{\rho}[(S_v)] = \sum_{j=1}^d \Pi_{\rho}[(E_{v,i,j})_{v,i}]$. Thus it follows that 
		\begin{align}
			\|T_{\rho} - S_{\rho}\|_{cb} \leq \sum_{j=1}^d \|\Pi_{\rho}[(E_{v,i,j})_{v,i}]\|_{cb}\leq
			\sum_{j=1}^d\prod_{i=1}^d \max_{v}C(E_{v,i,j})\\
			\leq
			d\left( \max_{v}\max\{C(T_{v}),C(S_{v})\}\right)^{d-1}\max_{v}C(T_{v}-S_v).
		\end{align}
	
	Then as
	$(T_d - S_d)\circ \lambda|_{\boldA_{d}} = \pi_d \circ D_d\circ (\widetilde{T}_d - \widetilde{S}_d) \circ \widetilde{\Theta}_d \circ \lambda|_{\boldA_{d}}$ we obtain
	$\|T_{d}-S_{d}\|_{cb} \leq \|\pi_d\|_{cb}\max_{\rho}\|T_{\rho}-S_{\rho}\|_{cb}$ which proves the bound.
	\end{proof}
\end{theorem}
Additionally we prove an analogue of \cref{thm:graph-product-cb-maps} for the Hilbert spaces.

\begin{theorem}\label{thm:graph-product-bounded-maps}
	Let $\Gamma$ be a finite graph and for $v\in V\Gamma$ let $(\boldA_{v},\varphi_v)$ and $(\boldB_{v},\psi_{v})$ be unital $\Cstar$-algebras with GNS-faithful states and consider the reduced graph products $\calA$ and $\calB$ respectively.
	Fix $d\geq 1$ and for $v\in V\Gamma$,  let $T_{v}:\boldA_{v}\to \boldB_{v}$ be state-preserving maps that extend to bounded maps from $L_2(\boldA_{v},\varphi_v)$ ($=\calH_{v}^\calA$) to  $L_2(\boldB_{v},\psi_v)$ ($=\calH_{v}^\calB$).  Then the map $T_d:\calF_{d}^{\calA}\to \calF_{d}^{\calB}$ defined for $\eta=\eta_1\otimes \cdots \otimes \eta_d\in \mathring{\calH}_{\vv}^\calA\subseteq \calF_{d}^{\calA}$ as
	\begin{align}
		T_d(\eta) = T_{v_1}(\eta_1)\otimes \cdots \otimes T_{v_d}(\eta_d)
	\end{align} 
	extends to a bounded map. Moreover, if $(S_v)_{v\in V\Gamma}$ are maps satisfying the same conditions as $(T_v)_{v\in V\Gamma}$ then
	\begin{align}\label{eq:difference-estimate-L2}
		\|&T_{d}-S_{d}\|_{\calB(\calF_{d}^{\calA},\calF_{d}^{\calB})} \\
		&\leq d(\max_{v}\max\{\|T_{v}\|_{\calB(\mathring{\calH}_{v}^{\calA},\mathring{\calH}_{v}^{\calB})} , \|S_{v}\|_{\calB(\mathring{\calH}_{v}^{\calA},\mathring{\calH}_{v}^{\calB})}\})^{d-1}\max_{v}\|T_{v}-S_{v}\|_{\calB(\mathring{\calH}_{v}^{\calA},\mathring{\calH}_{v}^{\calB})}
	\end{align}
	\begin{proof}
		Fix $d\geq 1$ and for $v\in V\Gamma$ and $1\leq i\leq d$ let $E_{v,i}:\boldA_{v}\to\boldB_{v}$ be state-preserving that extend to a map in $\calB(\calH_{v}^{\calA},\calH_{v}^{\calB})$. Then as $E_{v,i}$ is state-preserving we have $E_{v,i}(\mathring{\calH}_{v}^{\calA}) \subseteq \mathring{\calH}_{v}^{\calB}$ so that the map $\Pi[(E_{v,i})]:\calF_{d}^{\calA}\to \calF_{d}^{\calB}$ defined for $\eta=\eta_1\otimes \cdots \otimes \eta_d\in \mathring{\calH}_{\vv}^{\calA}\subseteq \calF_{d}^{\calA}$ as
		\begin{align}
			\Pi[(E_{v,i})_{v,i}](\eta) = E_{v_1,1}(\eta_1)\otimes \cdots \otimes E_{v_d,d}(\eta_d)
		\end{align}  is well-defined algebraically and maps $\mathring{\calH}_{\vv}^{\calA}$ to $\mathring{\calH}_{\vv}^{\calB}$ for $\vv\in W$. Hence, since these subspaces are mutually orthogonal for $\vv\in W$ we obtain  \begin{align}
		\|\Pi[(E_{v,i})]\|_{\calB(\calF_{d}^{\calA},\calF_{d}^{\calB})} 
		&= \max_{\vv\in W, |\vv|=d}\|\Pi[(E_{v,i})]\|_{\calB(\mathring{\calH}_{\vv}^{\calA},\mathring{\calH}_{\vv}^{\calB})}\\
		&=\max_{\vv\in W, |\vv|=d} \prod_{i=1}^d \|E_{v_i,i}\|_{\calB(\mathring{\calH}_{v_i}^{\calA},\mathring{\calH}_{v_i}^{\calB})} \\
		&\leq \prod_{i=1}^d \max_{v}\|E_{v,i}\|_{\calB(\mathring{\calH}_v^{\calA},\mathring{\calH}_v^{\calB})}
	\end{align}
	Now let $(T_v)$ and $(S_v)$ be maps satisfying the conditions from the theorem. We see that $T_{d} = \Pi[(T_{v})_{v,i}]$ (i.e. taking $E_{v,i}=T_{v}$ for all $1\leq i\leq d$) and $S_{d} = \Pi[(S_{v})_{v,i}]$ so these maps are indeed bounded.
	Now set $E_{v,i,j} = T_{v}$ for $i<j$, set $E_{v,i,j}=T_{v}-S_{v}$ for $i=j$ and set $E_{v,i,j} = S_{v}$ for $i>j$. It follows from cancellation that 
	\begin{align}
		\Pi[(T_{v})_{v,i}] - \Pi[(S_{v})_{v,i}] = \sum_{j=1}^d \Pi[(E_{v,i,j})_{v,i}]
	\end{align}
	Hence $\|T_{d}-S_{d}\|_{\calB(\calF_{d}^{\calA},\calF_{d}^{\calB})} \leq \sum_{j=1}^d \|\Pi[(E_{v,i,j})_{v,i}]\|_{\calB(\calF_{d}^{\calA},\calF_{d}^{\calB})}$ from which \eqref{eq:difference-estimate-L2} follows.	
	\end{proof}
\end{theorem}
\section{u.c.p extension for CCAP is preserved under graph products}
\label{section:ucp-extension-for-ccap-preserved}
We will introduce the following definition, originating from \cite[Section 4]{ricardKhintchineTypeInequalities2006a}.
Let $(\boldA, \varphi)$ be a unital $\Cstar$-algebra with GNS-faithful state $\varphi$. We will say that it \textit{has a u.c.p extension for the CCAP}, when the following are all satisfied: \begin{enumerate}
	\item There is a net $(V_{j})_{j\in J}$ of finite rank state-preserving maps on $\boldA$ that converge to the identity pointwise and with $\limsup\limits_{j} \|V_{j}\|_{cb} =1$.
	 \item There is a unital $\Cstar$-algebra $(\boldB,\psi)$ that contains $\boldA$ as a unital subalgebra, and s.t. $\psi$ is GNS-faithful and extends the state $\varphi$.
	 \item  There exist a net $(U_j)_{j\in J}$ of state-preserving u.c.p. maps $U_{j}:\boldA\to \boldB$ s.t. 
$\|V_{j}-U_j\|_{cb}$, $\|V_{j}-U_j\|_{\calB(L_2(\boldA,\varphi),L_2(\boldB,\psi))}$ and $\|V_j - U_j\|_{\calB(L_2(\boldA^{op},\varphi),L_2(\boldB^{op},\psi))}$ all converge to $0$ as $j\to \infty$.
\end{enumerate}
Note that by the first property $(\boldA,\varphi)$ must posses the CCAP.  It is clear that any finite dimensional $\Cstar$-algebra possesses the above property. In \cite [proof of Theorem 4.13]{ricardKhintchineTypeInequalities2006a} it was proven that the reduced group $\Cstar$-algebra of any discrete group that possess the CCAP, also satisfies above criteria. In \cite[proof of Theorem 4.2]{freslonNoteWeakAmenability2012} it was proven that the same is true for reduced $\Cstar$-algebra of a compact quantum group with Haar state whose discrete dual quantum group is weakly amenable with Cowling-Haagerup constant $1$.\\

We will now show in the next theorem that the property of having a u.c.p extension for the CCAP is being preserved under graph products, for finite simple graphs. The proof  imitates the proof method of \cite [Proposition 4.11]{ricardKhintchineTypeInequalities2006a}. We will use here \cref{prop:graph-product-of-ucp-maps-vNa}, \cref{prop:graph-products-ucp-maps-Cstar}, 
\cref{prop:semi-group-ucp-maps} and \cref{thm:graph-product-cb-maps} and 
\cref{thm:graph-product-bounded-maps}
\begin{theorem} \label{thm:ucp-extension-gives-CCAP-for-graph-products} Let $\Gamma$ be a finite simple graph and for $v\in V\Gamma$ let $(\boldA_v,\varphi_v)$ be unital $\Cstar$-algebras (with GNS-faithful states $\varphi_v$) that have a u.c.p. extension for the CCAP. Then the reduced graph product $(\calA,\varphi) = *_{\Gamma}(\boldA_{v},\varphi_{v})$ has a u.c.p. extension for the CCAP.\label{thm:ucp-extension-CCAP-preserved-under-graph-products}
	\begin{proof}
			We let $(V_{v,j})_{j\in J_v}$, $(\boldB_{v},\psi_v)$ and $(U_{v,j})_{j\in J_{v}}$ be a u.c.p extension for the CCAP for $(\boldA_{v},\varphi_v)$. As for all $v$ the algebras $\boldA_{v},\boldB_v$ have GNS-faithful states, their reduced graph products $(\calA,\varphi)$ and $(\calB,\psi)$ respectively are well-defined, and the states $\varphi$ and $\psi$ are GNS-faithful as well. Hence, by \cite[Proposition 3.12]{caspersGraphProductsOperator2017a} there exists a unital $*$-homomorphism $\pi:\calA\to \calB$ that intertwines the graph product states. Now for $a\in \ker\pi$ and $b\in \lambda(\boldA)$ we have $\varphi(b^*a^*ab) = \psi(\pi(b^*)\pi(a)^*\pi(a)\pi(b)) = 0$. By the faithfulness of the GNS-representation of $\calA$, this shows that $a=0$ and hence $\pi$ is injective. We will hence consider $\pi$ as an inclusion $\calA\subseteq \calB$.

			We construct a single directed set $\calJ =\prod_{v\in V\Gamma}J_v$ with partial order $(j_v)_{v\in V\Gamma}\prec (j_v')_{v\in V\Gamma}$ if and only if  $j_v\prec j_v'$ for all $v\in V\Gamma$.  We can now define nets $(V_{v,j})_{j\in \calJ}$, $(U_{v,j})_{j\in \calJ}$ as follows: for  $j=(j_v)_{v\in V\Gamma}$ we set $V_{v,j} := V_{v,j_v}$, and $U_{v,j} := U_{v,j_v}$.  Note that these nets still satisfy the assumptions of a u.c.p. extension for CCAP. For $v\in V\Gamma, j\in \calJ$ we will set
				$$\epsilon_{v,j} = \|V_{v,j}-U_{v,j}\|_{cb} +\|V_{v,j}-U_{v,j}\|_{\calB(L_2(\boldA_{v},\varphi),L_2(\boldB_{v},\psi))}  + \|V_{v,j}-U_{v,j}\|_{\calB(L_2(\boldA_v^{op},\varphi),L_2(\boldB_v^{op},\psi))}$$ and by restricting to a subnet we can assume $\epsilon_{v,j}<1$. Since the  maps $U_{v,j}$ are u.c.p and state-preserving we have that $U_{v,j}$ is a contraction from $L_2(\boldA_v,\varphi_v)$ to $L_2(\boldB_v,\psi_v)$ and from $L_2(\boldA_v^{op},\varphi_v)$ to $L_2(\boldB_v^{op},\psi_v)$. Hence we also obtain $$\|V_{v,j}\|_{cb},\|V_{v,j}\|_{\calB(L_2(\boldA_{v},\varphi),L_2(\boldB_{v},\psi))},\|V_{v,j}\|_{\calB(L_2(\boldA_{v}^{op},\varphi_v),L_2(\boldB_{v}^{op},\psi_v))}\leq 2$$

			We can now by \cref{thm:graph-product-cb-maps} construct for $j\in \calJ$, the finite rank c.b. maps $F_{d,j}=*_{v,\Gamma} V_{v,j}$ on $\calA_d$.  We then obtain completely bounded, finite rank maps $D_{N,j} = \sum_{d=0}^{N}(1-\frac{1}{\sqrt{N}})^dF_{d,j}\calP_d$ on $\calA$ that on the dense subset $\lambda(\boldA)$ tend in norm to the identity as $N,j\to \infty$.
			We can now by \cref{prop:graph-products-ucp-maps-Cstar} construct the state-preserving u.c.p maps $U_j:=*_{v,\Gamma}U_{v,j}$, and by \cref{prop:semi-group-ucp-maps} construct the u.c.p maps $\calT_{1-\frac{1}{\sqrt{N}}}$ and the c.b. maps $\calT_{(1-\frac{1}{\sqrt{N}}),N}$ on $\calA$. This gives us state-preserving u.c.p maps  $E_{N,j} = U_j\circ \calT_{1-\frac{1}{\sqrt{N}}}$ and state-preserving c.b. maps $\widetilde{D}_{N,j} = U_j\circ T_{1-\frac{1}{\sqrt{N}},N}$. Applying \cref{thm:graph-product-cb-maps} and using that $C(V_{v,j}),C(U_{v,j})\leq 2$ and $C(V_{v,j}-U_{v,j}) \leq \epsilon_{v,j}$ we obtain	\begin{align}\|F_{d,j} - U_{j}|_{\calA_d}\|_{cb}\leq (\#\Cliq(\Gamma))^3 d^2 2^{d-1}(\max_{v}\epsilon_{v,j})\to 0 \text{ as }j\to\infty.
				\end{align}  
			Similarly, by \cref{thm:graph-product-bounded-maps} we obtain
			\begin{align}\|F_{d,j} - U_{j}\|_{\calB(\calF_{d}^{\calA},\calF_{d}^{\calB})} \leq d 2^{d-1}(\max_{v}\epsilon_{v,j})\to 0 \text{ as }j\to\infty.
			\end{align}  
			Now
			\begin{align}
				\| E_{N,j} - D_{N,j}\|_{cb} &\leq \|E_{N,j} - \widetilde{D}_{N,j}\|_{cb} +  \|\widetilde{D}_{N,j} - D_{N,j}\|_{cb}\\
				&\leq \|\calT_{1-\frac{1}{\sqrt{N}}} - \calT_{1-\frac{1}{\sqrt{N}},N}\|_{cb} +  \sum_{d=0}^N \|U_j|_{\calA_{d}} - F_{d,j}\|_{cb}\|\calP_d\|_{cb}
			\end{align} 
		and similarly 
			\begin{align}
			\| E_{N,j}& - D_{N,j}\|_{\calB(\calF^{\calA},\calF^{\calB})}\\ &\leq \|E_{N,j} - \widetilde{D}_{N,j}\|_{\calB(\calF^{\calA},\calF^{\calB})} +  \|\widetilde{D}_{N,j} - D_{N,j}\|_{\calB(\calF^{\calA},\calF^{\calB})}\\
			&\leq \|U_j\|_{\calB(\calF^{\calA},\calF^{\calB})}\|\calT_{1-\frac{1}{\sqrt{N}}} - \calT_{1-\frac{1}{\sqrt{N}},N}\|_{\calB(\calF^{\calA},\calF^{\calA})} \\
			&+  \sum_{d=0}^N \|U_j|_{\calA_{d}} - F_{d,j}\|_{\calB(\calF_{d}^{\calA},\calF_{d}^{\calB})}\|\calP_d\|_{\calB(\calF^{\calA},\calF_{d}^{\calA})}
		\end{align}
		and we note that $\|\calP_d\|_{cb}\leq C_{\Gamma}d$ (by \cref{thm:projection-maps-linear-growth}), $\|\calP_d\|_{\calB(\calF^{\calA},\calF_{d}^{\calA})} \leq 1$ and $\|U_j\|_{\calB(\calF^{\calA},\calF^{\calB})}=1$. 
		We now obtain using \cref{prop:semi-group-ucp-maps} that
		\begin{align}
		\lim_N\lim_j	\| E_{N,j} - D_{N,j}\|_{cb}&\leq \lim_N  \|\calT_{1-\frac{1}{\sqrt{N}}} - \calT_{1-\frac{1}{\sqrt{N}},N}\|_{cb} \\
		&\leq \lim_{N}C_{\Gamma} N^2 (1-\frac{1}{\sqrt{N}})^N=0
		\end{align} so that in particular $\lim_{N}\lim_{j}\|D_{N,j}\|_{cb}=1$. Similarly we obtain 
		\begin{align}
		\lim_N\lim_j	\| E_{N,j} - D_{N,j}\|_{\calB(L_2(\calA,\varphi),L_2(\calB,\psi))} &\leq \lim_N  \|\calT_{1-\frac{1}{\sqrt{N}}} - \calT_{1-\frac{1}{\sqrt{N}},N}\|_{\calB(\calF^{\calA},\calF^{\calA})} \\
		&\leq \lim_{N}\sup_{d\geq N}(1-\frac{1}{\sqrt{N}})^d =0
	\end{align}
	and analogously $\lim_N\lim_j	\| E_{N,j} - D_{N,j}\|_{\calB(L_2(\calA^{op},\varphi),L_2(\calB^{op},\psi))}=0$ can be shown. Now the construction of $(D_{N,j})$,  $(\calB,\psi)$ and $(E_{N,j})$ shows that $(\calA,\varphi)$ has a u.c.p extension for the CCAP.
	\end{proof}
\end{theorem}

Reasoning similarly to \cite[Corollary 3.17]{caspersGraphProductsOperator2017a} we show for arbitrary (possibly infinite) simple graphs that, under the assumptions on the algebras $\boldA_{v}$, we have that the reduced graph product possesses the CCAP.
\begin{theorem}\label{theorem:CCAP-preserved-under-graph-products}
	Let $\Gamma$ be a simple graph and for $v\in V\Gamma$ let $(\boldA_v,\varphi_v)$ be unital $\Cstar$-algebras that have a u.c.p. extension for the CCAP. Then the reduced graph product $(\calA,\varphi) = *_{\Gamma}(\boldA_{v},\varphi_{v})$ has the CCAP.
\begin{proof}
		It follows from \cref{thm:ucp-extension-CCAP-preserved-under-graph-products} that for any finite subgraph $\Gamma_0\subseteq \Gamma$, the reduced graph product $*_{v,\Gamma_0}\boldA_{v}$ possesses the CCAP. As the reduced graph product over $\Gamma$ is the induced limit of all reduced graph products over finite subgraphs and as the CCAP is preserved under inductive limits, this shows the result
	\end{proof}
\end{theorem}

\begin{corollary}\label{corollary:CCAP-preserved-for-certain-algebras}
	Let $\Gamma$ be a simple graph and for $v\in V\Gamma$ let $\boldA_{v}$ be one of the following:
	\begin{enumerate}
		\item $(\boldA_{v},\varphi_v)$ is a finite-dimensional $\Cstar$-algebra with GNS-faithful state $\varphi_v$.
		\item $(\boldA_v,\varphi_v)$ is the reduced group $\Cstar$-algebra  of a discrete group with Plancherel state $\varphi_v$ that possesses the CCAP
		\item $(\boldA_v,\varphi_v)$ is the reduced $\Cstar$-algebra of a compact quantum group whose discrete dual quantum group is weakly amenable with Cowling-Haagerup constant $1$. Here $\varphi_v$ denotes the Haar state. \label{corollary:quantum-groups}
	\end{enumerate} Then the reduced graph product $\Cstar$-algebra $(\calA,\varphi) = *_{v,\Gamma} (\boldA_v,\varphi_v)$ has the CCAP.
\end{corollary}
We recall, that for a discrete group $G$ we have that $G$ is weakly amenable with constant $1$ if and only if the reduced group $\Cstar$-algebra $C_r^*(G)$ possesses the CCAP, if and only if the group von Neumann algebra $\calL(G)$ possesses the wk-$*$ CCAP. Using this we obtain the following result for von Neumann algebras.
\begin{corollary}\label{corollary:group-von-neumann}
	Let $\Gamma$ be a simple graph and for $v\in V\Gamma$ let $\boldM_v=\calL(G_v)$ be the group von Neumann algebra of a discrete group with the canonical state. If $\boldM_{v}$ has the wk-$*$ CCAP for all $v\in V\Gamma$, then the graph product von Neumann algebra $\calM = \overline{*_{v,\Gamma}}\boldM_v$ possesses the wk-$*$ CCAP as well.
	\begin{proof}
		Note that $\calM = \overline{*_{v,\Gamma}}\calL(G_v) = \calL(*_{v,\Gamma} G_v)$ has the wk-$*$ CCAP if and only if $C_r^*(*_{v,\Gamma}G_v) = *_{v,\Gamma}C_r^*(G_v)$ has the CCAP. The result then follows from \cref{corollary:CCAP-preserved-for-certain-algebras}
	\end{proof}
\end{corollary}
We note that \cref{corollary:group-von-neumann} was already known by \cite{reckwerdtWeakAmenabilityStable2017} where using different techniques it was shown that for discrete groups weak amenability with constant $1$ is preserved under graph products. However, \cref{corollary:CCAP-preserved-for-certain-algebras} does give new examples of algebras that posses the CCAP as you can consider graph products of the form $*_{v,\Gamma}(\boldA_{v},\varphi_v)$ where some of the algebras $(\boldA_v,\varphi_v)$ satisfy condition (1), some satisfy condition (2) and some satisfy condition (3).

\section*{Acknowledgements}
I want to thank Martijn Caspers for his helpful feedback and great supervision during the project. Furthermore, I want to thank the referee for the feedback that substantially helped to improve the exposition of the results. Also I want to thank Jacek Krajczok for comments on a preprint of this paper.

\printbibliography
\end{document}